\documentstyle[12pt]{article}
\begin{titlepage}
\title{Theta and zeta functions for locally symmetric spaces of rank one}

\author{Ulrich Bunke\thanks{Humboldt-Universit\"at zu Berlin, Institut f\"ur
Reine Mathematik (SFB288), Ziegelstr. 13a, Berlin 10099.
E-mail:ubunke@mathematik.hu-berlin.de
       } and
Martin
Olbrich\thanks{Humboldt-Universit\"at zu Berlin, Institut f\"ur Reine
Mathematik (SFB288), Ziegelstr. 13a, Berlin 10099.
E-mail:olbrich@mathematik.hu-berlin.de  }
}
\end{titlepage}

\begin{document}
\newcommand{\oh }{{\bf h}}
\newcommand{\om}{{\bf m}}
\newcommand{\kaaa}{{\bf k}}
\newcommand{\paaa}{{\bf p}}
\newcommand{\taaa}{{\bf t}}
\newcommand{\haaa}{{\bf h}}
\newcommand{\R}{{\bf R}}
\newcommand{\Z}{{\bf Z}}
\newcommand{\C}{{\bf C}}
\newcommand{\HR}{H{\bf R}}
\newcommand{\HC}{H{\bf C}}
\newcommand{\HH}{H{\bf H}}
\newcommand{\PR}{P{\bf R}}
\newcommand{\PC}{P{\bf C}}
\newcommand{\PH}{P{\bf H}}
\newcommand{\Saaa}{{\bf S}}
\newcommand{\G}{{\bf G}}
\newcommand{\A}{{\bf A}}
\newcommand{\Naaa}{{\bf N}}
\newcommand{\g}{{\bf g}}
\newcommand{\aaa}{{\bf a}}
\newcommand{\naaa}{{\bf n}}
\newcommand{\M}{{\bf M}}
\newcommand{\K}{{\bf K}}
\newcommand{\Oaaa}{{\cal O}}
\newcommand{\Haaa}{{\bf H}}
\newcommand{\db}{{\bar{\partial}}}
\newcommand{\Paaa}{{\bf P}}

\maketitle
\newtheorem{prop}{Proposition}[section]
\newtheorem{lem}[prop]{Lemma}
\newtheorem{ddd}[prop]{Definition}
\newtheorem{theorem}[prop]{Theorem}
\newtheorem{kor}[prop]{Corollary}
\newtheorem{ass}[prop]{Assumption}
\newtheorem{con}[prop]{Conjecture}

\tableofcontents
\section{Review of the results}

\subsection{Introduction}\label{intro}

Let $(M,g)$ be a locally symmetric even-dimensional Riemannian manifold with
negative
sectional curvature. Its universal covering $X\rightarrow M$
is a Riemannian symmetric space of rank one. The group of orientation
preserving isometries $\G^\prime$ of
$X$ is a connected semisimple Lie group of real rank one. Then
$X=\G^\prime/\K^\prime$, where
$\K^\prime$ is the maximal compact subgroup of $\G^\prime$. The fundamental
group of $M$ acts by covering transformations on $X$ and gives
rise to a discrete, co-compact subgroup $\Gamma\subset \G^\prime$ such that
$M:=\Gamma\backslash\G^\prime/\K^\prime$.
Let $\G$ be a linear connected finite covering of $\G^\prime$ such that the
embedding
$\Gamma\hookrightarrow\G^\prime$ lifts to an embedding
$\Gamma\hookrightarrow\G$ and let
$\K$ be the maximal compact subgroup of $\G$.
Then $X=\G/\K$ and $M=\Gamma\backslash\G/\K$.
Let $\g=\kaaa\oplus \paaa$ be the Cartan
decomposition  of $\g$ and $\aaa$ be a maximal
abelian subspace in $\paaa$. Let $\M=\G_{\aaa}\cap\K$ be the centralizer
of $\aaa$ in $\K$ and $\hat{\M}$ be the unitary dual.
To $M$ and the $\M$-type $\sigma\in\hat{\M}$ we associate a theta function
$\theta(t,\sigma)$ and a Selberg zeta function $Z_S(p,\sigma)$.
Both will be meromorphic functions on the complex plane.

The $\theta$-function is defined using the spectral data of
elliptic operators over $M$.
The singularities of the theta function
contain geometric
information about the closed geodesics of $M$ and their lifts to the locally
homogeneous bundle $\G\times_\M V_\sigma$ over the unit sphere bundle of $M$.

The Selberg zeta function is,
in some sense, the dual object to the theta function. It is defined,
using the geometric data provided by the closed geodesics of $M$.
This data appears as certain residues of the corresponding theta function.
The singularities of the Selberg zeta function,
in turn, contain the  spectral data employed
in the definition of the theta function.

The theta function and the Selberg zeta function satisfy
functional equations. The functional equation of the theta function
connects the theta function of $M$ with an analogous function
which is associated to the compact dual symmetric space.

The right hand side of the functional equation of the Selberg zeta
function involves Plancherel measures. One of the main ideas of the present
paper is to relate these Plancherel measures to the
harmonic analysis over the compact dual symmetric space.

Let $\G^d\subset \G^c$ be the compact dual group
of $\G$, where $\G^c$ is the complexification of $\G$. Then the
compact dual symmetric space is $X_d:=\G^d/\K$.

The Selberg zeta function $Z_S(p,\sigma)$ associated to the $\M$-type
$\sigma\in \hat{\M}$ is usually defined by an infinite
product converging on some half plan $Re(p)>\rho$.
The problem is to find an analytic continuation and to explain
the singularities (zeros and poles).

The Selberg zeta functions were first introduced by Selberg \cite{selberg56}.
Gangolli \cite{gangolli77} gave a definition in the spherical case, i.e. for
the trivial $\M$-type.
Fried \cite{fried86} defined Selberg zeta functions for general
$\sigma\in\hat{\M}$
in the rank one case.

Gangolli applies the Selberg trace formula in order
to obtain the meromorphic extension of the zeta function, the functional
equation and an explanation of its singularities.

Fried employs the connection of Selberg zeta functions with
the Ruelle zeta function in order to give an alternative proof
of the meromorphic extension based on dynamical systems.
He announced a detailed investigation of the Selberg zeta function in
\cite{fried2}
(unpublished).

The Selberg zeta functions appearing in the work of Wakayama \cite{wakayama85}
are associated to $\K$-types $\gamma\in\hat{\K}$. Using the Selberg
trace formula Wakayama provides a meromorphic extension, a functional equation
and a description of the singularities. Wakayama's zeta function
is closely related to the zeta functions $Z_S(p,\sigma)$, where
$\sigma\subset \gamma_{|\M}$.
The fact that Wakayama's Selberg zeta functions are associated to $\K$-types
causes a rather complicated description of their singularities and
some of them have no explicit  spectral interpretation at all.
We feel that the Selberg zeta function should naturally be
associated to $\M$-types.

Juhl \cite{juhl93} describes a direct approach to the description
of the singularities of the Selberg zeta function in terms
of $\naaa$-cohomologies of  unitary $(\g,\K)$-modules $V_\pi$ and the
multiplicities
$N_\Gamma(\pi)$.  The $\G$-module decomposition
$L^2(\Gamma\backslash\G)=\sum_{\pi\in\hat{\G}} N_\Gamma(\pi)V_\pi$
defines $N_\Gamma(\pi)$. Juhl employs his dynamical Lefschetz formula.

The theta function associated to the trivial $\M$-type was first studied
by Cartier-Voros \cite{cartiervoros90} in the case when $M$ is a Riemann
surface. For general $\sigma\in\hat{\M}$ the theta function $\theta(t,\sigma)$
was introduced by Juhl \cite{juhl931}. In his definition
he employs the data given by a subset of singularities of the
corresponding Selberg zeta function. Juhl's proof of the meromorphic
continuation
is based on the functional equation and the Euler product representation
for the Selberg zeta function
(for related results see Cram\'er \cite{cramer19} and Jorgenson-Lang
\cite{jorgensonlang93}).

The present work is a continuation of our previous paper
\cite{bunkeolbrich93}. It provides a self contained approach to
the theta functions, Selberg and Ruelle zeta functions.

We reprove and strengthen results about the Selberg zeta function which
are already known or at least folklore. Our motivation for writing this paper
is to introduce a new point of view leading to a
better understanding of the zeta functions.

The point of view in the previous papers dealing with Selberg zeta functions
was
\begin{itemize}
\item harmonic analysis on $\G$ (Gangolli \cite{gangolli77}, Wakayama
\cite{wakayama85})
\item the geodesic flow as a dynamical system (Fried \cite{fried86}, Juhl
\cite{juhl93}).
\end{itemize}
In our approach we try to relate everything to the analysis of elliptic
operators on locally
homogeneous bundles over $M$. The main tool is the fundamental solution of
the wave equation which we employ in order to prove a distributional
trace formula. The coefficients of  Hadamard's asymptotic expansion of the
local trace of this fundamental solution are related to the corresponding
coefficients of the asymptotic expansion of the local trace of the fundamental
solution
of the wave equation on the compact dual space. Taking the spectral analysis
on the compact dual space as input we derive the trace formula without
using any harmonic analysis over the non-compact symmetric space.
Once having established a suitable trace formula, our techniques to derive the
properties of the theta functions are the same as in \cite{bunkeolbrich93}.

As in the standard approach we relate the logarithmic derivative of Selberg
zeta function
to the product of resolvents. Instead of regularizing by considering a suitable
derivative of this logarithmic derivative we apply our $\Psi$-construction
(\cite{bunkeolbrich93}, compare with Parnovskij \cite{parnovskij92}). Thus, we
don't have problems with the integration constants.
Since we can relate the singularities of the logarithmic derivative of
the Selberg zeta function to eigenspaces of operators on $M$ and $X_d$
we can easily conclude that all its residues are integers.
This resolves the problem that previously in the bundle case one
only was able to prove (using harmonic analysis)
that the residues are rational.
This would lead the meromorphicity of a certain power of the Selberg zeta
function.
Of course, in the work of Fried the meromorphicity of the Selberg zeta function
was established. But Fried's approach does not provide a functional equation
and a nice interpretation of the  singularities.

An alternative description of the residues of the logarithmic derivative
of the Selberg zeta function in terms of data associated
to the geodesic flow was given by Juhl \cite{juhl93}.

The Ruelle zeta function can be expressed in terms of  Selberg
zeta functions as it was observed by Fried \cite{fried86}.
Thus, in principle, the singularities of the Ruelle zeta function
can be expressed in terms of the singularities of the corresponding
Selberg zeta functions and eventually in terms of eigenvalues
of operators on bundles over $M$ and $X_d$. In order to carry out this program
one has to understand quite explicitly how the representations in $\hat{\K}$,
which are connected
to differential forms, split when restricted to $\M$. A final answer we only
obtained
in the case of real and complex hyperbolic manifolds $M$.

Nevertheless, using a trick, we derive a functional equation for the
Ruelle zeta function. This functional equation was previously obtained by Juhl
up to an exponential polynomial, which turns out be trivial.
The functional equation implies the result of Juhl that the order of the
singularity of the Ruelle zeta function at $s=0$ is $1/2 dim(M)\chi(M)$
(compare also Seifarth \cite{seifarthdiss93} ).

Recall that we want to associate the theta and the zeta functions to
irreducible
representations $\sigma\in\hat{\M}$.
Moreover we want to find elliptic operators over $M$ which provide
the related spectral information.

There is no immediate way to
associate a bundle to $\M$-types $\sigma$ carrying the interesting elliptic
operators.

A bundle $V(\gamma)=\G\times_\K V_\gamma$ can be associated
to a $\K$-type $\gamma\in\hat{\K}$. Let $R(\K)$  be the representation
ring of $\K$ with integer coefficients. Then, more general, a $\Z_2$-graded
bundle $V(\gamma)$ can be associated to $\gamma\in R(\K)$ in the obvious way.

Let $r:R(\K)\rightarrow R(\M)$ be the restriction.
We show that $r$ is surjective (note that we require $dim(M)$ even).
Fix $\sigma\in\hat{\M}$.
We introduce the subset of $r^{-1}(\sigma)$ of $\sigma$-admissible
$\gamma\in R(\K)$.
Next we explain this notion of $\sigma$-admissibility.

To $\gamma\in R(\K)$ we associate
a bundle $V_d(\gamma):=\G^d\times_\K V_\gamma$ over $X_d$.
We introduce a constant $c(\sigma)$ in Definition \ref{d.constants} and
the operator $A_d(\gamma,\sigma):=\sqrt{\Omega_d(\gamma)+c(\sigma)}$,
where $\Omega_d(\gamma)$ is the Casimir of $\G^d$ acting on sections of
$V_d(\gamma)$.
A $\sigma$-admissible $\gamma$ will be characterized by the property
that the spectrum of $A_d(\gamma,\sigma)$ forms the positive part of a certain
lattice $L(\sigma)\subset \R$,
$L(\sigma)=T(\Z+\epsilon(\sigma))$, where $T\in\R$ is the period only
depending on $X_d$
and $\epsilon(\sigma)\in\{0,1/2\}$.
Moreover, the multiplicity of the eigenvalue $\lambda$, denoted by
$m_d(\lambda,\gamma,\sigma)$, is given by a polynomial $P(\lambda,\sigma)$.
Note that we understand the spectrum in the weighted sense using the
$\Z_2$-grading
of $V_d(\gamma,\sigma)$.
Thus the multiplicities may be negative and  $A_d(\gamma,\sigma)$
may have other eigenvalues not contained in the lattice,
but with zero multiplicity.

The fact that $\gamma$ is $\sigma$-admissible, fixes the weighted
multiplicities of the spectrum
of $A_d(\gamma,\sigma)$ completely.
Hence the local (super) trace $K_d(t,\gamma,\sigma)$ of the distributional
kernel
of $cos(tA_d(\gamma,\sigma))$ is determined, too.

On the bundle $V(\gamma)\rightarrow X$ we define
$A(\gamma,\sigma):=\sqrt{\Omega(\gamma)-c(\sigma)}$,
where $\Omega(\gamma)$ is the Casimir of $\G$ acting on sections of
$V(\gamma)$.
Let $K(t,\gamma,\sigma)$ be the local trace of the distributional kernel of
$cos(tA(\gamma,\sigma))$.
We establish that $K(t,\gamma,\sigma)=(-1)^{n/2}K_d(\imath t,\gamma,\sigma)$.
This observation furnishes the bridge between the analysis on $X,M$ and on
$X_d$.

The meromorphic function $K(t,\gamma,\sigma)$ can be expressed in terms of
Fourier transforms of the Plancherel
measures associated to $X$ and $\sigma^\prime$ occurring in $\gamma_{|\M}$.
An inspection of this relation shows that the $\sigma$-admissibility of
$\gamma\in r^{-1}(\sigma)$
can be characterized by the fact, that the contributions of all discrete series
representations of $\G$ in $L^2(X,V(\gamma))$ to $K(t,\gamma,\sigma)$ must
cancel out.
This provides sometimes an easy way to check the $\sigma$-admissibility, namely
if $L^2(X,V(\gamma))$  does not contain any discrete series representation at
all.

In order to derive the properties of the theta and zeta functions associated
to $\sigma\in\hat{\M}$ we use the analysis of the operators $A(\gamma,\sigma)$,
$A_M(\gamma,\sigma)$
and $A_d(\gamma,\sigma)$, where $A_M(\gamma,\sigma)$ is the square root
of $\Omega_M(\gamma)-c(\sigma)$ taken positive or with positive imaginary part.
$\Omega_M(\gamma)$ is the operator on $V_M(\gamma):=\Gamma\backslash V(\gamma)$
induced from $\Omega(\gamma)$. The order of the singularities of the Selberg
zeta functions are expressed in terms of weighted
multiplicities of the eigenvalues of these operators.
By comparison with the results of Juhl \cite{juhl93} it follows that these
weighted multiplicities
can be written in terms of Euler characteristics of $\naaa$-cohomologies.
It would be interesting to see this connection directly.

\subsection{The results}

In this subsection we give the definition of
the theta functions and the zeta functions.
We collect the structural results proved in the course of the
paper together in order to provide a self contained reference.

Let $M$ be an even-dimensional compact rank-one locally symmetric space of
non-compact type. Its universal covering $X$ is a Riemannian symmetric space
of negative sectional curvature with the  group $\G^\prime$ of orientation
preserving isometries.
Let $\Gamma\subset \G^\prime$ be isomorphic to the fundamental group of $M$
such that $M=\Gamma\backslash X$.
Choose a linear connected finite covering $\G$ of $\G^\prime$ such that the
embedding
$\Gamma\hookrightarrow\G^\prime$ lifts to an embedding
$\Gamma\hookrightarrow\G$.
Let $\K$ be the maximal compact subgroup of $\G$ and $\g=\kaaa\oplus\paaa$ be
the
Cartan decomposition of the Lie algebra $\g$ of $\G$.
Fix a one-dimensional  subspace $\aaa\subset \paaa$.
Let $\M:=\K\cap\G_\aaa$ be the centralizer of $\aaa$ in $\K$.
Consider an irreducible representation ($\M$-type) $\sigma\in\hat{\M}$ of $\M$.
Let $R(\K),R(\M)$ be the integer representation rings of $\K,\M$
Choose $\gamma\in R(\K)$ such that $r(\gamma)=\sigma$ and $\gamma$ is
$\sigma$-admissible , $r:R(\K)\rightarrow R(\M)$ being the
restriction. The notion of $\sigma$-admissibility was extensively discussed at
the end of Subsection
\ref{intro}.

We normalize the invariant scalar product on $\g$ such that
it induces a metric on $\aaa$ which has the following property.
Fix $H\in\aaa$ with $|H|=1$. Then $\R\ni t\rightarrow g\:exp(tH)\K\in X$
is a unit speed geodesic.
The scalar product on $\g$ induces an invariant one on $\om$.
The Lie algebra $\g$ splits as $\naaa^-\oplus(\om\oplus\aaa)\oplus\naaa^+$
with respect the negative, zero and positive eigenvalues of $H$.
Let $\rho=\frac{1}{2}tr(H)_{|\naaa^+}\in\R$.
Later we will omit the "$+$"-sign at $\naaa$. Let
$\sigma(\Omega_\M)\in \R$ be the
action on $V_\sigma$ of the Casimir operator of $\M$ defined by the invariant
scalar product on $\om$.
Define $c(\sigma):=\rho^2-\sigma(\Omega_\M)$.

Let $V(\gamma)\rightarrow X$ be the $\Z_2$-graded bundle $\G\times_\K V_\gamma$
and
$V_M(\gamma)$ be the push down to $M$. The Casimir operator $\Omega(\gamma)$
acts as an unbounded selfadjoint operator on $L^2(X,V(\gamma))$.
Define $A(\gamma,\sigma)=\sqrt{\Omega(\gamma)-c(\sigma)}$, where
we take the positive root or the root with positive imaginary part,
respectively.
Let $A_M^2(\gamma,\sigma)$ be the push down to $M$ of
$A(\gamma,\sigma)^2$  and define the square root $A_M(\gamma,\sigma)$ as above.
Similarly we define $A_d(\gamma,\sigma):=\sqrt{\Omega_d(\gamma)+c(\sigma)}$ on
$V_d(\gamma):=\G^d\times_\K V_\gamma\rightarrow X_d$,
where $G^d$ is the compact dual group of $\G$ and $X_d$ is
the compact dual space of $X$.
$X_d$ can be embedded into a complexification of $X$.
The metric on $X_d$ is
obtained from the metric on $X$ by analytic continuation to the imaginary
space.
In particular, this fixes the scale of the metric and thus the volume of $X_d$.

Let $m_d(\lambda,\gamma,\sigma)=:P(\lambda,\sigma)$ be
the multiplicity of the eigenvalue $\lambda>0$ of $A_d(\gamma,\sigma)$, i.e.
$$m_d(\lambda,\gamma,\sigma)=Tr E_{A_d(\gamma,\sigma)}\{\lambda\},$$
where the trace is the supertrace taking the $\Z_2$-grading into account.
$P(\lambda,\sigma)$ is a polynomial in $\lambda$ and does not depend on
$\gamma$
(by the very definition of the $\sigma$-admissibility of $\gamma$).
The eigenvalues of $A_d(\gamma,\sigma)$ with non-trivial multiplicity are
located
in a lattice $L(\sigma)=T(\Z+\epsilon(\sigma))$,
where $T$ is chosen maximal (and is independent of $\sigma$) and
$\epsilon(\sigma)\in\{0,1/2\}$.
Define $Q(\lambda,\sigma)$ by $P(\lambda,\sigma)=\lambda Q(\lambda,\sigma)$.
It turns out that $Q(\lambda,\sigma)$ is a polynomial.

We introduce the distributions
\begin{eqnarray*}
K(t,\sigma)&:=&-(-1)^{n/2}\frac{1}{Tvol(X_d)}\frac{d^2}{dt^2} Q(\frac{d }{dt
},\sigma) ln(|sinh(Tt/2)|),\quad \epsilon(\sigma)=0\\
K(t,\sigma)&:=&-(-1)^{n/2}\frac{1}{Tvol(X_d)}\frac{d^2}{dt^2 } Q(\frac{d }{dt
},\sigma) ln(|tgh(Tt/4)|), \quad \epsilon(\sigma)=1/2.
\end{eqnarray*}

For $g\in\Gamma\not=1$ let $[g]^\G$ be its conjugacy class in $\G$.
There is an element $ma\in[g]^\G$ such that $m\in\M$ and $a=exp(l(g)H)$
with $\R\ni l(g)>0$. Define
$$C(g,\sigma):=-\frac{l(g)e^{\rho l(g)}tr\sigma(m)}{2det(1-Ad(ma)_\naaa)}.$$
Let $n_\Gamma(g)$ be the number of classes in $\Gamma_g/<g>$, where
$\Gamma_g$ is the centralizer of $g$ in $\Gamma$ and $<g>$ is the group
generated by $g$.
By $C\Gamma$ we denote the set of conjugacy classes of $\Gamma$.

For $\phi\in C_c^\infty(\R)$ we define the trace class operator
$$K_{M,\phi}=\int_\infty^\infty \phi(t)cos(tA_M(\gamma,\sigma)).$$
\begin{prop}[The distributional trace formula]
The distribution $$C_c^\infty(\R)\ni \phi\rightarrow Tr K_{M,\phi}$$
is given by
$$vol(M) K(t,\sigma)+\sum_{[g]\in C \Gamma,[g]\not=[1]}
\frac{C(g,\sigma)}{n_\Gamma(g)}(\delta(t-l(g))+\delta(t+l(g))).$$
\end{prop}
Note that we take the trace in the graded sense.

Let $m(\lambda,\gamma,\sigma)$ be the (weighted) multiplicity of the eigenvalue
$\lambda$
of $A_M(\gamma,\sigma)$. The following Proposition is a  consequence
of the $A$-index theorem \ref{aaind}.
\begin{prop}[The well-definedness of $m(\lambda,\sigma)$]
The multiplicity $m(\lambda,\gamma,\sigma)$
is independent of the choice of the $\sigma$-admissible $\gamma$.
\end{prop}
We will write $m(\lambda,\sigma):=m(\lambda,\gamma,\sigma)$.
The multiplicities have polynomial growth and the theta function
$$\theta(t,\sigma):=\sum_{\lambda}m(\lambda,\sigma)e^{-t\lambda}$$
is well defined as a holomorphic function for $Re(t)>0$.
Defining the dual theta function by
$$\theta_d(t,\sigma):=\sum_{0\le \lambda \in L(\sigma)}
m_d(\lambda,\sigma)e^{-t\lambda}$$
we have
in the case $\epsilon(\sigma)=0$
$$\theta_d(t,\sigma)=-\frac{1}{T}\frac{d^2}{dt^2} Q(\frac{d }{dt },\sigma)
ln(|sinh(Tt/2)|)$$
and in the case $\epsilon(\sigma)=1/2$
$$\theta_d(t,\sigma)=-\frac{1}{T}\frac{d^2}{dt^2} Q(\frac{d }{dt },\sigma)
ln(|tgh(Tt/4)|).$$
\begin{theorem}[The theta function]
The theta function $\theta(t,\sigma)$ admits a meromorphic extension to the
whole
complex plane. It satisfies the functional equation
$$\theta(t,\sigma)+\theta(-t,\sigma)=2 vol(M)K(\imath
t,\sigma)=2\frac{\chi(M)}{\chi(X_d)}\theta_d(\imath t,\sigma),$$
where we view $K(\imath t,\sigma)$ as a meromorphic function.
The singularities of $\theta(t,\sigma)$ are
\begin{itemize}
\item first order poles at
$t=\pm \imath l(g)$, $g\in C\Gamma$ with residue
$$res_{t=\pm \imath l(g)}\theta(t)=\frac{C(g,\sigma)}{\pi n_\Gamma(g)},$$
\item poles of order $n$ at $t=2\pi k/T$, $k=-1,-2,\dots$
(at these points the singular part coincides with that of $2 vol(M)K(\imath
t,\sigma)$ )
and
\item a pole of order $n$ at $t=0$.
\end{itemize}
\end{theorem}
Note that $2\pi k/T$, $k=-1,-2,\dots$ are the lengths of the closed geodesics
of $X_d$.

We consider now the Selberg  zeta function.
It is defined for $s>\rho$ by the infinite product
$$Z_S(s,\sigma)=\prod_{[g]\in
C\Gamma,[g]\not=1,n_\Gamma(g)=1}\prod_{k=0}^\infty
\:det\left(1-e^{(-\rho-s)l(g)}S^k(Ad(ma)_\naaa^{-1})\otimes \sigma(m)\right).$$
\begin{theorem}[The Selberg zeta function]
The Selberg zeta function has a meromorphic continuation to the complex plane.
The singularities (zeros have positive and poles negative order) of
$Z_S(s,\sigma)$ are
\begin{itemize}
\item at
$p=\pm \imath \lambda$
of  order
$m(\lambda,\sigma)$,
where $\lambda\not=0$ is an eigenvalue of $A_M(\gamma,\sigma)$,
\item
at $p=0$ of order $2m(0,\sigma)$ if $0$ is an eigenvalue of
$A_M(\gamma,\sigma)$,
\item at
$p=-\lambda$, $\lambda\in T(\Naaa+\epsilon(\sigma))$  of order
$-2\frac{\chi(M)}{\chi(X_d)} m_d(\lambda,\sigma)$.
Then $\lambda>0$ is an eigenvalue of $A_d(\gamma,\sigma)$.
\end{itemize}
If two such points coincide, then the orders add up.
\end{theorem}
\begin{theorem}[Functional equation]
The Selberg zeta function satisfies the functional equation
$$\frac{Z_S(s,\sigma)}{Z_S(-s,\sigma)}= exp\left(
-\frac{\chi(M)}{\chi(X_d)}\frac{2\pi}{T}
\int_{0}^spQ(p,\sigma)\left\{\begin{array}{cc} tg(\pi p/T)
& \epsilon(\sigma)=1/2\\
                                                           - ctg(\pi p/T)&
\epsilon(\sigma)=0 \end{array}\right\}dp\right)
                                                           . $$
\end{theorem}
Note that the integrand is the Plancherel density (up to a constant factor) of
the contribution to $L^2(\G)$ of the principal series associated to $\sigma$.
\begin{prop}[Representation by regularized determinants]
The Selberg zeta function has the representation
\begin{eqnarray}Z_S(p,\sigma)
&=&C
det(p^2+A_M(\gamma,\sigma)^2)
(det(A_d(\gamma,\sigma)^2-p^2))^{-\chi(M)/\chi(X_d)}  \\
&.&exp\left(-\frac{\chi(M)}{\chi(X_d)}\frac{\pi}{T} \int_{p_0}^p
sQ(s,\sigma)
\left\{\begin{array}{cc} tg(\pi s/T)& \epsilon(\sigma)=1/2\\
-ctg(\pi s/T)& \epsilon(\sigma)=0 \end{array}\right\} ds
\right).\nonumber\end{eqnarray}
The constant $C$ is determined by the condition that $Z_S(p,\sigma)$ should
tend to one if $p$ tends to infinity on the positive real axis.
 \end{prop}
The determinants of the elliptic differential operators
are understood in the zeta-regularized sense and take the $Z_2$-grading into
account.

The Ruelle zeta function associated to $M$ is defined for $Re(s)>2\rho$ by
the infinite product
$$Z_R(s):=\prod_{[g]\in C\Gamma,[g]\not=1,n_\Gamma(g)=1} (1-e^{-sl(g)})^{-1}.$$
Using the representation of the Ruelle zeta function in terms of Selberg zeta
functions
we obtain a meromorphic continuation of $Z_R(s)$. In particular we have
\begin{theorem}[Functional equation]
The Ruelle zeta function satisfies the functional equation
$$Z_R(s)Z_R(-s)=  sin(\pi s/T)^{n\chi(M)}.$$
In particular, the order of the singularity at
$s=0$ is $\frac{n}{2}\chi(M)$ while the order of the singularities
at $s<-2\rho$, $s\in T\Z$  is $n\chi(M)$.
\end{theorem}
 We give a satisfactory spectral interpretation
of the singularities of the Ruelle zeta function in the case when $M$ is real
or complex hyperbolic.
In principle, the other cases can be solved in a similar way
but are still open due to computational difficulties.

In the Section \ref{exampleee} we give some more examples for the spectral
interpretation of the singularities of zeta functions.

We tried to state our final results in a self-contained form not involving
any deep machinery but only geometric and spectral data.

A part of the results of the present paper is "weakly"
contained in the previous literature.
Sometimes it was rather complicated to extract them in a clear way.
In many cases we were able to improve the published results.

We hope that our paper not only serves this purpose, but also
gives a new approach to the topic which is accessible to
readers with not much background in harmonic analysis and representation
theory.
Of course, some representation theoretic input is necessary.
Thus, certain facts about the harmonic analysis on compact symmetric
spaces are added as a separate section.
In particular, the generalization of the Cartan-Helgason theorem
giving the spectral decomposition of $A_d(\gamma,\sigma)$
seems not to be contained in the published literature.

In forthcoming paper we will discuss the odd-dimensional
real hyperbolic case.
We also plan to study the theta function in the higher rank case.

{\it We want to thank A. Juhl for his constant interest in this work
and for many discussions. His results and philosophy were helpful
for us to find our way through the theory developed in the present paper
and to check our final results.}
\section{Ideal ladders}

\subsection{Admissible lifts}

Let $R(\K)$, $R(\M)$ denote the representation rings with integer
coefficients of the groups $\K$, $\M$. The inclusion
$\M \subset \K$ induces a restriction map $r:R(\K)\rightarrow R(\M)$,
which is surjective by Proposition \ref{p.sur}.
For any $\gamma\in R(\K)$ let $V_d(\gamma)\rightarrow X_d$ be the associated
$\Z_2$-graded homogeneous vector bundle. More explicitly, let
$\gamma=\sum_{i=1}^l k_i \delta_i$ with $k_i\in\Z$ and $\delta_i\in\hat{\K}$.
Form
\begin{equation}\label{wdef}W^\pm:=\sum_{\{i=1\: :\:sign(k_i)=\pm
1\}}^l\sum_{m=1}^{|k_i|} W_{\delta_m},\end{equation}
where $W_\delta$ is the representation space of $\delta$.
Then $V_d(\gamma)^\pm$ is given by $V_d(\gamma)^\pm=\G_d\times_\K W^\pm$.

For any $\sigma\in\hat{\M}$ we define the shift constant
$c(\sigma)=|\rho|^2+|\rho_m|^2-|\Lambda_\sigma+\rho_m|^2$
(see Definition \ref{d.constants}  for the notation).
Let $\Omega_d(\gamma)$ be the Casimir operator of $\G^d$ acting on sections of
the bundle
$V_d(\gamma)$. We define $c(\sigma)$ and the Casimir operator using an
invariant
scalar product on $\g^d$. The normalization of this scalar product will be
fixed later.
$\Omega_d(\gamma)$ is a positive, selfadjoint operator on
$L^2(X_d,V_d(\gamma))$ with a discrete spectrum. We consider the
positive square root $A_d(\gamma,\sigma):=\sqrt{\Omega_d(\gamma)+c(\sigma)}$.
Let $E_{A_d(\gamma,\sigma)}(.)$ denote the family of spectral projections
of $A_d(\gamma,\sigma)$ and define the multiplicity
$$m_d(\lambda,\gamma,\sigma):=Tr \: E_{A_d(\gamma,\sigma)}({\lambda}),$$
where the trace is the super trace with respect to the $\Z_2$-grading. Thus,
$m_d(\lambda,\gamma,\sigma)$ is a weighted dimension of the eigenspace
of $A_d(\gamma,\sigma)$ corresponding to the eigenvalue $\lambda$.

For the definition of the polynomial $P(\lambda,\sigma)$   and
the lattice $L(\sigma)$ we refer to Definition  \ref{d.constants}.
{}From Proposition \ref{p.CH} immediately follows
\begin{kor} Let $\sigma\in\hat{\M}$ and $\gamma\in R(\K)$ such that
$r(\gamma)=\sigma$.
Then $m_d(\lambda,\gamma,\sigma)=P(\lambda,\sigma)$ for all but finitely many
$0\le \lambda\in L(\sigma)$.
\end{kor}
The exceptional set where $m(\lambda,\gamma,\sigma)\not=P(\lambda,\sigma)$
depends on $\gamma$.
\begin{ddd}
Let $\sigma\in \hat{\M}$. Then $\gamma\in R(K)$ is called $\sigma$-admissible
if $r(\gamma)=\sigma$ and $m_d(\lambda,\gamma,\sigma)=P(\lambda,\sigma)$
for all $0\le \lambda\in L(\sigma)$.
\end{ddd}
Thus, for $\sigma$-admissible $\gamma$ the weighted multiplicities of
spectrum of $A_d(\gamma,\sigma)$ have a particular nice regularity.
In order to construct $\sigma$-admissible elements $\gamma$ we start
with some $\gamma\in R(\K)$  such that $r(\gamma)=\sigma$.
Then we modify $\gamma$ by adding elements in the kernel of $r$.
We control the effect of that on the multiplicities using the $\G^d$-
index theorem.
Since any representation of $\G^d$ can be realized by the index of
homogeneous Dirac operators by Theorem \ref{th.ind} we can modify
$\gamma$ such $m_d(\lambda,\gamma,\sigma)$ changes its value
in a given finite number of places without changing the multiplicity
in the remaining lattice points.
\begin{kor}
Let $\sigma\in \hat{\M}$. Then there exist $\sigma$-admissible $\gamma\in
R(\K)$.
\end{kor}
In fact, there are infinitely many $\sigma$-admissible $\gamma\in R(\K)$.

\subsection{The dual theta function}

In the following we will fix a $\sigma\in \hat{\M}$ and a
$\sigma$-admissible $\gamma\in R(\K)$.
We simplify our notation by omitting  $\sigma$ and $\gamma$. Thus,
$m_d(\lambda):=m_d(\lambda,\gamma,\sigma)$
is independent of the choice of the $\sigma$-admissible $\gamma$.
\begin{ddd}
The dual theta function is defined by
$$\theta(t\sigma)=\theta(t):= Tr e^{-tA_d}=\sum_{\lambda} m_d(\lambda)
e^{-t\lambda} ,\quad Re(t)>0.$$
\end{ddd}
By construction we have
$$\theta_d(t)=\sum_{0\le \lambda\in L} P(\lambda) e^{-t\lambda}.$$
{}From this it can be seen using the Borel-Weil-Bott theorem that our
definition of $\theta_d$
coincides with that of Juhl \cite{juhl931}.

$\theta_d(t)$ is a holomorphic function on the right half plane. We will
provide
a meromorphic extension of $\theta_d$ to the whole complex plane.
Recall the notation $L:=L(\sigma)$, $P(\lambda):=P(\lambda,\gamma,\sigma)$,
the period $T$
of $L$ and $\epsilon=\epsilon(\sigma)\in \{0,1/2\}$ (Definition
\ref{d.constants}).
Since $P(-\lambda)=-P(\lambda)$ we can write $P(\lambda)=\lambda Q(\lambda )$
for some polynomial $Q$.
\begin{lem}\label{ll1}
In the case $\epsilon=0$ we have
$$\theta_d(t)=\frac{T}{4}Q(\frac{d }{dt }) \frac{1}{sinh^2(tT/2)}.$$
In the case $\epsilon=1/2$ we have
$$\theta_d(t)=\frac{T}{4} Q(\frac{d }{dt }) \frac{cosh(tT/2)}{sinh^2(tT/2)}.$$
\end{lem}
{\it Proof:}
Consider the case $\epsilon=0$.
In the half plane $Re(t)>0$ we have
\begin{eqnarray*}
\theta_d(t)&=&\sum_{0\le \lambda\in L} P(\lambda) e^{-t\lambda}\\
&=&\sum_{k=0}^\infty -\frac{d}{dt} Q(\frac{d }{dt })e^{-tTk}\\
&=&-Q(\frac{ d}{dt })\frac{d}{dt} \frac{1}{1-e^{-tT}}\\
&=&\frac{T}{4}Q(\frac{d }{dt })\frac{1}{sinh^2(tT/2)}.
\end{eqnarray*}
In the case $\epsilon=1/2$ for $Re(t)>0$ we have
\begin{eqnarray*}
\theta_d(t)&=&\sum_{0\le \lambda\in L} P(\lambda) e^{-t\lambda}\\
&=&\sum_{k=0}^\infty -\frac{d}{dt} Q(\frac{d }{dt })e^{-tT(k+1/2)}\\
&=&-Q(\frac{d }{dt })\frac{d}{dt} \frac{e^{-tT/2}}{1-e^{-tT}}\\
&=&\frac{T}{4}Q(\frac{d }{dt })\frac{cosh(tT/2)}{sinh^2(tT/2)}.
\end{eqnarray*}
$\Box$\newline
\begin{kor} The dual theta function $\theta_d(t)$ has a meromorphic
continuation to the whole complex plane with singularities on the
imaginary axis at the points $\frac{2\pi\imath}{T}\Z$.
\end{kor}
\section{The distributional trace formula}

\subsection{The distributions $K(t)$ and $K_d(t)$}

Let $X:=\G/\K$ be a rank-one symmetric space of non-compact type.
We normalize the metric of $X$ in the compatible way with
the invariant scalar product on $\g$.
Let $\Gamma\subset \G$ be a co-compact, torsion-free
discrete subgroup and $M:=\Gamma \backslash X$ be the
compact, locally symmetric space.
We fix some $\sigma \in\hat{\M}$ and choose some $\sigma$-admissible
$\gamma\in R(\K)$. We will again omit $\sigma$, $\gamma$ in our notation.
Recall the definition (\ref{wdef}) of $W^\pm$.
Let $V:=V(\gamma)=\G\times_{\K} (W^+\oplus W^-)$
be the associated $\Z_2$-graded homogeneous vector bundle. It induces a
corresponding bundle $V_M$ over $M$.
Let $-\Omega$ be the action of the Casimir operator of $\G$ on $V$ and
$\Omega_M$
be push down of $\Omega$ to $V_M$.
They give rise to selfadjoint operators on the Hilbert spaces $L^2(X,V)$ and
$L^2(M,V_M)$, respectively. We form the square roots
\begin{ddd}\label{hieradef}
$$A:=\sqrt{\Omega-c(\sigma)}, \quad A_M:=\sqrt{\Omega_M-c(\sigma)},$$
\end{ddd}
where take the positive root on the positive spectrum and the root
with positive imaginary part on the negative spectrum.

For $\phi\in C_c^\infty(\R)$ let $K_\phi$ be the operator
$$K_\phi:=\int_R \phi(t) cos(tA) dt.$$
This operator has a smooth integral kernel.
By the $\G$-invariance of $A$ the local (super) trace of this
kernel is a constant. Choose some $x\in X$. We define the distribution
$K\in C_c^\infty(\R)^\ast$ by
$$<K,\phi>:=tr K_\phi(x,x),$$
where the trace is taken on $End(V_x)$.
The definition of $K$ does not depend on $x\in X$.
A similar distribution $K_d$ can be defined on the compact dual side.
Choose some $x_d\in X_d$. Let $K_{d,\phi}$ be the operator
$$K_{d,\phi}:=\int_R \phi(t) cos(tA_d) dt.$$
This operator again has a smooth integral kernel and we define
$$<K_d,\phi>:=tr K_{d,\phi}(x_d,x_d),$$ where the trace is taken on
$End(V_{d,x_d})$.
Fix some $a\in (0,2\pi/T)$. On $(-a,a)\setminus \{0\}$ the distribution
$K_d(t)$ is regular and given by
$\frac{1}{2vol(X^d)}(\theta_d(\imath t)+\theta_d(-\imath
t))=\frac{1}{vol(X^d)}\theta_d(\imath t)$ because of the symmetry of
$\theta_d$. In particular, $K_d(t)$ has a convergent asymptotic expansion
$$K_d(t)\stackrel{t\to 0}{\sim} \sum_{k=-n}^\infty b_{d,k} t^k.$$
The numbers $b_{d,k}$ can be obtained from the explicit expression
Lemma \ref{ll1}. Alternatively, the coefficients of the asymptotic expansion
can be obtained by Hadamard's analysis. It also provides  an asymptotic
expansion
of the distribution $K(t)$
$$K(t)\stackrel{t\to 0}{\sim} \sum_{k=-n}^\infty b_k t^k.$$
Moreover, it gives a detailed information about the distributional
nature of the singularity of $K(t)$ and $K_d(t)$ at $t=0$.
We can compare $b_k$ and $b_{d,k}$ using the fact, that $A_d^2$ is
some sort of analytic continuation of $A^2$.

\subsection{Hadamard's analysis}

We start with recalling the relevant analytical results.
For $k\ge 0$ let $\tau^k$ be the regular distribution on $\R$ given by the
function $\tau^k(t):=t^k$. Define $\tau^{-1}$ by
$$\tau^{-1}(t):=\frac{d}{dt} ln|t|$$
and
$$\tau^{-k-1}:=-\frac{1}{k}\frac{d}{dt} \tau^{-k}.$$

Let $0\in U\subset \R^n$ ($n$ even) be an open subset. We consider an even
(with respect to the $\Z_2$-grading)
second order differential operator $H$ acting on functions on $U$
with values in some $\Z_2$-graded vector space $W$.
There are $End(W)^{ev}$-valued functions $a_{i,j},a_j,a$, $i,j=1,\dots,n$,
such that
$$H:=\sum_{i,j=1}^n a_{i,j} D_iD_j +\sum_{j=1}^n a_j D_j + a,$$
where $D_j:=\partial/\partial x^j$.
We assume, that $H$ is formally selfadjoint and that its principal symbol
is positive, i.e. there is a $C<\infty$ such that for all $x\in U$
and all $\xi=(\xi^1,\dots,\xi^n)\in \R^n$
$$C^{-1} |\xi|^2 \le -\sum_{i,j=1}^n a(x)_{i,j} \xi^i\xi^j\le C|\xi|^2.$$
Let $K(t,x,y)$ be the distributional kernel of the solution operator
of the Cauchy problem
\begin{eqnarray}
(\frac{d^2}{dt^2}+H)f(t,x)&=&0\nonumber \\
f(0,x)&=&f_0(x)\label{cauchy}\\
\frac{d}{dt}_{t=0}f(t,x)&=&0\:\:.\nonumber
\end{eqnarray}
If $f_0$ supported near $0$, the solution of problem (\ref{cauchy}) exists
and is unique for small $|t|$.
In some sense, $K$ is the kernel of $cos(tH)$.
One can use the Hadamard construction in order to obtain approximations
of $K(t,x,y)$ for small $|t|$. The immediate generalization of
H\"ormander \cite{hoermander853}, Thm. 17.5.5 to operators acting on vector
valued functions provides
\begin{kor} \label{kor1}
Let $x\in U$. The distribution $t\rightarrow tr K(t,x,x)$ has an asymptotic
expansion
$$ K(t,x,x):=\sum_{k=-n/2}^\infty c_{2k}(x) \tau^{2k},$$
where the functions $c_{2k}$ are differential polynomials
in the coefficients of the operator $H$.
\end{kor}
A differential polynomial in certain functions is a polynomial in
the various derivatives of these functions.
The differential polynomials in Corollary \ref{kor1}
are homogeneous with respect to a suitable grading.
Let $D_j$ have the degree $1$. We give $H$ the degree $2$. Thus $a_{i,j}$ have
degree $0$, $a_i$ have degree $1$ and $a$ has degree $2$.
The total degree of a monomial in a differential polynomial
is the sum of the degrees of the functions involved and the total number
of differentiations.
The $c_{2k}$ are related to the solutions of the transport equation
\cite{hoermander853}, Eq. 17.4.6. Tracing back degrees we obtain
\begin{kor}
The coefficient $c_{2k}(x)$ of the asymptotic expansion
of $Tr K(t,x,x)$ is given by
a homogeneous differential polynomial in the coefficients
of the operator $H$ of degree $2k+n$.
\end{kor}
Now Assume that the coefficients of the differential operator $H$ are real
analytic.
\begin{lem}
For $x\in U$ the local trace  $K(t,y,y)$ of the fundamental solution
of the problem (\ref{cauchy}) is real analytic in $t$ for small $t>0$
and has a convergent Laurent expansion in $t$ at $t=0$.
\end{lem}
{\it Proof:}
By the method of Hadamard we construct
an approximate fundamental solution $\tilde{K}(t,x,y)$.
It satisfies
the initial condition and
$$(\frac{d^2}{dt^2}+H)\tilde{K}(t,x,y)=:h(t,x,y)$$
is real analytic in $t,x$ for $(t,x)$ inside the forward light cone near the
tip $(0,y)$
(we consider a fixed $y\in U$) and, in fact, extends to a whole neighbourhood
of $(0,y)$. Let $h_{ex}(t,x,y)$ denote this extension.
The local trace of $\tilde{K}(t,x,y)$
is a Laurent polynomial.
We solve the Cauchy problem
$(\frac{d^2}{dt^2}+H)V(t,x,y)=-h_{ex}(t,x,y)$
with zero initial conditions.
Applying the Cauchy-Kowalevskaja theorem
we obtain a unique solution $V(t,x,y)$ which
is real analytic in $(t,x)$ near $(0,y)$.
Then $K(t,x,y)=\tilde{K}(t,x,y)+\chi(t,x,y) V(t,x,y)$, where $\chi$
is the characteristic function of the local forward light cone with tip
$(0,y)$.
In particular, the local trace
of $K(t,y,y)$ near $t=0$ is given by a convergent Laurent expansion
and is real analytic for small $t>0$.
$\Box$

\subsection{The relation between $A_d^2$ and $A^2$}\label{sqw}

We choose coordinates near $X\ni [1]\in \G/\K$ using the exponential map.
Let $\g=\kaaa\oplus\paaa$  be the Cartan decomposition of the Lie algebra $\g$.
We identify $\paaa=\R^{n}$ by specifying a basis in $\paaa$.
The composition of exponential map with the projection
$$E:\paaa\stackrel{exp}{\rightarrow} \G\rightarrow \G/\K=X$$ provides
coordinates
in a neighbourhood of $[1]$. The fibre of the bundle $V$ over $[1]$ can be
identified
in a canonical way with $W=W^+\oplus W^-$ (see (\ref{wdef})).
We trivialize the bundle $V$ in a  neighbourhood
of $[1]$ using the radial parallel transport, i.e. the fibre of $V$ over
$E(p)$, $p\in\paaa$,
is identified with $W$ using the parallel transport along the path $E(tp)$,
$t\in [0,1]$.
We can write the operator $A^2$ in the form
$$A^2:=\sum_{i,j=1}^n a_{i,j} D_iD_j +\sum_{j=1}^n a_j D_j + a,$$
where $a_{i,j}, a_j, a$ are functions on $\paaa$ with values in $End(W)^{ev}$.

The Lie algebra of $\G^d$ is $\g^d=\kaaa\oplus\imath\paaa$. We compose the
exponential
map with the projection in order to obtain coordinates on a neighbourhood of
$[1]\in X_d$:
$$\tilde{E}:\imath\paaa\stackrel{exp}{\rightarrow}\G^d\rightarrow X_d.$$
The fibre of $V_d$ over $[1]$ is again canonically isomorphic to $W$. We employ
the
radial parallel transport in order to trivialize the bundle $V_d$ in a
neighbourhood
of $[1]$. If we multiply the fixed basis of $\paaa$ by $\imath$ we specify a
basis
of $\imath \paaa$. We can write the operator $A_d^2$ in the form
$$A_d^2:=\sum_{i,j=1}^n \tilde{a}_{i,j} \tilde{D}_i\tilde{D}_j +\sum_{j=1}^n
\tilde{a}_j \tilde{D}_j + \tilde{a},$$
where $\tilde{D}_i$ is the derivative $\partial/\partial(\imath x^j)$ and
$\tilde{a}_{i,j}$, $\tilde{a}_j$ and $\tilde{a}$ are $End(W)^{ev}$-valued
functions on $\imath\paaa$.
\begin{lem}\label{extens}
The coefficients of $A^2$ extend as holomorphic functions
to a neighbourhood of $0$ in the complexification $\paaa^c$.
The following relations between the coefficients of $A^2$ and $A^2_d$ hold:
\begin{eqnarray*}
   \tilde{a}_{i,j}(\imath p)&=&a_{i,j}(\imath p)\\
   \tilde{a}_j(\imath p)&=&\imath a_j(\imath p)\\
   \tilde{a}(\imath p)&=&-a(\imath p).
\end{eqnarray*}
\end{lem}
{\it Proof:}
The lemma is a consequence of the fact, that $A^2$ and $-A_d^2$ are
restrictions
to real submanifolds of the same holomorphic differential operator.
Let $\G^c$ and $\K^c$ be the complexifications of the Lie groups. Then
$X^c:=\G^c/\K^c$
is a complex manifold. The representation of $\K$ on $W$ can be extended to
$\K^c$.
The Casimir $\Omega^c$ of $\G^c$ acts on sections of $\G^c\times_{\K^c} W=
V^c\rightarrow X^c$.
Let $(A^c)^2:=\Omega^c-c(\sigma)$.
We use the exponential map of $\G^c$ in order to provide holomorphic
coordinates near $[1]$.
$$E^c:\paaa^c\stackrel{exp}{\rightarrow}\G^c\rightarrow \G^c/\K^c=X^c.$$
We trivialize the bundle $V^c$ using the radial parallel transport.
Then $(A^c)^2$ can be written in the form
$$(A^c)^2:=\sum_{i,j=1}^n a^c_{i,j} D^c_iD^c_j +\sum_{j=1}^n a^c_j D^c_j +
a^c,$$
where $a^c_{i,j}$ are holomorphic, $End(W)^{ev}$-valued functions on $\paaa^c$
and
$D^c_j=\partial/\partial z^j$, where $z^j$ is the complex coordinate
corresponding
to $x^j$. On one hand, by definition $A^2$ is the restriction of $(A^c)^2$ to
$X$. Thus
we obtain
$$a_{i,j}(p)=a^c_{i,j}(p),\quad a_j(p)=a^c_j(p),\quad a(p)=a^c(p).$$
On the other hand, $A^2_d$ is the restriction of $-(A^c)^2$ to $X^d$.
Thus
$$\tilde{a}_{i,j}(\imath p)=a^c_{i,j}(\imath p),\quad \tilde{a_j}(\imath
p)=\imath a^c_j(\imath p),\quad \tilde{a}(\imath p)=-a^c(\imath p)$$
and the Lemma follows.
$\Box$\newline

\subsection{Computation of $K(t)$}

The coefficients $b_k$ of the asymptotic expansion near $t=0$ of
$K(t)$ are zero for odd $k$ and given by  homogeneous differential
polynomials in the coefficients $a_{i,j},a_j,a$ of order $k+n$ for even $k$.
The coefficients $b_{d,k}$ of the asymptotic expansion of $K_d(t)$
are zero for odd $k$ and given by the same differential polynomials
polynomial of degree $k+n$, but now in the coefficients
$\tilde{a}_{i,j},\tilde{a}_j,\tilde{a}$.
{}From Lemma \ref{extens} we obtain
\begin{lem}
The following relation between $b_k$ and $b_{d,k}$ holds:
$$b_{d,k}=i^{k+n} b_k.$$
The distribution $K(t)$ is regular on $\R\setminus \{0\}$ and given there
by
\begin{equation}\label{pots1} K(t)=(-1)^{n/2} K_d(\imath
t)=(-1)^{n/2}\frac{1}{vol(X_d)}\theta_d(t).\end{equation}
\end{lem}
{\it Proof:}
Since $K(t)$ has a convergent Laurent expansion near $t=0$ the equality
(\ref{pots1}) is true for small $|t|$.
Since $K_d(\imath t)$ is real analytic for $t\in\R\setminus\{0\}$ it is enough
to argue that $K(t)$ is real analytic for all $t\not=0$, too.
$\Box$\newline

The nature of $K(t)$ as a distribution is obtained by specifying
the regularization of the meromorphic function $K(t)$ at zero.
The type of regularization follows from Corollary \ref{kor1}.
Recall the proof of Lemma \ref{ll1}. In a similar manner one can prove
\begin{kor}\label{dstre}
In the case $\epsilon=0$ we have
$$\theta_d(t)=-\frac{1}{T}\frac{d^2}{dt^2} Q(\frac{d }{dt })
ln(|sinh(Tt/2)|).$$
In the case $\epsilon=1/2$ we have
$$\theta_d(t)=-\frac{1}{T}\frac{d^2 }{dt^2 } Q(\frac{d }{dt })
ln(|tgh(Tt/4)|).$$
\end{kor}
In view of Corollary \ref{kor1} we obtain
\begin{kor}\label{kkkk}
As a distribution on $\R$ we have in the case $\epsilon=0$
$$K(t)=-(-1)^{n/2}\frac{1}{Tvol(X_d)}\frac{d^2}{dt^2} Q(\frac{d }{dt })
ln(|sinh(Tt/2)|)$$
and in the case $\epsilon=1/2$
$$K(t)=-(-1)^{n/2}\frac{1}{Tvol(X_d)}\frac{d^2}{dt^2} Q(\frac{d }{dt })
ln(|tgh(Tt/4)|).$$
\end{kor}

\subsection{The distributional trace formula}

Let $\phi\in C_c^\infty(\R)$. Then
$$K_{M,\phi}:=\int_{-\infty}^\infty \phi(t) cos(tA_M) dt$$
is of trace class. The linear map
$$C_c^\infty(\R)\ni\phi\rightarrow Tr K_{M,\phi}$$
defines a distribution on $\R$ formally written as $Tr cos(tA_M)$.
Recall the Definition \ref{cd.g} of $C(g)=C(g,\sigma)\in \R$, $l(g)\in\R$
for hyperbolic $g\in G$.
Let $n_\Gamma(g)$ be the number of classes in $\Gamma_g/<g>$, where
$\Gamma_g$ is the centralizer of $g$ in $\Gamma$ and $<g>$ is the group
generated by $g$.
By $C\Gamma$ we denote the set of conjugacy classes of $\Gamma$.

\begin{prop}
The distribution $Tr cos(tA_M)$ is given by
\begin{equation}\label{distrtrace}
vol(M)K(t)+\sum_{[g]\in C\Gamma,[g]\not=[1]}
\frac{C(g)}{n_\Gamma(g)}(\delta(t-l(g))+\delta(-t-l(g))),
\end{equation}
\end{prop}
{\it Proof:}
The operator $K_{M,\phi}$ has a smooth integral kernel
$K_{M,\phi}(x,y)$.
It can be obtained by averaging the integral kernel $K_\phi(\hat{x},\hat{y})$:
$$K_{M,\phi}(x,y)=\sum_{g\in\Gamma} K_\phi(\hat{x},g\hat{y}),$$
where $\hat{x}, \hat{y}$ are lifts of $x,y$ with respect to the covering
$X\rightarrow M$.
We identify the fibre $V_{\hat{x}}$ with $V_{g\hat{x}}$, $g\in\Gamma$
using the homogeneous bundle structure of $V$.
Since $\phi$ has compact support and $cos(tA)$
has unit propagation speed, the operator $K_\phi$ has finite propagation
and thus the sum is uniformly finite, if $x$ and $y$ vary in
a fundamental domain $\Omega_\Gamma$ of $\Gamma$.
We compute
\begin{eqnarray*}
Tr K_{M,\phi}&=&\int_{\Omega_\Gamma}\sum_{g\in\Gamma} tr \: K_\phi(x,gx)dx\\
&=&\int_{\Omega_\Gamma}\sum_{[g]\in C\Gamma}\sum_{h\in[g]}tr\: K_\phi(x,hx)
dx\\
&=&\sum_{[g]\in C\Gamma}\sum_{[h]\in\Gamma/\Gamma_g}\int_{\Omega_\Gamma}tr\:
K_\phi(x,h^{-1}ghx) dx\\
&=&\sum_{[g]\in C\Gamma}\sum_{[h]\in\Gamma/\Gamma_g}\int_{h\Omega_\Gamma}tr \:
K_\phi(x,gx) dx.
\end{eqnarray*}
Since all $g\in\Gamma$ with $g\not=1$
are hyperbolic, $\Gamma_g=\{g^{z/n_\Gamma(g)}\:|\:z\in\Z\}$.
Thus,
$$\cup_{[h]\in \Gamma/\Gamma_g,z\in\Z}\:\:hg^{z/n_\Gamma(g)}\Omega_\Gamma=X$$
and
$$\cup_{[h]\in\Gamma/\Gamma_g}\:\:h\Omega_\Gamma
=\Omega_{<g^{1/n_\Gamma(g)}>},$$
i.e. the fundamental domain of the group generated by $g^{1/n_\Gamma(g)}$.
We obtain
\begin{equation}\label{nummer}Tr K_{M,\phi}=\int_{\Omega_\Gamma}tr \:
K_\phi(x,x)dx +\sum_{[g]\in
C\Gamma,[g]\not=1}\frac{1}{n_\Gamma(g)}\int_{\Omega_{<g>}}tr\:
K_\phi(x,gx)dx.\end{equation}
The first term of (\ref{nummer}) is nothing else than
$$<vol(M) K(t),\phi>.$$
The summands of the second term of (\ref{nummer}) were evaluated in Theorem
\ref{thehyp}.
The result is
$$\frac{1}{n_\Gamma(g)}\int_{\Omega_{<g>}}tr\:
K_\phi(x,gx)dx=\frac{C(g)}{n_\Gamma(g)}(\phi(l(g))+\phi(-l(g))).$$
The proposition follows.
$\Box$\newline
\section{The theta function}

\subsection{The $A$-index theorem}\label{aindex}

Let $\gamma\in R(\K)$ be such that its restriction to $\M$
vanishes, i.e. $r(\gamma)=0\in R(\M)$.
Let $V^c:=V^c(\gamma)\rightarrow X^c$ be the associated $\Z_2$-graded
holomorphic bundle.
By Proposition \ref{p.sur} there exists a $\G^c$-invariant, holomorphic Dirac
operator $D^c$ acting on holomorphic
sections of $V^c$. The bundle $V^c$ restricts to  bundles
$V\rightarrow X$ and $V_d\rightarrow X_d$, and $D^c$ restricts to Dirac
operators
$D$ on V and $\imath D_d$ on $V_d$. Let $L^c$ be an even, $\G^c$-invariant,
holomorphic differential
operator restricting to $L$ on $V$ and $\imath^{q} L_d$ on $V_d$. Here $q$
denotes the order of $L^c$.
We assume that $L^c$ commutes
with $D^c$ and that $L$ and $L_d$ are formally selfadjoint.
Let $\Gamma\subset \G$ be a co-compact subgroup
such that $M:=\Gamma\backslash X$ is a closed smooth manifold.
Let $V_M$ and $D_M$, $L_M$ be the induced bundle and operators on $M$.
Then we can define the $A$-index of $D_M$ and $D_d$.
\begin{ddd}
The $A$-index of $D_M$ is the set of pairs
$$ind_A(D_M,L_M)=\{(\lambda,n_\lambda)|\mbox{$\lambda$ eigenvalue of
$L_{M|ker(D_M)}$ with graded multiplicity $n_\lambda$}\}.$$
Similarly
$$ind_A(D_d,L_d)=\{(\lambda,n_\lambda)|\mbox{$\lambda$ eigenvalue of
$L_{d|ker(D_d)}$ with graded multiplicity $n_\lambda$}\}.$$
\end{ddd}
\begin{theorem}[The proportionality principle for the $A$-index]\label{aaind}
$$ind_A(D_M,L_M)=\{(\imath^{-q} \lambda,\frac{\chi(M)}{\chi(X_d)}n_\lambda)|
(\lambda,n_\lambda)\in ind_A(D_d,L_d)\}.$$
\end{theorem}
{\it Proof:}
For all $N\in \Naaa$ we have
$$\sum_{(\lambda,n_\lambda)\in ind_A(D_M,L_M)}n_\lambda \lambda^N = Tr L_M^N
E_{D_M}\{0\}.$$
Analogously,
$$\sum_{(\lambda,n_\lambda)\in ind_A(D_d,L_d)}n_\lambda \lambda^N = Tr L_d^N
E_{D_d}\{0\}.$$
For any $t>0$
\begin{eqnarray}
Tr L_M^N E_{D_M}\{0\} &=& Tr L_M^N e^{-tD_M^2},\label{ter}\\
Tr L_d^N E_{D_d}\{0\} &=& Tr L_d^N e^{-tD_d^2}\nonumber
\end{eqnarray}
holds.
Now we argue as in Section \ref{sqw}.
We use the Cartan decompositions $\g=\kaaa\oplus\paaa$ and
$\g_d=\kaaa\oplus\imath \paaa$
and the exponential map in order to provide local charts near the origin.
We trivialize the bundles using the radial parallel transport.
Then the operators $D_M^2$ and $D_d^2$ are given by
\begin{eqnarray*}
D_M^2&=&\sum_{i,j=1}^n a_{i,j} D_iD_j +\sum_{j=1}^n a_j D_j + a\\
D_d^2&=&\sum_{i,j=1}^n \tilde{a}_{i,j} \tilde{D}_i\tilde{D}_j +\sum_{j=1}^n
\tilde{a}_j \tilde{D}_j + \tilde{a}.
\end{eqnarray*}
Let $D_i$ have degree one and fix the degrees of the coefficients such that
$D_M^2,D_d^2$ have degree two.
The analogue of Lemma \ref{extens} holds:
\begin{eqnarray*}
   \tilde{a}_{i,j}(\imath p)&=&a_{i,j}(\imath p)\\
   \tilde{a}_j(\imath p)&=&\imath a_j(\imath p)\\
   \tilde{a}(\imath p)&=&-a(\imath p).
\end{eqnarray*}
Thus if $a_I(p)$ is a coefficient of $D^2$, then $\tilde{a}_I(\imath
p)=\imath^{deg(a_I)} a(\imath p)$
is the corresponding coefficient of $D_d^2$.
The operators $L_M$ and $L_d$ also have analytic coefficients and are of the
form
\begin{eqnarray*}
L_M&=&\sum_{|I|\le q} l_I D^I \\
L_d&=&\sum_{|I|\le q} \tilde{l}_I \tilde{D}^I,
\end{eqnarray*}
where $I$ runs over a set of multi-indices.  We fix the degree of the
coefficients
of $L_M,L_d$ such that the total degree of $L_M,L_d$ is $q$.
We have
$$\tilde{l}_I(\imath p)=(\imath)^{q-|I|}l_I(\imath p)=\imath^{deg(l_I)}
l_I(\imath p).$$
The local traces of the operators on the r.h.s. of (\ref{ter})
have asymptotic expansions near $t=0$ (\cite{gilkey84}).
Let $x$ denote some point in $X,X_d$. Then
\begin{eqnarray}
tr L_M^N e^{-tD_M^2}(x,x) &\stackrel{t\to 0}{\sim}& \sum_{k=-n/2-Nq}^\infty t^k
b_k\label{bre}\\
tr L_d^N e^{-tD_d^2}(x,x) &\stackrel{t\to 0}{\sim}& \sum_{k=-n/2-Nq}^\infty t^k
\tilde{b_k}\nonumber,
\end{eqnarray}
where $b_k$ and $\tilde{b_k}$ are the values at $0\in\paaa$ of
universal homogeneous differential polynomials in the coefficients
$a_{i,j},a_i,a,l_I$ and $\tilde{a}_{i,j},\tilde{a}_i,\tilde{a},\tilde{l}_I$
of degree $2k+n+Nq$. It follows
$$\tilde{b_k}=\imath^{2k+n+Nq}\tilde{b_k}.$$
Since the l.h.s. of (\ref{ter}) does not depend on $t$ we obtain
\begin{eqnarray*}
Tr L_M^N E_{D_M}\{0\}&=&vol(M) b_0\\
Tr L_d^N E_{D_d}\{0\}&=&\imath^{n+Nq} vol(X_d) b_0
\end{eqnarray*}
by integrating the local trace over $X,X_d$, respectively.
Using
$$\frac{vol(M)}{vol(X_d)}=\imath^n\frac{\chi(M)}{\chi(X_d)}$$
we conclude
$$Tr L_M^N E_{D_M}\{0\}=\frac{\chi(M)}{\imath^{Nq} \chi(X_d)}Tr L_d^N
E_{D_d}\{0\}$$
and
$$\sum_{(\lambda,n_\lambda)\in ind_A(D_M,L_M)}n_\lambda
\lambda^N=\frac{\chi(M)}{\imath^{Nq} \chi(X_d)} \sum_{(\lambda,n_\lambda)\in
ind_A(D_d,L_d)}n_\lambda \lambda^N.$$
This equation holds for all $N\ge 0$ and thus implies
$$ind_A(D_M,L_M)=\{(\imath^{-q} \lambda,\frac{\chi(M)}{\chi(X_d)}n_\lambda)|
(\lambda,n_\lambda)\in ind_A(D_d,L_d)\}.$$
$\Box$\newline

\subsection{The multiplicities}

In this section we define the theta function $\theta(t)=\theta(t,\sigma)$
associated to an irreducible representation $\sigma\in \hat{\M}$.
We will omit $\sigma$ in our notation but keep $\gamma$ in certain places.
Let $\gamma\in R(\K)$ be admissible for $\sigma$.
We have defined the operator
$A_M(\gamma)$ acting on sections of the graded bundle $V_M(\gamma)\rightarrow
M$.
Let $E_{A_M(\gamma)}$ be its family of spectral projections.
\begin{ddd}
The integers
$$m(\lambda,\gamma):=Tr E_{A_M(\gamma)}\{\lambda\}$$
are called the multiplicities.
\end{ddd}
The reason for introducing the notion of $\sigma$-admissible $\gamma$ is to
obtain the following
\begin{prop}
The multiplicities $m(\lambda,\gamma)$ do not depend on the
choice of the $\sigma$-admissible $\gamma$.
\end{prop}
{\it Proof:}
The proposition is a consequence of the $A$-index theorem that was
proved in Subsection \ref{aindex}.
Let $\gamma,\tilde{\gamma}\in R(\K)$ be two $\sigma$-admissible elements.
Then $r(\gamma-\tilde{\gamma})=0$.
There exists an odd Dirac operator $D_M:=D_M(\gamma,\tilde{\gamma})$ acting
on sections of $V_M(\gamma)\oplus V_M(\tilde{\gamma})^{op}$  (Proposition
\ref{p.sur}).
This Dirac operator commutes with the action of $B_M:=A_M(\gamma)\oplus
A_M(\tilde{\gamma})$.
The eigenspace of $B_M$ corresponding to $\lambda$ splits as
$H(\lambda)=H(\lambda)^0\oplus H(\lambda)^\perp$,
where $H(\lambda)^0$ is in the kernel of $D_M$ and $H(\lambda)^\perp$
is orthogonal to the kernel of $D_M$.
On $H(\lambda)^\perp$ the operator $D_M$ induces an odd isomorphism. Thus
$Tr 1_{H(\lambda)^\perp}=0$ (note that we take the super trace).
Let $ind_A(D_M,B_M^2)=\{(\lambda,n_\lambda)\}.$
Then $Tr 1_{H(\lambda)^0}=n_{\lambda^2}$. We must show that
$ind_A(D_M,B_M^2)=0$.
The operator $D_M$ is a push down of a $\G$-invariant Dirac operator
$D:=D(\gamma,\tilde{\gamma})$ acting on sections of $V(\gamma)\oplus
V(\tilde{\gamma})^{op}$.
The operator $D$ extends to a $\G^c$-invariant, holomorphic operator on
holomorphic
sections of
$V^c(\gamma)\oplus V^c(\tilde{\gamma})^{op}$ which in turn restricts to
$D_d$ on $V_d(\gamma)\oplus V_d(\tilde{\gamma})^{op}$.
Let $B_d=A_d(\gamma)\oplus A_d(\tilde{\gamma})$.
Then $-B_d^2$ is the analytic continuation of $B^2$, where $B^2$ is the lift of
$B_M^2$.
$D_d$ commutes with $B_d^2$.
The $A$-index theorem states that
$$ind_A(D_M,B_M^2)=\{(-\lambda,\frac{\chi(M)}{\chi(X_d)}n_\lambda) |
(\lambda,n_\lambda)\in ind_A(D_d,B_d^2)\}.$$
It is enough to show that $ind_A(D_d,B_d^2)=0$, i.e. all $n_\lambda$ vanish.
Let $(\lambda,n_\lambda)\in ind_A(D_d,B_d^2)$. The eigenspace
$H_d(\lambda)$ of $B_d$ splits into $H_d(\lambda)^0\oplus H_d(\lambda)^\perp$
with respect to $D_d$. On $H_d(\lambda)^\perp$ $D_d$ induces an odd
isomorphism. Thus
$Tr 1_{H_d(\lambda)}=Tr 1_{H_d(\lambda)^0}=n_{\lambda^2}$.
Since $\gamma$ as well as $\tilde{\gamma}$ are $\sigma$-admissible,
it follows $n_\lambda=0$. Hence,
$ind_A(D_d,B_d^2)=0$ and the proposition follows.
$\Box$\newline
Let $\gamma\in R(\K)$ be any $\sigma$-admissible element. We write
$m(\lambda):=m(\lambda,\gamma)$, where  $m(\lambda)$ does not
depend on the choice of $\gamma$.
Since $A_M^2$ is elliptic and of second order, we have the estimate
$$\sum_{\lambda<R} |m(\lambda)| < C R^n, \quad R>>0$$
with $C<\infty$ independent of $R$.
\begin{ddd}
For $Re(t)>0$ define the theta function $\theta(t):=\theta(t,\sigma)$
by
$$\theta(t):=\sum_{\lambda} m(\lambda)e^{-t\lambda}=Tr e^{-tA_M}.$$
\end{ddd}
$\theta(t)$ is holomorphic for $Re(t)>0$.

\subsection{Meromorphic continuation of the theta function}

Recall that we have fixed a $\sigma\in \hat{\M}$.
In this subsection we provide a meromorphic continuation of the theta function
to the whole complex plane and discuss its singularities.
\begin{theorem}
The theta function $\theta(t)$ admits a meromorphic continuation to the whole
complex plane. It satisfies the functional equation
$$\theta(t)+\theta(-t)=2 vol(M)K(\imath
t)=2\frac{\chi(M)}{\chi(X_d)}\theta_d(\imath t).$$
The singularities of $\theta(t)$ are
\begin{itemize}
\item first order poles at
$t=\pm \imath l(g)$, $g\in C\Gamma$ with residue
$$res_{t=\pm \imath l(g)}\theta(t)=\frac{C(g)}{\pi n_\Gamma(g)},$$
\item poles of order $n$ at $t=2\pi k/T$, $k=-1,-2,\dots$
(at these points the singular part coincides with that of $2 vol(M)K(\imath t)$
) and
\item a pole of order $n$ at $t=0$.
\end{itemize}
\end{theorem}
Note that $2\pi k/T$, $k=-1,-2,\dots$ are the lengths of the closed geodesics
of $X_d$.\newline
{\it Proof:}
We use the functional equation in order to define $\theta(t)$ for $Re(t)<0$.
For $Re(t)<0$ let $\theta(t):=-\theta(-t)+2 vol(M) K(\imath t)$.
Then $\theta(t)$ is meromorphic in the half plane $\{Re(t)<0\}$ and has the
singularities claimed in the theorem. Recall the explicit formula for the dual
theta function given in Lemma \ref{ll1}
and that $2vol(M) K(t)=2 \frac{\chi(M)}{\chi(X_d)}\theta_d(t)$.
Obviously, $2 vol(M) K(\imath t)$ is meromorphic.
We have to show, that the two pieces on the left and the right half plane
glue nicely at the imaginary axis to give a global meromorphic function.

We define $\theta$ as a distribution on $U:=\C\setminus \{\frac{-2\pi}{T}
\Naaa\}$.
Let $\phi\in C_c^\infty(U)$. Then by definition
$$<\theta,\phi>=\lim_{\epsilon\to 0} \int_{|Re(v)|\ge \epsilon}
\theta(v)\phi(v) dv.$$
We compute the distributional derivative $\bar{\partial}\theta$.
Note that $\theta(\imath t+\epsilon)$ converges
(considered as a distribution with respect to $t$)
to a tempered distribution on $\R$ when $\epsilon \to 0$.
\begin{eqnarray*}
<\bar{\partial}\theta,\phi>&=&<\theta,\bar{\partial}^\ast \phi>\\
&=&lim_{\epsilon\to 0}\int_\epsilon^\infty\int_{-\infty}^\infty
\frac{1}{2}[-\frac{\partial}{\partial u}-\frac{\partial}{\partial(\imath
t)}](\phi(\imath t+u)+\phi(-\imath t-u))\:\:\theta(\imath t+u)du\:\:dt\\
&+&lim_{\epsilon\to 0}vol(M)\int_\epsilon^\infty\int_{-\infty}^\infty
[-\frac{\partial}{\partial u}-\frac{\partial}{\partial(\imath t)}]\phi(-\imath
t-u)\:\:K(-t+\imath u)dt\:\:du\\
&=&<\frac{1}{2}(\theta(\imath t+0)+\theta(-\imath t+0))-vol(M)K(t),\phi(\imath
t) >.
\end{eqnarray*}
The distributional trace formula states
\begin{eqnarray*}
\theta(\imath t+0)+\theta(-\imath t+0)&=&2 Tr\: cos(t A_M)\\
&=&2 vol(M)K(t)\\
&+&2 \sum_{[g]\in C\Gamma,[g]\not=[1]}
\frac{C(g)}{n_\Gamma(g)}(\delta(t-l(g))+\delta(-t-l(g))).
\end{eqnarray*}
Thus,
$$ <\bar{\partial}\theta,\phi>=\sum_{[g]\in C\Gamma,[g]\not=[1]}
\frac{C(g)}{n_\Gamma(g)}(\phi(\imath l(g)+\phi(-\imath l(g))).$$
Since $\bar{\partial}\frac{1}{z}=\pi \delta(z)$ we obtain that
$\theta(t)$ is holomorphic near the imaginary axis except
at $t=\pm \imath l(g)$ $g\in C\Gamma$, where it has first order
poles with the residues claimed in the Theorem.
The pole at $t=0$ can be discussed as in Duistermaat-Guillemin
\cite{duistermaatguillemin75}.
$\Box$\newline
\section{The logarithmic derivative of the Selberg zeta function}\label{the555}

Let $M:=\Gamma\backslash X$ be a compact rank-one even-dimensional locally
symmetric space of non-compact type.
In this and the next section we develop the theory of the Selberg zeta function
of $M$
from scratch. We first define its logarithmic derivative
as a convergent series on some right half-plane. Employing the
relation of the logarithmic derivative of the Selberg zeta function
with the trace of resolvents we provide a meromorphic extension
and a functional equation. In order to obtain a meromorphic
extension of the Selberg zeta function itself we have to
show, that the residues of the poles of the logarithmic
derivative are integers. This is done by interpreting these
residues as dimensions of eigenspaces of suitable operators.
The Selberg zeta function  $Z_S(t,\sigma)$  will be associated to an
irreducible
representation $\sigma\in\hat{\M}$.
In the present subsection we fix
a $\M$-type $\sigma\in\hat{\M}$ and a $\sigma$-admissible $\gamma$ and omit
$\sigma,\gamma$ in our notation.

\subsection{Definition}

In this subsection we define the logarithmic derivative $D(p)$ of the Selberg
zeta function.
Recall the Definition \ref{cd.g} of the length $l(g)$ and
of the contributions
$$C(g)=-\frac{l(g)e^{\rho l(g)}tr\sigma(m)}{2det(1-Ad(g)_\naaa)},$$
where $\naaa$ comes from the Iwasawa decomposition $\G=\K\A\Naaa$ compatible
with $g$,
and $m\in\M$ is defined by the representation $g=ma$, $a\in\A^+$.
Recall the multiplicity $n_\Gamma(g)$ defined for $g\in\Gamma$.

There is a constant $C<\infty$ such that
$$C(g)<C\:e^{-\rho l(g)}l(g),\quad \forall g\in\Gamma\:\:\:\mbox{hyperbolic}.$$
Moreover, using the volume growth on $X$, one has the estimate
$$\sharp\{[g]\in C\Gamma\:|\: l(g)<N\}\le Ce^{2\rho N}.$$
\begin{ddd} For $Re(p)>\rho$ define the logarithmic derivative of the Selberg
zeta function
by
$$D(p):=\sum_{[g]\in C\Gamma,[g]\not=1}\frac{2C(g)e^{-pl(g)}
}{n_\Gamma(g)}=-\sum_{[g]\in
C\Gamma,[g]\not=1}\frac{l(g)e^{(\rho-p)l(g)}tr\sigma(m)}
{det(1-Ad(g)_\naaa)n_\Gamma(g)}.$$
\end{ddd}
The function $D(p)$ is holomorphic for $Re(p)>\rho$.

\subsection{Products and linear combinations of resolvents}

Let $A$ be some operator. We set $$R(p)=\frac{1}{p+A}.$$
\begin{lem}\label{resprod}
Let $p_1,\dots,p_N\in\C$ be in the resolvent set of $-A$ and $p_i\not= p_j$ for
all $i\not= j$.
Then
$$\prod_{i=1}^NR(p_i)=\sum_{i=1}^N\left(\prod_{j=1,j\not=
i}^N\frac{1}{p_j-p_i}\right)R(p_i).$$
\end{lem}
{\it Proof:}
We prove this equation by induction on $N$ using the resolvent identity
$$(b-a)R(a)R(b)=R(a)-R(b).$$
For $N=1$ the assertion is obvious.
Assume the assertion for $N$. Then
\begin{eqnarray*}
\prod_{i=1}^{N+1}R(p_i)&=&\sum_{i=1}^N
\left(\prod^n_{j=1,j\not=i}\frac{1}{p_j-p_i}
\right)R(p_i)R(p_{n+1})\\
&=&\sum_{i=1}^N\left(\prod^N_{j=1,j\not=i}
\frac{1}{(p_j-p_i)(p_{N+1}-p_i)}\right)(R(p_i)-R(p_{n+1}))\\
&=&\sum_{i=1}^N\left(\prod^{N+1}_{j=1,j\not=i}
\frac{1}{p_j-p_i}\right)R(p_i)-
\sum_{i=1}^N\left(\prod^{N+1}_{j=1,j\not=i}
\frac{1}{p_j-p_i}\right) R(p_{n+1})
\end{eqnarray*}
It remains to show that
$$\sum_{i=1}^N\prod^{N+1}_{j=1,j\not=i}
\frac{1}{p_j-p_i}=-\prod_{j=1}^N\frac{1}{p_j-p_{N+1}}$$
or, equivalently,
$$\sum_{i=1}^{N+1}\prod^{N+1}_{j=1,j\not=i}\frac{1}{p_j-p_i}=0.$$
This is a special case of the Lemma \ref{zwro}. $\Box$\newline
\begin{lem}\label{zwro}
For any $N$-tuple $p_1,\dots,p_N\in\C$ of complex numbers with $p_i\not= p_j$
for all $i\not= j$
and $l=0,1,\dots, N-2$ we have
$$\sum_{i=1}^{N+1}p_i^l\prod^{N+1}_{j=1,j\not=i}\frac{1}{p_j-p_i}=0.$$
\end{lem}
{\it Proof:}
Consider
$$f(p_1):=\sum_{i=1}^{N+1}p_i^l\prod^{N+1}_{j=1,j\not=i}\frac{1}{p_j-p_i}=0$$
as a rational function of $p_1$. It can have poles of first order at $p_1=p_j$,
$j=1,\dots,N$.
Actually, we show that the residues vanish and, hence, $f$ is entire.
The assertion of the lemma follows, since for $l\le N-2$ the function
$f(p_1)$ is bounded and tends to zero at infinity.
We compute
\begin{eqnarray*}
res_{p_1=p_j}f(p_1)&=&lim_{p_1\to p_j}f(p_1)(p_1-p_j)\\
&=&\left(-p_1^l\prod^{N+1}_{k=2,k\not=j}
\frac{1}{p_k-p_1}+p_j^l\prod_{k=2,k\not=j}^{N+1}\frac{1}
{p_k-p_j}\right)_{|p_1=p_j} \\
&=&0\end{eqnarray*}
$\Box$\newline

\subsection{The trace of the product of resolvents and meromorphic continuation
of $D(p)$.}

Recall the Definition \ref{hieradef} of the operator $A_M$.
Let $L<\infty$ be such that $-L^2<A_M^2$ (note that $A_M$ may have imaginary
eigenvalues).
Then for $Re(p)>L$
$$R(p^2)=\frac{1}{p^2+A_M^2}=\int_0^\infty e^{-pt}\frac{sin(tA_M)}{A_M} dt.$$
Let $p_1,\dots,p_N\in\C$ be complex numbers with $p_i\not= p_j$ for all $i\not=
j$ and $Re(p_i)>L$.
Employing Lemma \ref{resprod}
we can write
\begin{eqnarray*}\prod_{i=1}^NR(p_i^2)&=&\sum_{i=1}^N\left(\prod_{j=1,j\not=
i}^N\frac{1}{p^2_j-p^2_i}\right)R(p^2_i)\\
&=&\int_0^\infty\sum_{i=1}^N\left(\prod_{j=1,j\not=
i}^N\frac{1}{p^2_j-p^2_i}e^{-p_it}\right)\frac{sin(tA_M)}{A_M} dt.
\end{eqnarray*}
In order to simplify the notation we will write
$$\Psi(e^{-pt}):=\sum_{i=1}^N\left(\prod_{j=1,j\not=
i}^N\frac{1}{p^2_j-p^2_i}e^{-p_it}\right).$$
The symbol $\Psi$ stands for forming the corresponding linear combination of
the function
being the argument of $\Psi$ with $p$ replaced by the $p_i$.
For example, consider the function $sinh(pt)$. Then
$$\Psi(sinh(pt)):=\sum_{i=1}^N\left(\prod_{j=1,j\not=
i}^N\frac{1}{p^2_j-p^2_i}sinh(p_it)\right).$$
Thus
$$\prod_{i=1}^NR(p_i^2)=\int_0^\infty\Psi(e^{-pt})\frac{sin(tA_M)}{A_M} dt.$$
{}From Lemma \ref{zwro} we deduce:
\begin{lem}\label{taykoe}
Let $a_0+ta_1+t^2a_2+\dots$ be the Taylor expansion of $\Psi(f(pt))$ as a
function of $t$
at $t=0$. Then $a_{2n}=0$ for $n=0,\dots, N-1$.
\end{lem}
Recall that $n=dim(M)$.
For $N\ge n/2+1$ the product
$$\prod_{i=1}^NR(p_i^2)$$
is of trace class.
Recall the distribution $K(t)$ obtained in Corollary \ref{kkkk}
and the Definition \ref{d.constants} of
$\epsilon:=\epsilon(\sigma)\in\{0,1/2\}$.
Let
\begin{eqnarray*}F(t)&:=&ln(|sinh(Tt/2)|),\quad\mbox{case $\epsilon=0$}\\
F(t)&:=&ln(|tgh(Tt/4)|),\quad  \mbox{case $\epsilon=1/2$}.
\end{eqnarray*}
\begin{prop}[$D$ and the trace of resolvents]\label{traceprod} For $N\ge n/2+1$
and
$p_1,\dots,p_N\in \C$ with $Re(p_i)>L$ and $p_i\not= p_j$ for all $i\not= j$
\begin{equation}\label{tool}Tr\prod_{i=1}^NR(p_i^2)
=-\frac{\chi(M)}{T\chi(X_d)}\int_0^\infty
\Psi(pQ(p)e^{-pt})F(t)dt+\Psi(\frac{D(p)}{2p})\end{equation}
holds.
\end{prop}
{\it Proof:}
Let $\phi_\epsilon(t)$ be a smooth cut-off function being zero near $t=0$ and
for $t<0$ and $\phi_\epsilon(t)=1$ for $ t \ge \epsilon>0$.
\begin{lem}\label{tracwe}
\begin{equation}\label{poi}\prod_{i=1}^NR(p_i^2)=lim_{\epsilon\to 0}
\int_0^\infty\phi_\epsilon(t)\Psi(e^{-pt})\frac{sin(tA_M)}{A_M}
dt\end{equation}
in the sense of trace class operators.
\end{lem}
{\it Proof:} Let $Q$ be the projection onto the complement of the
kernel of $A_M$.
It is enough to show that
$$\prod_{i=1}^NR(p_i^2)=lim_{\epsilon\to 0}
\int_0^\infty\phi_\epsilon(t)\Psi(e^{-pt})\frac{sin(tA_M)}{A_M}Q dt$$
in the sense of trace class operators.
Integrating partially $2R$-times, we obtain
$$\int_0^\infty(-1)^{R}\phi_\epsilon(t)\Psi(e^{-pt})\frac{sin(tA_M)}{A_M}Q
dt=\int_0^\infty\left(\phi_\epsilon(t)\Psi(e^{-pt})\right)^{(2R)}
\frac{sin(tA_M)}{(A_M)^{2R+1}}Qdt .$$
Now we can carry out the limes $\epsilon\to 0$ and obtain
$$lim_{\epsilon\to 0}\int_0^\infty
(-1)^{R}\phi_\epsilon(t)\Psi(e^{-pt})\frac{sin(tA_M)}{A_M}Q
dt=\int_0^\infty\left( \Psi(e^{-pt})\right)^{(2R)}
\frac{sin(tA_M)}{(A_M)^{2R+1}}Qdt$$
in the trace class.
Since the application of $\Psi$ kills the even Taylor coefficients of
$\Psi(e^{-pt})$ by Lemma
\ref{taykoe}, if $2R<N$ we can partially integrate back to obtain the desired
result.
$\Box$ (Lemma \ref{tracwe}) \newline
Interchanging the trace and the limes $\epsilon\to 0$ in (\ref{poi}) and
integrating
partially once we obtain
\begin{eqnarray*}
Tr\prod_{i=1}^NR(p_i^2)&=&lim_{\epsilon\to 0} Tr
\int_0^\infty\phi_\epsilon(t)\Psi(e^{-pt})\frac{sin(tA_M)}{A_M} dt\nonumber\\
&=&lim_{\epsilon\to 0} \int_0^\infty Tr \left( \int_t^\infty
\phi_\epsilon(s)\Psi(e^{-ps})ds \right)cos(tA_M )dt.\end{eqnarray*}
We apply the distributional trace formula and obtain
\begin{eqnarray} &&Tr\prod_{i=1}^NR(p_i^2)\nonumber\\
&=&lim_{\epsilon\to 0}  vol(M) \int_0^\infty K(t) \left(\int_t^\infty
\phi_\epsilon(s)\Psi(e^{-ps})ds \right) dt\label{secth}\\
&+&lim_{\epsilon\to 0}  \sum_{[g]\in C\Gamma,[g]\not=1}\frac{C(g)}{n_\Gamma(g)}
\left(\int^\infty_{l(g)}
\phi_\epsilon(s)\Psi(e^{-ps})ds\right).\nonumber\end{eqnarray}
The second term  gives
\begin{eqnarray*}lim_{\epsilon\to 0}\sum_{[g]\in
C\Gamma,[g]\not=1}\frac{C(g)}{n_\Gamma(g)}\left(\int_{l(g)^\infty}
\phi_\epsilon(s)\Psi(e^{-ps})ds\right)
&=&\sum_{[g]\in
C\Gamma,[g]\not=1}\frac{C(g)}{n_\Gamma(g)}\Psi(\frac{e^{-pl(g)}}{p})\\
&=&\Psi(\frac{D(p)}{2p}).
\end{eqnarray*}
The function $f(t):=\int_t^\infty \phi_\epsilon(s)\Psi(e^{-ps})ds$
is not compactly supported. In order to apply the distributional
trace formula one first cuts off at large $|t|$ and extends the function
symmetrically to the negative real axis. Since $f(t)$ is  constant near $t=0$
one obtains a smooth function with compact support. Now the
trace formula will be applied.
In the resulting equation one lets the cut-off tend to infinity.

Inserting the formula \ref{kkkk} for $K(t)$ into the first term of
(\ref{secth})
and using
$$\frac{vol(M)}{vol(X_d)}=(-1)^{n/2}\frac{\chi(M)}{\chi(X_d)}$$
we obtain
\begin{eqnarray*}&&lim_{\epsilon\to 0}  vol(M)\int_0^\infty
K(t)\left(\int_t^\infty \phi_\epsilon(s)\Psi(e^{-ps})ds \right) dt\\
&=&-\lim_{\epsilon\to 0}\frac{\chi(M)}{T\chi(X_d)}\int_0^\infty
\frac{d^2}{dt^2} Q(\frac{d }{dt })\left(\int_t^\infty
\phi_\epsilon(s)\Psi(e^{-ps})ds \right)F(t) dt\\
&=&-\frac{\chi(M)}{T\chi(X_d)}\int_0^\infty  \Psi(pQ(p)e^{-pt})F(t) dt.
\end{eqnarray*}
$\Box$\newline
Next we compute the integral
$$I:=-\frac{\chi(M)}{T\chi(X_d)} \int_0^\infty  \Psi(pQ(p)e^{-pt})F(t) dt.$$
This integral converges for $Re(p_i)>0$, $i=1,\dots,N$ and defines a
holomorphic function.
Integrating by parts we obtain
$$I:=-\frac{\chi(M)}{2\chi(X_d)}\int_0^\infty
\Psi(Q(p)e^{-pt})\left\{\begin{array}{cc}
\frac{cosh(tT/2)}{sinh(tT/2)},&\epsilon=0\\

           \frac{1}{sinh(tT/2)},&\epsilon=1/2 \end{array}\right\}dt .$$
In the case $\epsilon=0$ we use $ctgh(x)=1+\frac{e^{-x}}{sinh(x)}$ in oder to
write
\begin{equation}\label{comb}I
=-\frac{\chi(M)}{2\chi(X_d)}
\Psi(\frac{Q(p)}{p})-\frac{\chi(M)}{2\chi(X_d)}\int_0^\infty
\Psi(Q(p)e^{-t(p+T/2)})\frac{1}{sinh(tT/2)}dt.
\end{equation}
Thus $I$ is a meromorphic function up to $Re(p_i)>-T/2$, $i=1,\dots,N$.

We decompose $I$  into an even and an odd part:
$$I_{ev}:=\frac{1}{2}(I(p_1,\dots,p_N)+I(-p_1,\dots,-p_N)),
I_{odd}:=\frac{1}{2}(I(p_1,\dots,p_N)-I(-p_1,\dots,-p_N)).$$
\begin{lem}
$$I_{odd}=\frac{\chi(M)}{\chi(X_d)}\frac{\pi}{2T}
\Psi\left(Q(p)\left\{\begin{array}{cc} tg(\pi p/T),& \epsilon=1/2\\

   - ctg(\pi p/T),& \epsilon=0 \end{array}\right\}\right).$$
\end{lem}
{\it Proof:}
We first consider the case $\epsilon=1/2$.
We have
$$I_{odd}=\frac{\chi(M)}{2\chi(X_d)}\int_0^\infty
\Psi(Q(p)sinh(pt))\frac{1}{sinh(tT/2)}  dt.$$
Assuming $|p_i|<T/2$, $i=1,\dots,N$ for a moment we can apply the formula
Ryshik/Gradstein \cite{ryshikgradstein57}, 3.311, in order to obtain
$$I_{odd}=\frac{\chi(M)}{\chi(X_d)}\frac{\pi}{2T}\Psi\left (Q(p)tg(\pi
p/T)\right).$$
Now we consider the case $\epsilon=0$.
We must compute the odd part of
$$-\frac{\chi(M)}{2\chi(X_d)} \int_0^\infty\Psi(Q(p)
e^{-t(p+T/2)})\frac{1}{sinh(tT/2)}dt,$$
where we can apply formula Ryshik/Gradstein \cite{ryshikgradstein57}, 3.331,
and obtain
\begin{eqnarray*}&&\frac{\chi(M)}{2\chi(X_d)} \int_0^\infty
e^{-Tt/2}\Psi(Q(p)sinh(pt))\frac{1}{sinh(tT/2)}dt\\&
=&\frac{\chi(M)}{2\chi(X_d)}\Psi(\frac{Q(p)}{p})-
\frac{\chi(M)}{2\chi(X_d)}\frac{\pi}{2T}
\Psi(Q(p)ctg(\pi p/T)).
\end{eqnarray*}
Combining this with equation (\ref{comb}) gives the desired result
$\Box$\newline

We have not found an explicit  formula for the even part $I_{ev}$ though it is
likely that there is a way to carry out the computation.
\begin{lem}
The function $I_{ev}$ is global meromorphic and given by the trace of a
holomorphic
family of trace class operators
$$I_{ev}:=Tr \frac{\chi(M)}{\chi(X_d)}\prod_{i=1}^N \frac{1}{p^2_i-A_d^2}.$$
\end{lem}
{\it Proof:}
Using Corollary \ref{dstre} we can compute:
\begin{eqnarray*}I_{ev}&=&\left(-\frac{\chi(M)}{T\chi(X_d)} \int_0^\infty
\Psi(\frac{d^2}{dt^2}Q(\frac{d }{dt })e^{-pt}/p)F(t) dt\right)_{ev}\\
&=&\left(\frac{\chi(M)}{\chi(X_d)}\int_0^\infty  \Psi(\frac{e^{-pt}}{p}
\theta_d(t)) dt\right)_{ev}\\
&=&\left(\frac{\chi(M)}{\chi(X_d)}\int_0^\infty Tr \Psi(\frac{e^{-(p+A_M)}}{p})
dt\right)_{ev}\\
&=&\frac{\chi(M)}{\chi(X_d)}Tr
\Psi\left(\frac{1}{2p}(\frac{1}{A_d+p}-\frac{1}{A_d-p})\right)\\
&=&\frac{\chi(M)}{\chi(X_d)} Tr \Psi(\frac{1}{p^2-A_d^2})\\
&=&\frac{\chi(M)}{\chi(X_d)} \prod_{i=1}^N \frac{1}{p^2_i-A_d^2}
\end{eqnarray*}
$\Box$\newline

Thus we have found a meromorphic extension of the integral $I=I_{odd}+I_{ev}$
to all of $\C$.
\begin{prop}\label{parta}
The logarithmic derivative of the Selberg zeta function has a meromorphic
extension to all of $\C$ and is given by
\begin{eqnarray}
\Psi(\frac{D(p)}{2p})&=&Tr\prod_{i=1}^N \frac{1}{p_i^2+A_M^2}-I\nonumber\\
&=&Tr\prod_{i=1}^N \frac{1}{p_i^2+A_M^2}\\
&-&\frac{\chi(M)}{\chi(X_d)}\frac{\pi}{2T}
\Psi\left(Q(p)\left\{\begin{array}{cc} tg(\pi p/T)& \epsilon=1/2\\
                                                               -ctg(\pi p/T)&
\epsilon=0 \end{array}\right\}\right)\label{besser}\\
&-&\frac{\chi(M)}{\chi(X_d)} Tr \prod_{i=1}^N \frac{1}{p^2_i-A_d^2}. \nonumber
\end{eqnarray}
This equation is viewed as an equation of functions in $p:=p_1$ while the other
$p_i$, $i>2$, are kept fixed with $Re(p_i)>>0$, $p_i\not\in spec(A_d)$.
The singularities of $D(p)$ are :
\begin{itemize}
\item first order poles at
$p=\pm \imath \lambda$
with residue
$m(\lambda)$ if $\lambda\not=0$ is an eigenvalue of $A_M$,
\item a first order pole at $p=0$ with residue $2m(0)$ if $0$ is an eigenvalue
of $A_M$,
\item first order poles at $p=-\lambda\in T(\Naaa+\epsilon)$ with
residue $ -2\frac{\chi(M)}{\chi(X_d)} m_d(\lambda)$.
Then $\lambda>0$ is an eigenvalue of $A_d$.
\end{itemize}
If two such points coincide, then the residues add up.
\end{prop}
{\it Proof:}
We only have to discuss the singularities.
The poles of the integral $I$ are located at $T(\Naaa+\epsilon)$.
$A_d$ may have more eigenvalues, in fact even imaginary ones, but
we have chosen the $\gamma\in R(\K)$ such that those eigenvalues
cancel out each other. The weighted multiplicity of the eigenvalues of $A_d$
are $m_d(\lambda)=\lambda Q(\lambda)$.
In particular, $I_{ev}$ is regular at $p=0$.
The poles of $I_{ev}$ and $I_{odd}$ cancel each other on the positive real
axis and add up on the negative real axis.
Consider the case $\epsilon=0$.
$I_{odd}$ has poles at $p=-kT$, $k=1,2,\dots$ with
the residue
$$\frac{\chi(M)}{\chi(X_d)} \prod_{i=2}^N \left(\frac{1}{p^2-p_i^2}\right)
\frac{m_d(-kT)}{-2kT}.$$
Thus $D(p)$ has a pole at $p=-kT$ and $I$ contributes to the residue by
$$-2\frac{\chi(M)}{\chi(X_d)}m_d(-kT).$$
In the case $\epsilon=1/2$ the poles of $I_{odd}$ are located at
$(-k+1/2)T$, $k=1,2,\dots$.
$I$ contributes to the residue of $D(p)$ at $p=(-k+1/2)T$ by
$$-2\frac{\chi(M)}{\chi(X_d)}m_d((-k+1/2)T).$$
Let now $\lambda \not=0$ be an eigenvalue of $A_M$ with $m(\lambda)\not=0$.
Then the singular term at $p=\pm\imath\lambda$ in
$$Tr\prod_{i=1}^N \frac{1}{p_i^2+A^2_M}$$
is
$$m(\lambda)\prod_{i=2}^N \frac{1}{\lambda^2+p_i^2} \frac{1}{p^2+\lambda^2}.$$
It contributes to the residue of $D(p)$ by $m(\lambda)$.
The case $\lambda=0$ is discussed similarly.
$\Box$\newline

\begin{prop}[All residues integer]\label{uuuuu}
At all singularities of the function $D(p)$ the residues are integers.
\end{prop}
{\it Proof:}
By Proposition \ref{parta} the residues of $D(p)$
are sums of multiplicities $m(\lambda)$ and
expressions
$$ -2\frac{\chi(M)}{\chi(X_d)}m_d(\lambda).$$
Hence, it is enough to show
\begin{lem}\label{zet}
$$\frac{\chi(M)}{\chi(X_d)}\in \Z.$$
\end{lem}
{\it Proof:}
It is enough to provide an elliptic homogeneous (with respect to $\G^d$)
differential
operator $D_d$ on $X_d$ with $index(D_d)=1$. Then there is a corresponding
homogeneous operator on $X$ inducing $D$ on $M$.
By the proportionality principle (a special case of Theorem \ref{aaind})
$$\frac{\chi(M)}{\chi(X^d)}=\frac{index(D)}{index(D_d)}=index(D)$$
is an integer.
In fact, we show in Theorem \ref{th.ind} that any representation of
$\G^d$ can be realized as the $\G^d$-index of a homogeneous Dirac operator.
In particular, we can realize the trivial representation, which provides an
index
one operator.
$\Box$ (Lemma and Proposition)\newline
\section{The Selberg zeta function}

\subsection{Singularities}

We adopt the same conventions as in Section \ref{the555}.
There, we constructed
the logarithmic derivative $D(p):=D(p,\sigma)$ of the Selberg zeta function
as a global meromorphic function such that all its residues are integers.
Let $p_0\in (\rho,\infty)$ and define
$$C(p_0)=exp(-\int_{p_0}^\infty D(p)dp).$$
The integral converges since for some $C<\infty$ depending on $p_0$
the estimate
$$|D(p)|\le C e^{-(p-\rho)/2},\quad  p\in(p_0,\infty)$$
holds.
Since all residues of $D(p)$ are integers the Selberg zeta function is well
defined by:
\begin{ddd}
The Selberg zeta function associated to $M$ and $\sigma\in\hat{\M}$ is defined
by
$$Z_S(p):=Z_S(p,\sigma):=C(p_0)exp(\int_{p_0}^p D(s) ds).$$
\end{ddd}
The constant $C(p_0)$ is chosen such that $Z_S(p)$ tends to one if $p$
approaches
infinity on the positive real axis.
An immediate consequence of Proposition \ref{parta} is
\begin{theorem}\label{zetasing}
The singularities (zeros have positive and poles negative order) of the Selberg
zeta function are:
\begin{itemize}
\item at
$p=\pm \imath \lambda$ of order
$m(\lambda)$ if $\lambda\not=0$, where $\lambda$ is an eigenvalue of $A_M$,
\item at $p=0$ of order
$2m(0)$ if $0$ is an eigenvalue of $A_M$,
\item at
$p=-\lambda$, $\lambda\in T(\Naaa+\epsilon)$ of order
$-2\frac{\chi(M)}{\chi(X_d)} m_d(\lambda)$.
Then  $\lambda>0$ is an eigenvalue of $A_d$.
\end{itemize}
If two such points coincide, then the orders add up.
\end{theorem}

\subsection{Euler product}
The Selberg zeta function is usually defined by an Euler product.
The following proposition establishes that our Selberg zeta
function coincides with those considered by Fried \cite{fried86}, Gangolli
\cite{gangolli77} and Juhl \cite{juhl93}.

We say that $[g]\in C\Gamma$ is primitive, if $n_\Gamma(g)=1$.
\begin{prop}[Euler product]
For $Re(s)>\rho$ we have
$$Z_S(s)=\prod_{[g]\in C\Gamma,[g]\not=1,{\rm primitive}}\prod_{k=0}^\infty
\:det\left(1-e^{(-\rho-s)l(g)}S^k(Ad(g)_\naaa^{-1})\otimes \sigma(m)\right).$$
 \end{prop}
{\it Proof:}
Let $W\in End(\C^n)$ be some matrix such that $1-W$ is invertible.
Let $a_1,\dots,a_n$ be the eigenvalues of $W$.
Define the $k$'th symmetric power of $W$ by
$$S^kW=W\underbrace{\otimes\dots\times\dots \otimes }_{k}W.$$
Then
$$tr\:S^kW=\sum_{i_1+\dots i_n=k,i_j\ge 1} a_1^{i_1}\dots a_n^{i_n}.$$
Assuming for the moment $a_i\in[0,1)$, $\forall i=1,\dots, n$, we get
\begin{eqnarray*}
\frac{1}{det(1-W)}&=&\prod_{j=1}^n\frac{1}{(1-a_j)}\\
&=&\prod_{j=1}^n\sum_{k=0}^\infty a_j^k\\
&=&\sum_{k=0}^\infty\sum_{i_1+\dots+i_n=k,i_j\ge 1} a_1^{i_1}\dots a_n^{i_n}\\
&=&\sum_{k=0}^\infty tr\: S^kW.
\end{eqnarray*}
Note, that $dim(\naaa)$ is  odd.
We can compute
\begin{eqnarray*}
ln(Z_S(s))&=&\sum_{[g]\in C\Gamma,[g]\not=1}\frac{e^{(\rho-s)l(g)}tr
\sigma(m)}{det(1-Ad(g)_\naaa)n_\Gamma(g)}\\
&=&-\sum_{[g]\in
C\Gamma,[g]\not=1}
\frac{e^{(-\rho-s)l(g)}tr\sigma(m)}
{det(1-Ad(g)^{-1}_\naaa)n_\Gamma(g)}\\
&=&-\sum_{[g]\in C\Gamma,[g]\not=1}
\sum_{k=0}^\infty
e^{(-\rho-s)l(g)}n^{-1}_\Gamma(g) tr S^k(Ad(g)^{-1}_\naaa)tr\sigma(m) \\
&=&-\sum_{[g]\in C\Gamma,[g]\not=1,{\rm primitive}}\sum_{k=0}^\infty
tr\sum_{j=1}^\infty e^{(-\rho-s)jl(g)}j^{-1}  (S^k(Ad(g)_\naaa)^{-j}) \otimes
\sigma(m)^j\\
&=&\sum_{[g]\in C\Gamma,[g]\not=1,{\rm primitive}}\sum_{k=0}^\infty
tr\:ln\left(1-e^{(-\rho-s)l(g)}S^k(Ad(g)_\naaa^{-1})\otimes \sigma(m)\right) \\
&=&\sum_{[g]\in C\Gamma,[g]\not=1,{\rm primitive}}\sum_{k=0}^\infty ln
\:det\left(1-e^{(-\rho-s)l(g)}S^k(Ad(g)_\naaa^{-1})\otimes \sigma(m)\right) \\
\end{eqnarray*}
Exponentiating this equation we obtain the assertion.
$\Box$\newline

\subsection{Functional equation}

\begin{theorem}
The Selberg zeta function satisfies the functional equation
$$\frac{Z_S(s)}{Z_S(-s)}=exp\left(
-\frac{\chi(M)}{\chi(X_d)}\frac{2\pi}{T}\int_{0}^spQ(p)\left\{\begin{array}{cc}
tg(\pi p/T)& \epsilon=1/2\\
                                                           - ctg(\pi p/T)&
\epsilon=0 \end{array}\right\}dp\right)
                                                           . $$
\end{theorem}
{\it Proof:}
If one takes the difference of the representations
of the logarithmic derivative $D(p)$ given in Proposition
\ref{parta} for $p$ and $-p$, then the terms involving resolvents cancel out.
One obtains
$$\frac{D(p)+D(-p)}{2p}=
-\frac{\chi(M)}{\chi(X_d)}\frac{\pi}{T}Q(p)\left\{\begin{array}{cc} tg(\pi
p/T)& \epsilon=1/2\\
                                                           - ctg(\pi p/T)&
\epsilon=0 \end{array}\right\}.
$$
Multiplication by $2p$ and integration from $p_0$ to $p$ leads
the desired formula.
$\Box$\newline

The existence of a functional equation of this type up to an exponential
polynomial
was already pointed out by Juhl \cite{juhl93} and was established by Br\"ocker
\cite{broecker94}.
There is a similar functional equation for Wakayama's Selberg zeta function
\cite{wakayama85}.

\subsection{Representation by regularized determinants}

For $|Re(p_i)|$ large we consider the function $L_M$ of $p_i$, $i=1,\dots,N$ :
$$L_M(p_1,\dots,p_N):=Tr\:\Psi(\frac{1}{p^2+A^2_M}).$$
We can write
\begin{eqnarray*}L_M(p_1,\dots,p_N) &=&Tr\Psi(\int_0^\infty
e^{-t(p^2+A^2_M)}dt)\\
&=&Tr\Psi(\int_0^\infty -\frac{d}{2pdp}e^{-t(p^2+A^2_M)}\frac{dt}{t}) \\
&=&\Psi(-\frac{d}{2pdp} P.F.\int_0^\infty Tr e^{-t(p^2+A^2_M)}\frac{dt}{t}),
\end{eqnarray*}
where $P.F. \int_0^\infty\dots$ stands for taking the finite part of
$\int_\epsilon^\infty\dots $.

Applying Soule-Abramovioch-Burnol-Kramer, \cite{souleabramovichburnolkramer92},
Theorem 5.1.1, we obtain
\begin{equation}\label{ext333}L_M(p_1,\dots,p_N)
=\Psi(\frac{d}{2pdp}(ln(det(p^2+A^2_M))+Eu\: a_M(p^2))).\end{equation}
The determinant is the zeta regularized (super) determinant of $p^2+A^2_M$
defined by
$$ln(det(p^2+A^2_M)):=- P.F.\int_0^\infty Tr e^{-t(p^2+A^2_M)} \frac{dt}{t} -
Eu\: a_M(p^2),$$
where $a_M(p^2)$ is the constant coefficient in the asymptotic expansion of $Tr
e^{-t(p^2+A^2_M)}$
at $t=0$ and $Eu$ is the Euler constant.
Let $\gamma=\sum_{i=1}^lk_i\delta_i\in R(\K)$, $\delta_i\in \hat{\K}$, $k_i\in
\Z$, be
$\sigma$-admissible. Then $A_M(\gamma)$ splits as
$\oplus_{i=1}^l\oplus_{j=1}^{|k_i|}A_M(\delta_i,\sigma)$
and
$$det(A^2_M+p^2)=\prod_{i=1}^l det(A_M^2(\delta_i,\sigma)+p^2)^{k_i}.$$

Equation (\ref{ext333}) extends analytically  to all of $\C$.
In an analogous manner one obtains for $|Im(p_i)|$ large
\begin{eqnarray*}
L_d(p_1,\dots,p_N)&:=&-Tr\Psi(\frac{1}{A^2_d-p^2})\\
&=&\Psi(\frac{d}{2pdp}(ln(det(A^2_d-p^2))+ Eu\:a_d(p^2))),
\end{eqnarray*}
where $a_d(p^2)$ is the constant term in the asymptotic expansion
of $Tr e^{-t(A^2_d-p^2)}$ at $t=0$.
Note that $a_M(p^2)=\frac{\chi(M)}{\chi(X_d)} a_d(p^2)$.
In the following equation derived from (\ref{besser})
the terms involving $a_M$ and $a_d$ cancel out.
\begin{eqnarray*}
\Psi(\frac{D(p)}{2p})&=&-\frac{\chi(M)}{\chi(X_d)}\frac{\pi}{2T}
\Psi\left(Q(p)\left\{\begin{array}{cc} tg(\pi p/T)& \epsilon=1/2\\
                                                               -ctg(\pi p/T)&
\epsilon=0 \end{array}\right\}\right)\\
&+&\Psi(\frac{d}{2pdp} ln(det(p^2+A^2_M))) \\
&-& \frac{\chi(M)}{\chi(X_d)}\Psi(\frac{d}{2pdp}ln(det(A^2_d-p^2))).
\end{eqnarray*}
It follows
\begin{eqnarray*}
D(p)&=&R(p)^\prime-\frac{\chi(M)}{\chi(X_d)}\frac{\pi}{T}
pQ(p)\left\{\begin{array}{cc} tg(\pi p/T)& \epsilon=1/2\\
                                                              - ctg(\pi p/T)&
\epsilon=0 \end{array}\right\}\\
&+&\frac{d}{dp} ln(det(A^2_M+p^2)) \\
&-&\frac{\chi(M)}{\chi(X_d)} \frac{d}{dp}ln(det(A^2_d-p^2)),
\end{eqnarray*}
where $R(p)^\prime$ is a certain odd polynomial.
Integrating once and exponentiating we obtain
that the Selberg zeta function has the representation
\begin{eqnarray}\label{fghf1}Z_S(p)
&=&e^{R(p)}det(p^2+A^2_M)(det(A^2_d-p^2))^{-\chi(M)/\chi(X_d)}  \\
&.&exp\left(-\frac{\chi(M)}{\chi(X_d)}\frac{\pi}{T} \int_{0}^p    sQ(s)
\left\{\begin{array}{cc} tg(\pi s/T)& \epsilon=1/2\\
-ctg(\pi s/T)& \epsilon=0 \end{array}\right\} ds
\right).\nonumber\end{eqnarray}
Note that the exponential of the integral on the r.h.s. of (\ref{fghf1}) is
well
defined since the residues of the integrand are integers by Lemma \ref{zet}.
\begin{lem}\label{lemmm3}
$R(p)=0$
\end{lem}
{\it Proof:}
We thank U. Br\"ocker for suggesting the following argument.
Let $Z^\pm(s,\lambda):=Tr(A_d\pm\lambda)^{-s}$ for $s>>0$ and $\pm\lambda>0$.
Since $$\theta_d(t)=Tr\:e^{-tA_d}\stackrel{t\to 0}{\sim}\sum_{k=-n}^\infty
d_kt^k$$
(note that $d_k=0$ for $k$ odd)
$Z^\pm(s,\lambda)$ have meromorphic continuations with respect to $s,\lambda$
and
$D^\pm(\lambda):=e^{-\frac{d}{ds}Z^\pm(s,-\lambda)}=det(A_d\mp\lambda)$
are well defined entire functions of finite order having zeros at
$\pm\lambda\in L=spec(A_d)$
of multiplicity $m_d(\pm\lambda)$.
The Weierstra\ss{} products
$$\Delta^\pm:=\prod_{\mu\in
L}\left[(1\mp\frac{\lambda}{\mu})exp\left(\sum_{r=1}^n
\frac{(\pm\lambda)^r}{r\mu^r}\right)\right]^{m_d(\mu)}$$
are other entire function of the same order and with the same zeros such that
they differ from $D^\pm(\lambda)$ by exponential polynomials of order $n$,
only.
By results of Voros \cite{voros87} these polynomials are known.
Let $Z^\pm(s):=Z^\pm(s,0)$. Then
\begin{eqnarray}D^\pm(\lambda)&=&\Delta^\pm(\lambda)
\nonumber\\
&&exp\left(-(Z^\pm)^\prime(0)-\sum_{m=1}^n P.F.
Z^\pm(m)\frac{(\pm\lambda)^m}{m}-\sum_{m=1}^nd_{-m}\left(\sum_{r=1}^{m-1}
\frac{1}{r}\right) \frac{\lambda^m}{m!}\right).\label{formel3}\end{eqnarray}
Now
\begin{eqnarray}
\Delta(\lambda)&:=&\Delta^+(\lambda)\Delta^-(\lambda)\label{formel1}\\
&=&\prod_{\mu\in
L}\left[(1-\frac{\lambda^2}{\mu^2})
exp\left(2\sum_{r=2,\mbox{even}}^n
frac{\lambda^r}{r\mu^r}\right)\right]^{m_d(\mu)}\nonumber\\
&=&\prod_{\mu\in
L}\left[(1-\frac{z}{\mu^2})exp\left(\sum_{r=1}^{n/2}\frac{z^r}
{r\mu^{2r}}\right)\right]^{m_d(\mu)} \label{formel2}
\end{eqnarray}
is the Weierstra\ss{} product corresponding to the eigenvalues of $A_d^2$,
where $z=\lambda^2$.
Note that
$$Tr\:e^{-tA^2_d}\stackrel{t\to 0}{\sim}\sum_{k=-n/2}^\infty c_kt^{k}.$$
In general one would expect half interger exponents, too.
But these do not occur in our special case.
Let $Z(s,z)=Tr (A_d^2+z)^{-s}$ and $D(z):=e^{-\frac{d}{ds}
Z(s,-z)}=det(A_d^2-z)$.
Then again by the results of Voros
\begin{eqnarray}D(z)&=&\Delta(z)\nonumber\\
&&exp\left(-Z^\prime(0)-\sum_{m=1}^{n/2} P.F. Z(m)\frac{z^m}{m}
-\sum_{m=1}^{n/2}c_{-m}\left(\sum_{r=1}^{m-1} \frac{1}{r}\right)
\frac{z^m}{m!}\right).\label{formel4}\end{eqnarray}
Note that $Z^+(s)=Z^-(s)=Z(s/2)$.
Combining (\ref{formel1}), (\ref{formel2}), (\ref{formel3}) and (\ref{formel4})
we obtain
\begin{eqnarray*}
\frac{D(\lambda^2)}{D^+(\lambda)D^-(\lambda)}&=&exp\left(
-\sum_{m=1}^{n/2}c_{-m}\left(\sum_{r=1}^{m-1} \frac{1}{r}\right)
\frac{\lambda^{2m}}{m!}\right.\\
&&\left.+2\sum_{m=1,\mbox{even}}^nd_{-m}\left(\sum_{r=1}^{m-1}
\frac{1}{r}\right) \frac{\lambda^m}{m!}\right).
\end{eqnarray*}
We recall the reflection formula proved by Br\"ocker \cite{broecker94}, Lemma
4:
$$\frac{D^+(\lambda)}{
D^-(\lambda)}=
exp\left(-\frac{
\pi}{T} \int_{0}^p
sQ(s)
\left\{
\begin{array}{cc} tg(\pi s/T)& \epsilon=1/2\\
-ctg(\pi s/T)& \epsilon=0 \end{array}\right\} ds \right).$$
We get
\begin{eqnarray*}
D(\lambda^2)&=&D^-(\lambda)^2 exp\left(
-\frac{\pi}{T} \int_{0}^p
 sQ(s)
 \left\{\begin{array}{cc} tg(\pi s/T)&
\epsilon=1/2\\
-ctg(\pi s/T)& \epsilon=0
\end{array}\right\} ds \right.\\
&&-\sum_{m=1}^{n/2}c_{-m}\left(\sum_{r=1}^{m-1} \frac{1}{r}\right)
\frac{\lambda^{2m}}{m!}\\
&&\left. +2\sum_{m=1,\mbox{even}}^nd_{-m}\left(\sum_{r=1}^{m-1}
\frac{1}{r}\right) \frac{\lambda^m}{m!}\right)
\end{eqnarray*}
We now insert this into (\ref{fghf1}) and obtain
\begin{eqnarray}
Z(p)&=&det(p^2+A^2_M)det(A_d+p)^{-2\chi(M)/\chi(X_d)}\label{formel5}\\
&&exp\left(R(p)+ \frac{\chi(M)}{\chi(X_d)}
\sum_{m=1}^{n/2}c_{-m}\left(\sum_{r=1}^{m-1} \frac{1}{r}\right)
\frac{p^{2m}}{m!}
-
2\frac{\chi(M)}{\chi(X_d)}\sum_{m=1,\mbox{even}}^nd_{-m}\left(\sum_{r=1}^{m-1}
\frac{1}{r}\right)
\frac{p^m}{m!}\right).    \nonumber
\end{eqnarray}
We determine $R(p)$ from the asymptotic expansion of the logarithm both sides
as $p\to\infty$.
In fact, $log\:Z(p)$ vanishes exponentially as $p\to \infty$.
Let $$Tr^{-tA_M^2}\stackrel{t\to 0}{\sim}\sum_{k=-n/2}^\infty e_kt^k$$
define the coefficients $e_k$. Note that $e_k$ are connected with
the $c_k$ by the proportionality principle, i.e. $e_k=0$ for $k$ odd and
$$e_k=\frac{\chi(M)}{\chi(X_d)}c_k$$ for $k$ even.
We have
$$log\:det(A_M^2+p^2)\stackrel{t\to 0}{\sim}\sum_{k=n/2}^1
e_{-k}\left(log(p^2)-\left(\sum_{r=1}^{k} \frac{1}{r}\right)\right)
\frac{p^{2k}}{k!}+e_0log(p^2)+o(1)$$
and
$$log\:det (A_D+p)\stackrel{t\to 0}{\sim}\sum_{k=n}^1
d_{-k}\left(log(p)-\left(\sum_{r=1}^{k} \frac{1}{r}\right)\right)
\frac{p^{k}}{k!}+d_0log(p)+o(1).$$
Inserting these expansions into (\ref{formel5}) we get
\begin{eqnarray*}
-R(p)&=&\frac{\chi(M)}{\chi(X_d)}
\left(c_0log(p^2)-2d_0log(p)\right.\\&&+
\left.\sum_{m=1}^{n/2}c_{-m}(log(p^2)-\frac{1}{m})\frac{p^{2m}}{m!}-
2\sum_{m=1,\mbox{even}}^nd_{-m}
(log(p)-\frac{1}{m})\frac{p^m}{m!}\right).
\end{eqnarray*}
Now we claim $\frac{c_{-m}}{m!}=\frac{d_{-2m}}{(2m)!}$. Assuming this claim
$R(p)=0$
immediately follows. To see the claim note that $Z(s)=Z^\pm(2s)$,
$\frac{c_{-m}}{(m-1)!}=res_{s=-m}Z(s)$ and
$\frac{d_{-2m}}{(2m-1)!}=res_{s=-2m}Z^\pm(s)$.
We thank U. Br\"ocker for suggesting this argument.
An immediate consequence of Lemma \ref{lemmm3} is
\begin{prop}\label{k.l}
The Selberg zeta function has the representation
\begin{eqnarray}\label{fghf}Z_S(p)
&=&det(p^2+A^2_M)(det(A^2_d-p^2))^{-\chi(M)/\chi(X_d)}  \\
&.&exp\left(-\frac{\chi(M)}{
\chi(X_d)}\frac{\pi}{T} \int_{0}^p    sQ(s)
\left\{\begin{array}{cc}
tg(\pi s/T)& \epsilon=1/2\\
-ctg(\pi s/T)& \epsilon=0 \end{array}\right\} ds
.\right),\nonumber\end{eqnarray}
 \end{prop}
Proposition \ref{k.l} implies that the Selberg zeta function is completely
determined by the spectral data. This is not surprising in view of the fact
that the geometric data
defining the Selberg zeta function occurs as the residues of the theta
function.


\section{The Ruelle zeta function}\label{ruelleee}

\subsection{Definition and relation with the Selberg zeta function}

We discuss the Ruelle zeta function of compact
rank-one locally symmetric spaces $M$ of non-compact type.
Ruelle \cite{ruelle76} introduced a zeta function for
Anosov flows coding the length spectrum of closed orbits.
In our situation the corresponding flow is the geodesic flow
on the unit sphere bundle of $M$. The closed orbits correspond
to closed geodesics and thus to the conjugacy classes of $\Gamma$.
\begin{ddd}
The Ruelle zeta function of $M$
is defined by the infinite product
$$Z_R(s):=\prod_{[g]\in C\Gamma,[g]\not=1,primitive} (1-e^{-sl(g)})^{-1}$$
converging for $Re(s)>2\rho$.
\end{ddd}
Consider a conjugacy class $[g]\in C\Gamma$.
To $g\in [g]$ we associate
an adapted Iwasawa decomposition $\G=\K\A\Naaa$ of $\G$.
Recall that $\M:=\G_\aaa\cap \K$.
Then $g=ma$ with $m\in \M$ and $a\in \A^+$.
We identify $\A=\R$ using the exponential map.
The Riemannian metric of $X$
induces a metric $A$ and $l(g)$ is the length of $a$ with respect to this
metric,
i.e. the length of the closed geodesic of $M$ corresponding to $[g]$.

It was discovered by Fried \cite{fried86} that the Ruelle
zeta function can be expressed
in terms of Selberg zeta functions.
In the present subsection we provide a short derivation of this fact.
In the following subsection we use this representation in order to discuss
properties of the Ruelle zeta function.

Let $\naaa^c$ be the complexification of the Lie algebra of $\Naaa$.
For $p\ge 0$ consider $\Lambda^p\naaa^c$ as a representation of $\M\A$.
Any $\lambda\in \R$ induces
a one-dimensional representation of $\A$ on $\C_\lambda:=\C$
by $\A\ni a \rightarrow e^{\lambda log(a)}.$
There are sets
$$I_p=\{(\sigma,\chi)|\sigma\in\hat{\M},\chi\in\R\}$$
such that $\Lambda^p\naaa^c$ decomposes with respect to $\M\A$ as
$$\Lambda^p\naaa^c = \sum_{(\sigma,\chi)\in I_p} V_\sigma\otimes \C_\lambda,$$
where $V_\sigma$ is the space of the representation $\sigma$.
We define
$$S(s,p):=\prod_{(\sigma,\chi)\in I_p}Z_S(s+\rho-\chi,\sigma).$$
\begin{prop}\label{prewq}
The Ruelle zeta function has the representation
$$Z_R(s)=\prod_{p=0}^{n-1} S(s,p)^{(-1)^p}.$$
\end{prop}
{\it Proof:}
We employ the following fact.
Let $W$ be a vector space and $A\in End(W)$. Then
\begin{equation} \sum_{l=0}^\infty tr (-1)^l\Lambda^lA =det(1-A),\end{equation}
where $\Lambda^lA$ is the induced operator on $\Lambda^lW$.

Consider $[g]\in C\Gamma$ and adapt the Iwasawa decomposition of $\G$ to
$g=ma$.
Combine the contributions of $[g]$ to all zeta functions occurring in $S(s,p)$.
Then the contribution of $[g]$ to $ln S(s,p)$
is
$$L(s,g,p):=\frac{e^{-sl(g)}}{det(1-Ad(g)_{\naaa^c})n_\Gamma(g)}
\sum_{(\sigma,\chi)\in I_p} e^{\chi l(g)} tr \sigma(m).$$
We can also write
$$L(s,g,p)=\frac{e^{-sl(g)}}{det(1-Ad(g)_{\naaa^c})n_\Gamma(g)} tr \Lambda^p
Ad(g)_{\naaa^c}.$$
Forming the alternating sum over all $p$ we obtain
$$\sum_{p=0}^{n-1} (-1)^p L(s,g,p)=\frac{e^{-sl(g)}}{n_\Gamma(g)}.$$
Thus for $s>>0$
\begin{eqnarray*}
\sum_{p=0}^{n-1} (-1)^p ln S(s,p)&=&\sum_{[g]\in
C\Gamma,[g]\not=1}\frac{e^{-sl(g)}}{n_\Gamma(g)}\\
&=&\sum_{[g]\in C\Gamma,[g]\not=1,primitive} \sum_{k=1}^\infty
\frac{e^{-skl(g)}}{k}\\
&=&-\sum_{[g]\in C\Gamma,[g]\not=1,primitive}   ln(1-e^{-sl(g)})\\
&=&ln Z_R(s)
\end{eqnarray*}
$\Box$\newline

The representation of the Ruelle zeta function provides
its meromorphic continuation to the whole complex plane.
That the Ruelle zeta function has a meromorphic
continuation was observed in Ruelle's original work by dynamical methods
(\cite{ruelle76} ).

\subsection{The functional equation for the Ruelle zeta function}

We employ the functional equations for the Selberg zeta functions
in order to provide a functional equation for the Ruelle zeta function.
First we write
$$Z_R(s)=\prod_{p=0}^{n/2-1}\left(\prod_{(\sigma,\chi)\in
I_p}\frac{Z_S(s+\rho-\chi,\sigma)}{Z_S(s+\chi-\rho,\sigma)}\right)^{(-1)^p}.$$
Here we used the Poincar\'e duality, i.e. for $p\le n/2-1$
$$I_{n-1-p}=\{(\sigma,2\rho-\chi)|(\sigma,\chi)\in I_p\}.$$
Using the functional  equations for the Selberg zeta functions we obtain
\begin{eqnarray}
&&Z_R(s)Z_R(-s)=\prod_{p=0}^{n/2-1}\left(\prod_{(\sigma,\chi)\in
I_p}\frac{Z_S(s+\rho-\chi,\sigma)Z_S(-s+\rho-\chi,\sigma)}
{Z_S(-s+\chi-\rho,\sigma)Z_S(s+\chi-\rho,\sigma)}
\right)^{(-1)^p}\label{fd}\\
&=&exp\left(-\frac{2\pi}{T}\frac{\chi(M)}{\chi(X_d)}
\sum_{p=0}^{n/2-1} (-1)^p
\sum_{(\sigma,\chi)\in I_p}
\int_{s-(\rho-\chi)}^{s+(\rho-\chi)}
P(r,\sigma)\left\{\begin{array}{cc}
tg(\pi r/T)& \epsilon(\sigma)=1/2\\
 -ctg(\pi r/T)& \epsilon
(\sigma)=0\end{array}\right\}
 dr\right). \nonumber
\end{eqnarray}
Let $h(s)$ denote the derivative of the exponent (up to the prefactor), i.e.
\begin{eqnarray*}
h(s)&=&\sum_{p=0}^{n/2-1} (-1)^p \sum_{(\sigma,\chi)\in
I_p}(P(s+(\rho-\chi),\sigma)-P(s-(\rho-\chi),\sigma))\\
&&\quad   \left\{\begin{array}{cc}tg(\pi(s+(\rho-\chi)) /T)&
\epsilon(\sigma)=1/2\\
                               - ctg(\pi(s+(\rho-\chi))/T)&
\epsilon(\sigma)=0\end{array}\right\},
\end{eqnarray*}
where we have used the periodicity of $tg$ and $ctg$ and that
$\rho,\chi\in\frac{1}{2}T\Z$.
Since $\Lambda^p\naaa^c$ is a representation of $\M\A$
we have, using the
definition of $\epsilon(\sigma)$,
$$T\epsilon(\sigma)+\chi+\rho\in T\Z,\quad
\forall (\sigma,\chi)\in I_p.$$
Hence,
$$h(s)=-ctg(\pi s/T)\sum_{p=0}^{n/2-1} (-1)^p \sum_{(\sigma,\chi)\in
I_p}(P(s+(\rho-\chi),\sigma)-P(s-(\rho-\chi),\sigma)) .$$
Let $H(s)$ be the polynomial
$$H(s):=\sum_{p=0}^{n/2-1} (-1)^p \sum_{(\sigma,\chi)\in
I_p}(P(s+(\rho-\chi),\sigma)-P(s-(\rho-\chi),\sigma)).$$
\begin{lem}
The polynomial $H(s)$ is a constant
$$H(s)=\frac{n}{2}\chi(X_d).$$
\end{lem}
{\it Proof:}
The claim is a consequence of identities in the Weyl polynomials.
In order to establish them we employ a geometric argument.
It is enough to evaluate $H(s)$ at a large number of points.

Recall the Borel-Weil-Bott theorem.
Let $\A^d=exp(\imath\aaa)\subset \G^d$. Extend $\C_\lambda$ to a representation
of $\A^d$
for $\lambda\in\frac{1}{2}T\Z$.
Let $\sigma\in\hat{\M}$. If $\lambda\in\R$ is such that
$T\epsilon(\sigma)+\rho+\lambda\in T\Z$, then the representation
$V_\sigma\otimes \C_\lambda$ exists as a representation of $\M\A^d$.
Form the holomorphic vector bundle
$$W(\lambda,\sigma)=\G^d\times_{\M\A^d}(V_\sigma\otimes \C_\lambda)\rightarrow
\G^d/\M\A^d=:B.$$
The Borel-Weil-Bott theorem asserts that the representation
of $\G^d$ with the highest weight $\Lambda=\mu_\sigma+\lambda$ (see Subsection
\ref{chchc} for conventions)
can be realized as the space of holomorphic sections of $W(\lambda,\sigma)$ and
all higher cohomologies of $W(\lambda,\sigma)$ vanish.
A consequence of the theorem of Borel-Weil-Bott is
$$P(\lambda+\rho,\sigma)=\chi_a(B,W(\lambda,\sigma)), \quad
\lambda+\rho+T\epsilon(\sigma)\in\Z,$$
 where $\chi_a$ is the analytic genus of the bundle $W(\lambda,\sigma)$,
i.e. the
Euler characteristic of the complex given by the Dolbeault resolution
of $W(\lambda,\sigma)$.
Thus
$$P(\lambda+\rho,\sigma)=index(\bar{\partial}+\bar{\partial}^\ast),$$
i.e.
$P(\lambda+\rho,\sigma)$ is the index of the Dirac type operator
$\bar{\partial}+\bar{\partial}^\ast$ on a $\Z_2$-graded vector bundle
$\Lambda^{0,\ast}T^\ast B\otimes W(\lambda,\sigma)$
(the grading given by even and odd form degree).
For $p\le n/2-1$ let
$$H_p(s)=\sum_{(\sigma,\chi)\in I_p}P(s+\rho-\chi,\sigma).$$
For $s$ such that
$s-\chi+T\epsilon(\sigma)\in T\Z$, i.e. $s \in T\Z$,
we have
$$H_p(s)=\chi_a(B,\oplus_{(\sigma,\chi)\in I_p} W(-\chi+s,\sigma)).$$
But
$$\G^d\times_{\M\A^d}(\oplus_{(\sigma,\chi)\in I_p} V_\sigma\otimes
\C_{-\chi}\otimes \C_{s })
       =\G^d\times_{\M\A^d}(\Lambda^p(\naaa^c)^\ast\otimes \C_{s }).$$
Now
$$\G^d\times_{\M\A^d}\Lambda^p(\naaa^c)^\ast=\Lambda^{p,0}T^\ast B,$$
and we obtain for $p<n/2-1$
$$H_p(s)=\chi_a(B,\Lambda^{p,0} T^\ast B\otimes W(s ,1)).$$

For $p>n/2-1$ let
$$H_p(s)=\sum_{(\sigma,\chi)\in I_{n-p-1}}P(s-\rho+\chi,\sigma).$$
Applying the consequence of the Borel-Weil-Bott we obtain
$$H_p(s)=\chi_a(B,\oplus_{(\sigma,\chi)\in I_{n-p-1}}W(\chi+s-2\rho,\sigma)).$$
But
\begin{eqnarray*}
\G^d\times_{\M\A^d}(\oplus_{(\sigma,\chi)\in I_{n-p-1}} V_\sigma\otimes
\C_{\chi}\otimes \C_{s-2\rho})
&=&\G^d\times_{\M\A^d}\Lambda^{n-p-1}\naaa^c\otimes \C_{s-2\rho})\\
&=&\G^d\times_{\M\A^d}\Lambda^p(\naaa^c)^\ast\otimes \C_{s })\\
&=&\Lambda^{p,0}T^\ast B\otimes  W(s ,1).
\end{eqnarray*}
Since $H(s)=\sum_{p=0}^{n-1}(-1)^pH_p(s)$ we obtain that $H(s)$ is the
index of $(\bar{\partial}+\bar{\partial}^\ast)$
on the bundle
$$\Lambda^{\ast,\ast}T^\ast B\otimes W(s ,1)=\Lambda^\ast T^\ast B\otimes W(s
,1).$$
Now the differential operator $(\bar{\partial}+\bar{\partial}^\ast)$
can be deformed to $D=\nabla+\nabla^\ast$, where $\nabla$ is induced by the
Levi-Civita connection of $B$ and the homogeneous connection on $W(s+\rho,1)$.
Here $\nabla$ acts by alternating differentiation.
The deformation is given by
$$D_t=(\bar{\partial}+\bar{\partial}^\ast)+t(\partial+\partial^\ast)$$
and stays inside the elliptic operators.
$\partial$ is again defined using the connection.
The index of $D_t$ is independent of $t$ and for $t=1$ it is independent
of the twisting line bundle and equal to $index D_1=\chi(B)$.
The Euler characteristic of $B$ can be expressed in terms of $X_d$
(see Juhl \cite{juhl93})
$\chi(B)=\frac{n}{2}\chi(X_d)$.
$\Box$\newline

Now we come back to discuss the functional equation of the Ruelle zeta
function.
Note that $h(s)=-n/2 \chi(X_d) ctg(\pi s/T)$ is the derivative
of $C-\frac{Tn}{2\pi} \chi(X_d) ln(|sin(\pi s/T)|)$ for any constant $C$.
Hence
\begin{eqnarray}&&\label{p.klkl}\sum_{p=0}^{n/2-1} (-1)^p
\sum_{(\sigma,\chi)\in I_p}  \int_{s-(\rho-\chi)}^{s+(\rho-\chi)}
P(r,\sigma)\left\{\begin{array}{cc}tg(\pi r/T)& \epsilon(\sigma)=1/2\\
 -ctg(\pi r/T)& \epsilon(\sigma)=0\end{array}\right\}
dr\nonumber\\ &=&C-\frac{Tn}{2\pi} \chi(X_d) ln(|sin(\pi s/T)|).\end{eqnarray}
We argue that $C=0$.
We consider the asymptotics of both sides of this equation for
$s\to\imath\infty$.
The r.h.s. behaves as $C-\frac{ n}{2 } \chi(X_d) \imath s + O(e^{\imath cs})$
for some $c>0$.
A first approximation for large $Im(s)$ of the l.h.s is
$$\sum_{p=0}^{n/2-1} (-1)^p \sum_{(\sigma,\chi)\in I_p}
\int_{s-(\rho-\chi)}^{s+(\rho-\chi)} P(r,\sigma) dr  + O(e^{\imath cs}).$$
The integral $$\int_{s-(\rho-\chi)}^{s+(\rho-\chi)} P(r,\sigma) dr$$ is an odd
polynomial in $s$ since   $P(r,\sigma)$ is odd. In particular, its constant
term vanishes.
Thus, the constant term $C$ of the r.h.s. of (\ref{p.klkl}) must vanish, too.

Inserting (\ref{p.klkl})   into  (\ref{fd}) we obtain
\begin{theorem}
The Ruelle zeta function satisfies the functional equation
$$Z_R(s)Z_R(-s)=  sin(\pi s/T)^{n\chi(M)}.$$
In particular, the order of the singularity at
$s=0$ is $\frac{n}{2}\chi(M)$ while the order of the singularities
at $s<-2\rho$, $s\in T\Z$  is $n\chi(M)$.
\end{theorem}

Juhl \cite{juhl93} obtained this functional equation up to an exponential
polynomial.
As in the real and complex hyperbolic case (see Section \ref{exampleee})
it would be interesting to know the set of operators which contribute to the
divisor of the Ruelle
zeta function in the quaternionic hyperbolic and exceptional case.
\section{Examples}\label{exampleee}

\subsection{Forms on the real hyperbolic space}

Let $\G=SO(n,1)$, $\K=SO(n)$ and $\M=SO(n-1)$ for $n=2m$ even.
The corresponding symmetric space of non-compact type is the
real hyperbolic space $H\R^n$ of sectional curvature $-1$.
It's compact dual space is the unit sphere
$S^n$. Let $M:=\Gamma\backslash X$ be a smooth,
closed hyperbolic manifold. We consider the Selberg zeta function associated
to the $\M$-types related to the differential forms on $M$.

Let $\lambda^p:=\Lambda^p\R^n$ be the  exterior powers of the standard
representation of $SO(n)$
and let $\sigma^p:=\Lambda^p\R^{n-1}$ be similarly defined for $SO(n-1)$.
We have the decomposition rules
$$r(\lambda^p)=\sigma^p+ \sigma^{p-1},0<p<n/2,\:\:\:\:\:\:
r(\lambda^0)=\sigma^0,$$
where $r:R(\K)\rightarrow R(\M)$ is the restriction.
We can represent the $\M$-types $\sigma^p$ by restricted $\K$-types as follows:
\begin{eqnarray}
\sigma^0&=&r(\lambda^0)\nonumber\\
\sigma^1&=&r(\lambda^1-\lambda^0)\nonumber\\
\sigma^2&=&r(\lambda^2-\lambda^1+\lambda^0)\nonumber\\
\sigma^p&=&r(\sum_{l=0}^p(-1)^{l-p}\lambda^l),\quad 0<p<n/2\label{dsf}.
\end{eqnarray}
These representations are admissible. In fact, in $L^2(\HR^n,\Lambda^\ast
T^\ast\HR^n)$
discrete series representations only occur on the forms of middle degree.
But these forms are not involved in (\ref{dsf}).
As remarked in the introduction, the absence of the discrete series
(in the weighted sense) on the non-compact side is equivalent to
the fact that the representation of the $\M$-type is admissible.
Thus the singularities of the Selberg zeta functions
$\Z_S(p,\sigma^p)$ are related to the eigenvalues of the Laplace operators
$\Delta_M^j$ on the $j$-forms on $M$ and $S^n$,$ 0\le j\le p<n/2$.
We have $c(\sigma^p)=(\frac{n-1-2p}{2})^2$.
\begin{prop}\label{resel}
The Selberg zeta function $Z_S(s,\sigma^p)$, $p<n/2$ has the following
singularities:
\begin{itemize}
\item a zero at $0\not= s=\imath \lambda\in\imath \R\cup
(-\frac{n-1-2p}{2},\frac{n-1-2p}{2})$ of order
$$dim\{\Delta_M^p\omega=(\lambda^2
+(\frac{n-1-2p}{2})^2)\omega,\:\:\:\delta\omega=0\},$$
\item a zero at $s=0$ of order
$$2dim\{\Delta_M^p\omega
=(\frac{n-1-2p}{2})^2\omega,\:\:\:\delta\omega=0\},$$
\item at $s=\frac{n-1-2p}{2}$ of order
$$\sum_{l=0}^p(-1)^{l-p}b_l(M),$$
\item at $s=-\frac{n-1-2p}{2}$ of order
\begin{equation}\label{tyty}
\sum_{l=0}^p(-1)^{l-p}b_l(M)+(-1)^{p+1}\chi(M)
=\sum_{l=0}^{n-p-1}(-1)^{l+1}b_l(M),
\end{equation}
\item at $\R\ni s=\lambda<-\frac{n-1-2p}{2}$ of order
$$-\chi(M)dim\{\Delta_d^p\omega
=(\lambda^2-(\frac{n-1-2p}{2})^2)\omega,\:\:\:\delta\omega=0\}.$$
\end{itemize}
\end{prop}
There are related results on the Selberg zeta function for differential forms
on real hyperbolic manifolds by Schuster \cite{schuster92}, Fried
\cite{fried862}, Juhl \cite{juhl931}.\newline
{\it Proof:}
We employ Theorem \ref{zetasing}. The Casimir operator $\Omega(\lambda^k)$
coincides with the form
Laplacian $\Delta^k$. Thus
$A(\lambda^k,\sigma^p)^2=\Delta^k-(\frac{n-1-2p}{2})^2$.
The order of the singularity at $s\in\imath \R\cup
(-\frac{n-1-2p}{2},\frac{n-1-2p}{2}]$ is given by
$$\sum_{k=0}^p dim\{\Delta_M^k\omega=(-s^2+(\frac{n-1-2p}{2})^2)\omega\}.$$
The differential $d$ induces isomorphisms of eigenspaces corresponding to
non-zero eigenvalues of different
Laplacians. Taking this into account we obtain the assertion about  the
spectral singularities.
A similar argument works for the topological singularities at
$s<-\frac{n-1-2p}{2}$.
The point $s=\pm \frac{n-1-2p}{2}$ corresponds to harmonic forms. Here $d$
vanishes.
Hence the result for $s=\frac{n-1-2p}{2}$ follows.
At $s=-\frac{n-1-2p}{2}$ we also have the contribution of
the 0'th Betti number of the dual space.  $\Box$\newline

We now combine the results about the Selberg zeta function in order to
discuss the Ruelle zeta function.
Let $\g=\kaaa\oplus\aaa\oplus\naaa$ be an Iwasawa decomposition of $\g$.
Let $\alpha$ be the positive root of  $(\g,\aaa)$ which has the multiplicity
$n-1$.
We have $|\alpha|=1$ and $|\rho|=\frac{n-1}{2}$.
By $\C_r$ we denote the one-dimensional representation of  $\A$ corresponding
to the weight $r\alpha$.
As a representation of $\M\A$ the Lie algebra $\naaa$ can be identified with
$\R^{n-1}\otimes \C_1$.
We obtain $\Lambda^p\naaa=\sigma^p\otimes \C_p$, $p<n/2$ and
$\Lambda^p\naaa=\sigma^{n-p-1}\otimes \C_p$ for $p\ge n/2$.
For $p<n/2$ let $I_p=\{p\}$.
By Proposition \ref{prewq} we obtain
$$Z_R(z)=\prod_{p=0}^{n/2-1}
\left(\frac{Z_S(z+\frac{n-1}{2}-p,\sigma^p)}
{Z_S(z-\frac{n-1}{2}+p,\sigma^p)}\right)^{(-1)^p}.$$
Using Proposition \ref{resel} we can describe the singularities of the
Ruelle zeta function. We demonstrate the contribution of the spectrum
of the form Laplacians.
If $ Im(z) \not=0$, then the order of the singularity of $Z_R(z)$
at the point $z$ is
\begin{eqnarray*}
&&\sum_{p=0}^{n/2-1} (-1)^p \left(ord_{z^\prime=z}
(Z_S(z^\prime+\frac{n-1}{2}-p,\sigma^p) -
Z_S(z^\prime-\frac{n-1}{2}+p),\sigma^p)\right)    \\
&=&\sum_{p=0}^{n/2-1} (-1)^p
dim\{\Delta_M^p\omega
=(\frac{n-1-2p}{2})^2-(z+\frac{n-1-2p}{2})^2)
\omega,\:\:\:\delta\omega=0\}\\
&-&\sum_{p=0}^{n/2-1} (-1)^p
dim\{\Delta_M^p\omega
=(\frac{n-1-2p}{2})^2-(z-\frac{n-1-2p}{2})^2)
\omega,\:\:\:\delta\omega=0\}\\
&=&\sum_{p=0}^{n/2-1} (-1)^p
dim\{\Delta_M^p\omega=(-z^2-(n-1-2p)z)\omega,\:\:\:\delta\omega=0\}\\
&-&\sum_{p=0}^{n/2-1} (-1)^p
dim\{\Delta_M^p\omega=(-z^2+(n-1-2p)z)\omega,\:\:\:\delta\omega=0\}\\
&=&\sum_{p=0}^{n-1} (-1)^p
dim\{\Delta_M^p\omega=(-z(z+n-1-2p))\omega,\:\:\:\delta\omega=0\}.
\end{eqnarray*}
Thus the eigenvalue $\lambda>\frac{(n-1-2p)^2}{4}$ of the Laplacian $\Delta^p$,
restricted to the
co-closed forms, gives a contribution at
$z=-\frac{n-1-2p}{2}\pm\sqrt{(\frac{n-1-2p}{2})^2-\lambda }.$
Here $Re(z)=- \frac{n-1-2p}{2}$.

Finally we discuss the singularity at the point $z=0$. We apply again
Proposition \ref{resel} and
obtain
\begin{eqnarray*}ord_{s=0} Z_R(s) &=&
\sum_{p=0}^{n/2-1}(-1)^p \left( \sum_{j=0}^p(-1)^{j-p}
b_j(M)-\sum_{j=0}^{n-p-1}(-1)^{j+1}b_j(M) \right)\\
&=&\sum_{p=0}^{n/2-1} \left( \sum_{j=0}^p(-1)^jb_j(M)+\sum_{j=p+1}^{n}(-1)^j
b_j(M) \right)\\
&=&\frac{n}{2} \chi(M).
\end{eqnarray*}
This reproves in the real hyperbolic case that the order of the
singularity of $Z_R(z)$ at $z=0$ is $\frac{n}{2} \chi(M)$.

\subsection{The Dirac operator on the real hyperbolic space}\label{dirara}

In order to consider the Dirac operator we take
$\G=Spin(n,1)$, $\K=Spin(n)$ and $\M=Spin(n-1)$, $n$ even.
If $\Gamma\subset \G$  is a discrete co-compact
subgroup  without elliptic elements, then $M=\Gamma\backslash\G/\K$
is a closed  spin manifold. Let $S^{n-1}$ be the spinor module
of $\M$. The spinor module $S^n$ of
$\K$ is reducible and splits $S^n=\Delta^n_+\oplus S^n_-$.
We have $r(S^n_\pm)=S^{n-1}$, where $r:R(\K)\rightarrow R(\M)$
is the restriction map. It is known that discrete series representations
of $\G$ can not occur in  $L^2(\HR^n,V(S^n_\pm))$ ( see Seifarth
\cite{seifarth93}).
Hence $S_\pm^n$ are $S^{n-1}$-admissible. The square of the
Dirac operator on $X$ can be expressed in terms of the Casimir
operator of $\G$ and the scalar curvature $\tau$ of $X$ using the Parthasarathy
formula \cite{parthasarathy72}:
$$D^2=\Omega(S^n)+\frac{\tau}{8}.$$
Let $D^2_\pm$ be the restrictions of $D^2$ to $C^\infty(X,V(S_\pm^n))$.
Then $D^2_\pm$ coincides with $A^2(S_\pm^n, S^{n-1})$
up to an additive constant.
By the calculation  in \cite{bunke902} we have in fact
$$A^2(S_\pm^n,S^{n-1})=D^2_\pm.$$
In fact, for both operators the essential spectrum starts at $0$ and
has constant multiplicity with respect to the principal series representations
of $\G$
occurring in $L^2(\HR^n,V(S_\pm^n))$.

Now we discuss the singularities of the Selberg zeta function associated
to the $\M$-type $ S^{n-1}$.
We normalize the metric of $X=\HR^n$ such that is has constant sectional
curvature $-1$.
\begin{kor}
The singularities of $Z_S(p,S^{n-1})$ are
\begin{itemize}
\item zeros at $p= \imath\lambda$ of order $dim E_{D }\{\lambda\}$,
where $ \lambda$ is an eigenvalue of $D$,
 \item singularities of order $-\chi(M)\left(\begin{array}{c} n-1-j
\\j\end{array}\right)2^{n/2-1}$
at $p=-(\frac{n}{2}+j)$, $j=0,1,2,\dots$.
\end{itemize}
\end{kor}
For the topological singularities on the negative real axis we employed
the computation of the spectrum of the Dirac operator on the sphere $S^n$
by Sulanke \cite{sulanke79}.

Note that the Selberg zeta function $Z_S(p,S^{n-1})$ for
even-dimensional hyperbolic spin manifolds $M$
is completely determined by the spectrum of the Dirac operator of $M$.

\subsection{Differential forms on complex hyperbolic manifolds}\label{formske}
Let $\G:=SU(n,1)$,
$\K:=S(U(1)\times U(n))$ and $\M:=U(n-1)_2\stackrel{2-fold}{\rightarrow}
U(n-1)$.
Consider $SU(n,1)$ as the subgroup of matrices of $GL(n+1,\C)$ with determinant
one fixing the form
$|z_0|^2-|z_1|^2-\dots-|z_n|^2$.
$S(U(1)\times U(n))$ consists of matrices $diag(q,A)$,
$A\in U(n) $, $q\in U(1)$ with $q\:det(A)=1$ and $\M$ is the subgroup
of $\K$ of matrices of the form $(q,q,B)$ with $B\in U(n-1)$ and $q^2det(B)=1$.
The two-fold covering of $U(n-1)$ is given by $(q,q, B)\rightarrow B$.

Let ${}^\M V $ be the representation of $\M$ on $\C^{n-1}$ given
$(q,q,B)\mapsto q^{-1}B$.
The complexification of ${}^\M V$ splits as ${}^\M V^c={}^\M V^{1,0}\oplus
{}^\M V^{0,1}$.
Taking appropriate alternating powers we obtain the representations
${}^\M \Lambda^{p,q}$.

Let $\omega=\sum_{i=1}^{n-1} z^i\wedge \bar{z}^i\in{}^\M\Lambda^{1,1}$ be the
K\"ahler "form".
The orthogonal complement of
$\omega{}^\M\Lambda^{p-1,q-1}\subset{}^\M\Lambda^{p,q}$
consists of the primitive $(p,q)$-forms and is denoted by $\Lambda_\M^{p,q}$.
The representations $\Lambda_\M^{p,q}$ are irreducible.
Similarly one defines ${}^\K\Lambda^{p,q}$ and $\Lambda_\K^{p,q}$.
Let $r:R(\K)\rightarrow R(\M)$ be the restriction map.
\begin{lem}\label{gaill}
For $p,q$ with $p+q\le n-1$
$$r(\Lambda_\K^{p,q})=\sum_{p-1\le s\le p,q-1\le r\le q} \Lambda_\M^{r,s}.$$
\end{lem}
(compare Gaillard \cite{gaillard88}).
By convention let $\Lambda^{p,q}_.=\{0\}$ if $p$ or $q$ is negative.
\begin{lem}\label{rerere}
For $p,q\in\Naaa$ such that $p+q \le n-1$ we have
$$\Lambda_\M^{p,q}=r(\sum_{r,s\ge 0} (-1)^{r+s} \Lambda_\K^{p-s,q-r}).$$
\end{lem}
This Lemma is an easy consequence of Lemma \ref{gaill}.
It is known that discrete series representation of $\G$ on
$L^2(\HC^n,\Lambda^\ast T^\ast\HC^n)$
can only occur on the forms of middle degree.
Thus for $p+q\le n-1$ the $\K$-type given in Lemma \ref{rerere} representing
$\Lambda_\M^{p,q}$
is admissible.

Let $\taaa\subset \om$ be the standard Cartan algebra given by diagonal
matrices of the form
$diag(q,q,x_1,\dots,x_{n-1})$, $q,x_i\in\imath\R$, $2q+\sum_{i=1}^{n-1}x_i=0$.
We can write weights as $n+1$-tuples of complex
numbers such that the first two entries coincide and the total sum vanishes.
The positive roots of $\M$ are $\{E_{ij}\}$, $2\le i<j \le n$, where the only
non-trivial entries of $E_{ij}$ are $1$ at the $i$'th and $-1$ at in $j$'th
place.
The half-sum $\rho_\om$ of the positive roots is
$\rho_\om=\frac{1}{2}(0,0,n-2,n-4,\dots,4-n,2-n)$.

The weight of the $i$'th coordinate $z^i$ of $V^{1,0}$ is
$(-1/2,-1/2,0,\dots, 0,1,0,\dots,0 )$, where $1$ is at the $(i+2)$'th  entry.
Analogously the weight of $\bar{z}^i$ is $(1/2,1/2,0,\dots,0,-1,0,\dots,0)$.
For $p+q\le n-1$
the highest weights of ${}^\M\Lambda^{p,q}$ and $\Lambda_\M^{p,q}$ coincide
and are given by
$$\mu_{\Lambda_\M^{p,q}}=(\frac{q-p}{2},\frac{q-p}{2},
\underbrace{1,\dots,1}_{p-times},0,\dots,0,
\underbrace{1,\dots,1}_{q-times}).$$
We compute the shift constants $c(\Lambda_\M^{p,q})$.
We use the invariant scalar product on $\g$  given by the trace norm
$(X,Y):=\frac{1}{2}tr\: XY$.
The metric induced on $\taaa^\ast$ is simply
$|(x_0,\dots,x_n)|^2=2\sum_{i=0}^n |x_i|^2$.
We obtain
$$|\rho_\om|^2-|\rho_\om+\mu_{\Lambda_\M^{p,q}}|^2= (p+q)^2-2n(p+q).$$
The algebra $\aaa$ can be  chosen as the set of matrices
$$diag(\left(\begin{array}{cc}
0&r\\r&0\end{array}\right),\underbrace{0,\dots,0}_{n-1\:\:- times}),\quad
r\in\R.$$
A linear form on $\aaa$ can be written in a compatible way as $(n+1)$-tuple
$(x,-x,0,\dots,0)$, $x\in \C$. Fixing a positive direction in $\aaa$ the short
root of $(\g,\aaa)$ is $\alpha=\frac{1}{2}(1,-1,0,0,\dots,0)$.
It has multiplicity $2n-2$ while the multiplicity of $2\alpha$ is one.
The half-sum $\rho$ of the positive roots of $(\g,\aaa)$
is $\frac{1}{2}(2(n-1)\alpha+2\alpha)=\frac{1}{2}(n,-n,0,\dots,0)$ and we have
$|\rho|^2=n^2$.
\begin{lem}
For $p+q\le n-1$ we have
$$c(\Lambda_\M^{p,q})= (p+q-n)^2.$$
 \end{lem}
It follows that
\begin{equation}\label{aass}
A(\Lambda_\K^{p,q},\Lambda_\M^{r,s})^2=\Delta^{p,q}_\circ-
(r+s-n)^2,\end{equation}
where $\Delta^{p,q}_\circ$ is the restriction of the $p+q$-form Laplacian to
the primitive $(p,q)$-forms.

Let $M$ be a closed locally symmetric complex hyperbolic manifold.
The choice of the trace norm on $\g$ induces a scale of the Riemannian
metric of $M$ such that the sectional curvatures are between
$-1 $ and $-4$.
Let $\omega$ denote the K\"ahler form of $M$ and $\delta$
be the adjoint of the deRham differential $d$.
By $h^{r,s}_\circ(M)$ we denote the dimension of the primitive
$(p,q)$-cohomology of $M$, i.e.
$$h^{r,s}_\circ(M)=dim(H^{p,q}(M)/\omega
H^{p-1,q-1}(M))=dim\:ker\Delta^{p,q}_\circ$$
Combining (\ref{aass})  with Lemma \ref{rerere} and Theorem \ref{zetasing}
we obtain
\begin{prop}\label{kaesel}
The Selberg zeta function $Z_S(z,\Lambda_\M^{p,q})$, $p+q\le n-1$,
has the following singularities :
\begin{itemize}
\item at $z=\pm \imath \lambda\not=0$, $\lambda\in\R$ or $\lambda\in\imath\R$
with $|\lambda|<n-p-q$, of order
      $$dim\{\alpha\in \C^\infty(M,\Lambda_\K^{p,q}T^\ast M)| \Delta_\circ
\alpha=(\lambda^2+ (p+q-n)^2 )\alpha, \delta\alpha=0\},$$
\item at $z=0$ of order
      $$2dim\{\alpha\in \C^\infty(M,\Lambda_\K^{p,q}T^\ast M)|
\Delta_\circ\alpha= (p+q-n)^2 \alpha,\delta\alpha=0\},$$
\item at $z= n-q-p$ of order $$\sum_{r,s}(-1)^{r+s}h^{p-r,q-s}_\circ(M) ,$$
\item at $z= q+p-n$ of order
$$\sum_{r,s}(-1)^{r+s}h^{p-r,q-s}_\circ(M)-2 (-1)^{p+q}\frac{\chi(M)}{n+1} ,$$
\item at $z=-\lambda$, $\R\ni\lambda\ge n-(p+q)$ of order
      $$dim\{\alpha\in \C^\infty(\PC^n,\Lambda_\K^{p,q}T^\ast \PC^n)|
\Delta_\circ\alpha=(\lambda^2- (p+q-n)^2 )\alpha,\delta\alpha=0\}.$$
\end{itemize}
\end{prop}
We have employed the isomorphisms induced by $d$
of the non-zero eigenspaces of
the different Laplacians on the bundles corresponding to the $\K$-types
occurring in Lemma \ref{rerere}.
Note that if $\alpha$ is co-closed and primitive, then $d\alpha$ is primitive,
too.
If $\alpha$ is primitive, then so is $\delta\alpha$ (see Wells \cite{wells73}
4.3.10.)

We consider the Ruelle zeta function.
The Lie algebra $\naaa$ splits as $\naaa=\naaa_\alpha\oplus\naaa_{2\alpha}$
with respect to the action of $\M\A$. We can identify $(\naaa_\alpha)^c={}^\M
V^c\otimes \C_\alpha$ and
$(\naaa_{2\alpha})^c$ with $\C_{2\alpha}$, where $\C_\beta$ is the
one-dimensional representation
of $\A$ with weight $\beta$.
We abbreviate $\C_{p\alpha}=\C_p$.
For $k\le n-1$  we have
\begin{eqnarray*}
\Lambda^k(\naaa)^c&=&{}^\M\Lambda^k\otimes \C_k\oplus {}^\M\Lambda^{k-1}\otimes
\C_{k+1}\\
&=&\sum_{p+q=k,k-1}\sum_{r,s\ge 0} \Lambda_\M^{p-r,q-s}\otimes \C_{2k-p-q} \\
\Lambda^{2n-1-k}(\naaa)^c&=&(\Lambda^{k}(\naaa)^c)^\ast\otimes \C_{2n}\\
&=&\sum_{p+q=k,k-1}\sum_{r,s\ge 0} \Lambda_\M^{p-r,q-s}\otimes \C_{2n-2k+p+q}.
\end{eqnarray*}
Let
$$I_k:=\{(p,q,r,s) |p,q,r,s\in\Naaa, p+q=k\: \:or\:\: p+q=k-1, 0\le r\le p,0\le
s\le q\}$$
By Proposition \ref{prewq} we obtain the following representation of the Ruelle
zeta
function in terms of Selberg zeta functions:
$$Z_R(z)=\prod_{k=0}^{n-1}\left(\frac{\prod_{(p,q,r,s)\in I_k}
Z_S(z+n-2k+p+q,\Lambda_\M^{p-s,q-r})}{\prod_{(p,q,r,s)\in I_k}
Z_S(z+2k-n-p-q,\Lambda_\M^{p-s,q-r})}\right)^{  (-1)^k}.$$

Using  Proposition \ref{kaesel} one can describe the singularities
of the Ruelle zeta function in terms of the spectrum of the form Laplace
operators.
In order to demonstrate the contribution of the differential form spectrum
we consider the singularities at $z$ for $ Im(z) \not=0$ large enough.

The singularity at $z$ has the order
\begin{eqnarray*}
&&\sum_{k=0}^{n-1} (-1)^k (\sum_{(p,q,r,s)\in I_k} ord_{z^\prime=z}
Z_S(z^\prime+ n-2k+p+q,\Lambda_\M^{p-s,q-r})\\
&&\quad -ord_{z^\prime=z} Z_S(z^\prime+ 2k-n-p-q ,\Lambda_\M^{p-s,q-r} ))\\
&=&\sum_{k=0}^{n-1} (-1)^k \sum_{(p,q,r,s)\in I_k}
dim\{\Delta^{p-r,q-s}_\circ\alpha=( (p+q-n)^2 -(z+ n-2k+p+q
)^2)\alpha,\delta\alpha=0\}\\
&-&\sum_{k=0}^{n-1} (-1)^k \sum_{(p,q,r,s)\in I_k}
dim\{\Delta^{p-r,q-s}_\circ\alpha=( (p+q-n)^2 -(z+ 2k-n-p-q
)^2)\alpha,\delta\alpha=0\}\\
&=&\sum_{k=0}^{n-1} (-1)^k \sum_{p+q=k,k-1} dim\{\Delta^{p,q} \alpha=(
(p+q-n)^2 -(z+ n-2k+p+q )^2)\alpha,\delta\alpha=0\}\\
&-&\sum_{k=0}^{n-1} (-1)^k \sum_{p+q=k,k-1} dim\{\Delta^{p,q} \alpha=(
(p+q-n)^2 -(z+ 2k-n-p-q )^2)\alpha,\delta\alpha=0\}\\
&=&\sum_{k=0}^{n-1} (-1)^k
dim\{\Delta^k\alpha=((n-k)^2-(z+(n-k)^2))\alpha,\delta\alpha=0\}\\
&-&\sum_{k=0}^{n-2} (-1)^k
dim\{\Delta^k\alpha=((n-k)^2-(z+(n-k-1)^2))\alpha,\delta\alpha=0\}\\
&-&\sum_{k=0}^{n-1} (-1)^k
dim\{\Delta^k\alpha=((n-k)^2-(z-(n-k)^2))\alpha,\delta\alpha=0\}\\
&+&\sum_{k=0}^{n-2} (-1)^k
dim\{\Delta^k\alpha=((n-k)^2-(z-(n-k-1)^2))\alpha,\delta\alpha=0\}
\end{eqnarray*}
Thus, an eigenvalue $\lambda>(n-k)^2$ of the $k$-form Laplacians, $k\le n-1$,
contributes to the singularity of the Ruelle zeta function
at $z=- \epsilon (n-k) \pm\sqrt{ (n-k)^2 -\lambda }$ with sign
$(-1)^k\epsilon\in\{-1,1\}$. In particular, $Re(z)=-\epsilon ( n-k) $.
If $k<n-1$ then there is also a contribution to the singularity at
$z=-\epsilon (n-k-1) \pm\sqrt{ (n-k)^2 -\lambda}$ with sign
$(-1)^{k+1}\epsilon$.
Here $Re(z)=-\epsilon (n-k-1) $.
Note that the pattern of the contribution of the $k$-form spectrum to the
singularities of the Ruelle zeta function differs from the real hyperbolic
case.
Note that Fried \cite{fried88} employs zeta functions associated to
$(p,q)$-forms
on complex hyperbolic manifolds
in order to express the analytic $\bar{\partial}$-torsion.

\subsection{The bundle $K_M^{-1/2}$ over complex hyperbolic manifolds}

We study the Selberg zeta function associated to
a square root $\sqrt{K_M}$ of the canonical bundle of a complex hyperbolic
manifold.

Let $\G,\K,\M$ as in Subsection \ref{formske}.
Assume that $n$ is odd and $n>1$.
The surface case $n=1$ is
contained in the Dirac operators case on real hyperbolic manifolds considered
in Subsection \ref{dirara}.
We consider the $\M$-type $k/2$ being the one-dimensional representation
$(q,q,B)\rightarrow q^{(n+1)/2}$.
Let $K/2 $ be the one-dimensional representation $(q,A)\rightarrow
q^{(n+1)/2}$.
Then $r(K/2)=k/2$, $r:R(\K)\rightarrow R(\M)$.
Let $M$ be a locally symmetric complex hyperbolic manifold.
The line bundle $V_M(K/2)$ is a  square root of the canonical bundle
$K_M=\Lambda^{max}T_c^\ast M$. It defines a spin structure and an associated
Dirac operator over $M$.
This Dirac operator is nothing else than the total Dolbeault operator
$D:=\bar{\partial}+\bar{\partial}^\ast$ on the spinor bundle
$S_M=V_M(K/2)\otimes(\oplus_{p=0}^n\Lambda^{0,p}T^\ast M)$.
It is known that $L^2(\HC^n,S_{\HC^n})$ does not contain discrete series
representations
of $\G$ (see Seifarth, \cite{seifarth93}).
In particular, $L^2(\HC^n,V(K/2))$ does not contain discrete series
representations
and hence $K/2$ is $k/2$-admissible.
The Dirac operator $D^2$ splits into a direct sum of $D^2_p$, $p=0,\dots,n$,
with respect to the form degree.
The operator $D_0^2$ is responsible for the singularities of the Selberg zeta
function associated to $k/2$.
We use again the Parthasarathy formula
$D^2=\Omega(S_M)+\tau/8$
in order to compare $D_0^2$ with
$A(K/2,k/2)^2$.

On $\g$ we again fix the invariant scalar product given by
$(X,Y)=\frac{1}{2}tr(XY)$
inducing on $M$ a metric with sectional curvatures between $-1$ and $-4$.
The scalar curvature of $M$ is then $\tau=-4n(n+1)$.

We fix the Cartan algebra of $\g$ as in Subsection \ref{formske}
and employ the notation for weights explained there.
The weight of the representation $k/2$ is given by
$$\mu=\frac{1}{2}(\frac{n-1}{2},\frac{n-1}{2},-1,\dots,-1).$$
We obtain $|\mu|^2=\frac{n^2-1}{4}$ and $<\rho_\om,\mu>=0$.
Using $|\rho|^2= n^2 $ it follows
$c(k/2)=\frac{3n^2+1}{4}$
We obtain
$$
A^2(K/2,k/2)=D_0^2-\frac{(n-1)^2}{4}.
$$
\begin{kor}
The singularities of the Selberg zeta function $Z_S(p,k/2)$ are
\begin{itemize}
\item zeros of order $dim E_{D_0^2}\{\lambda\}$ at
$0\not= p\in [-\frac{(n-1)}{2 },\frac{(n-1)}{2 }]\cup\imath\R$ with
$p^2=\frac{(n-1)^2}{4}-\lambda $
\item a zero of order $2dim E_{D_0^2}\{\frac{(n-1)^2}{4}\}$ at $p=0$
\item  zeros of order $-2\frac{\chi(M)}{n+1}dim
E_{D^2_{d,0}}\{p^2-\frac{(n-1)^2}{4}\}$
at the points $p<0$ with $p= 2 l+\frac{3n+1}{2 }$, $l\ge \frac{n+1}{ 2 } $.
\end{itemize}
\end{kor}
In order to describe the topological singularities of the zeta function
we used the computation of the spectrum of
$D_{d,0}^2=\bar{\partial}^\ast\bar{\partial}$
on $V_d(K/2)\rightarrow \PC^2$,
 i.e. of the Dolbeault operator on the square root the canonical bundle of
$\PC^n$, found in
Seifarth/Semmelmann \cite{seifarthsemmelmann93}.

\subsection{One-forms on  quaternionic hyperbolic manifolds}

In the case of quaternionic hyperbolic manifolds let
$\G=Sp(1,n)$. It is the subgroup of $Gl(n+1,\Haaa)$
fixing the quaternionic scalar product $|q_0|^2-|q_1|^2-\dots-|q_n|^2$.
Its maximal compact subgroup is $\K=Sp(1)\times Sp(n)$.
Then $X=\HH^n=\G/\K$ is the quaternionic hyperbolic space. We consider only the
case $n\ge 2$, since $n=1$ is contained in the hyperbolic case.
On $\HH^n$ fix the metric which has sectional curvatures between $-1$ and $-4$.
Then $\rho=2n+1$ and the short root $\alpha$ of $(\g,\aaa)$ has length one.
Here $\aaa$
is defined by an Iwasawa decomposition of $\g=\kaaa\oplus\aaa\oplus\naaa$ as
usual.
We consider $\K\subset \G$ as the subset of matrices  $diag(q,A)$,
$q\in Sp(1)$, $A\in Sp(n)$  of $\G$. Moreover, we let $\M=Sp(1)\times Sp(n-1)$.
We consider an $\M$-types occurring in the representation of $\M$ on $\naaa$.
This $\M$-type is interesting since it contributes to the singularities of the
Ruelle zeta function.
Already this example shows a rather complex
behavior. This complexity eventually prevents us from a complete understanding
of
the singularities of the Ruelle zeta function.

Let $\lambda^1$ be the isotropy representation of $\K$. Thus $\lambda^1$ is
the representation of $\K$ on $\Haaa^n$ given by $\Haaa^n\ni  x\rightarrow
Ax\bar{q}$,
$diag(q,A)\in\K$. Then $r(\lambda^1)=\sigma^1\oplus\sigma^\prime\oplus 1$,
$r:R(\K)\rightarrow R(\M)$ being the restriction map.
Here $1$ is the trivial representation, $\sigma^\prime$ is the $3$-dimensional
representation of $\M$ on the root space $\naaa_{2\alpha}$ corresponding to the
long root
and $\sigma^1$ is the representation of $\M$ on the root space $\naaa_\alpha$
corresponding to the short root.
 Identify $\naaa_\alpha$ with $\Haaa^{n-1}$
and $\naaa_{2\alpha}$ with $\R^\perp\subset \Haaa=(\Haaa^{n-1})^\perp\subset
\Haaa^n$.
Then the representation $\sigma^1$ is given by
$x\rightarrow Bx\bar{q}$, $x\in\Haaa^{n-1}$, $diag(q,q,B)\in\M$ and
$\sigma^\prime$ is given by $z\rightarrow qz\bar{q}$,
$z\in\R^\perp$. There is an analogous representation $\lambda^\prime$ on
$\R^\perp$
of $\K$ given by $z\rightarrow qz\bar{q}$, $z\in\R^\perp$, $diag(q,A)\in\K$.
Obviously, $r(\lambda^\prime)=\sigma^\prime$.

The bundle $V(\lambda^1)$  is the bundle of one-forms on $\HH^n$. The bundle
$V(\lambda^\prime)$ can be viewed as a subbundle of $\Lambda^2T^\ast\HH^n$.
The fibre of $V(\lambda^\prime)$  is the space
of K\"ahler forms distinguished by the quaternionic K\"ahler structure.

\begin{lem}
The representation $\lambda^1-\lambda^\prime-1$ is $\sigma^1$-admissible.
 \end{lem}
{\it Proof:} It is known that  discrete series representations of $\G$
in $L^2(\HH^n,\Lambda^pT^\ast\HH^n)$ only occur on the harmonic forms in the
middle degree.
$\Box$\newline

In the normalization of the metric with sectional curvatures in $[-4,-1]$
the shift constant of $\sigma^1$ is $c(\sigma^1)=4n^2$.
Hence,  $A(1,\sigma^1)^2=\Delta_0-4n^2$, where $\Delta_0$ is the function
Laplacian, $A(\lambda^1,\sigma^1)^2=\Delta_1-4n^2$, where
$\Delta_1$ is the one-form Laplacian and
$A(\lambda^\prime,\sigma^1)=\Delta_{2|V(\lambda^\prime)}-4n^2$,
where $\Delta_{2|V(\lambda^\prime)}$ is the restriction of the two-form
Laplacian
to $V(\lambda^\prime)$. We have used the fact that the $\G$-Casimir operator
on differential forms coincides with the Laplacian.

We define the differential operator
$B:C^\infty(\HH^n,\Lambda^1T^\ast\HH^n)\rightarrow
C^\infty(\HH^n,V(\lambda^\prime))$
by $B\alpha=\pi d\alpha$, where $\pi:\Lambda^2 T^\ast \HH^n\rightarrow
V(\lambda^\prime)$
is the orthogonal projection.

Let $M$ be a closed locally symmetric quaternionic hyperbolic manifold.
Then we have the induced bundle $V_M(\lambda^\prime)\rightarrow M$
and $B:C^\infty(M,\Lambda^1T^\ast M)\rightarrow
C^\infty(M,V_M(\lambda^\prime))$.
\begin{prop} The singularities of the Selberg zeta function $Z_S(p,\sigma^1)$
are
\begin{itemize}
\item zeros of order
$$dim\{ \alpha\in C^\infty(M,\Lambda^1T^\ast M) | \Delta_1
\alpha=(\lambda^2+4n^2)\alpha,\delta\alpha=0,B\alpha=0\}$$
at $  p= \imath\lambda$, $\lambda\in (\R\cup \imath (-2n,2n))\setminus
\imath\{2n-4,0,-2n+4\}$,
\item  a zero of order
$$2dim\{\alpha\in C^\infty(M,\Lambda^1T^\ast M) | \Delta_1
\alpha=(\lambda^2+4n^2)\alpha,\delta\alpha=0,B\alpha=0\}$$
at $p=0$,
\item a singularity at $p=\pm(2n-4)$ of order
\begin{eqnarray*}
&&dim\{\alpha\in C^\infty(M,\Lambda^1T^\ast M) | \Delta_1  \alpha= 16(n-1)
\alpha,\delta\alpha=0\}\\&&-dim\{\omega\in C^\infty(M,V_M(\lambda^\prime)) |
\Delta_2 \omega=16(n-1)\omega\},\end{eqnarray*}
\item a pole of order $b_0(M)-b_1(M) =1$ at $p=2n$
\item a zero of order
$b_1(M)-b_0(M)+\frac{2\chi(M)}{n+1}=\frac{2\chi(M)}{n+1}-1$ at $p=-2n$
\item poles of order
$$2\frac{\chi(M)}{n+1} dim\{ \alpha\in C^\infty(\PH^n,\Lambda^1T^\ast \PH^n) |
\Delta_1  \alpha=(\lambda^2-4n^2)\alpha,\delta\alpha=0,B\alpha=0\}$$
at $\lambda=\{-2n-4,-2n-6,\dots\}$.
\end{itemize}
\end{prop}

\begin{lem}\label{nokoge}
$ker\: B^\ast\subset \{\omega\in C^\infty(M,V_M(\lambda^\prime)) | \Delta_2
\omega=16(n-1)\omega\}$
\end{lem}
{\it Proof:}
First we discuss the behaviour of $B^\ast$ on $\HH^n$.
Let $H^{\sigma^\prime,z}$ be the hyperfunction globalization of the principal
series representation of $\G$ associated to
$\sigma^\prime$ and $z\in\C\cong \aaa^\ast_c$.
Then there are unique Poisson transforms (Olbrich, \cite{olbrichdiss94})
$$P^{\lambda^1}_{\sigma^\prime,z}:
H^{\sigma^\prime,z}\rightarrow
\{\alpha\in C^\infty(\HH^n,\Lambda^1T^\ast \HH^n) | \Delta_1  \alpha=(
(2n+1)^2+8-z^2 ) \alpha\}$$
and
 $$P^{\lambda^\prime}_{\sigma^\prime,z}:H^{\sigma^\prime,z}\rightarrow
\{\alpha\in C^\infty(\HH^n,V(\lambda^\prime) ) | \Delta_1  \alpha=(
(2n+1)^2+8-z^2 ) \alpha\}.
$$
 We have
$$
B^\ast P^{\lambda^\prime}_{\sigma^\prime,z}  = (z+ 2n-3)
P^{\lambda^1}_{\sigma^\prime,z} .
 $$
For $Re(z)\ge 0$ the Poisson transform $P^{\lambda^\prime}_{\sigma^\prime,z}$
is surjective.
Hence all eigenforms of $\Delta_{2|V(\lambda^\prime)}$ arise as Poisson
transforms.
Moreover, $P^{\lambda^1}_{\sigma^\prime,z}$ is injective for $Re(z)\ge 0$,
$z\not=2n-3$.
Therefore, if $\omega$ is an eigenform of $\Delta_2$ and $B^\ast \omega=0$,
then $\Delta_2\omega= 16(n-1)\omega$.
In particular
\begin{equation}\label{noharmo}
\{\omega\in C^\infty(\HH^n,V(\lambda^\prime)) | \delta\omega=d\omega=0\}=\{0\}
\end{equation}
Note that $C^\infty(M,V_M(\lambda^\prime))$ is the space
of $\Gamma$-invariant sections of $C^\infty(\HH^n,V(\lambda^\prime))$.
Any form in $C^\infty(M,V_M(\lambda^\prime))$ can be decomposed
into eigenforms of $\Delta_2$. The claim of the Lemma follows. $\Box$\newline

{\it Proof of the Proposition:}
We first discuss the spectral singularities at $p\in  \imath\R\cup (-2n,2n)$.
The operator $d:C^\infty(M)\rightarrow C^\infty(M,\Lambda^1T^\ast M)$
induces an isomorphic embedding of the eigenspaces of $A(1,\sigma^1)$ into the
corresponding eigenspaces of $A(\lambda^1,\sigma^1)$ for all eigenvalues
$\lambda\not= \imath 2n$.
The kernel of $d$ is the space of constant functions and corresponds to
the eigenvalue $\imath 2n$ of $A(1,\sigma^1)$. The orthogonal complement
of the image of $d$ is characterized by $\delta\alpha=0$.

By Lemma \ref{nokoge}  the operator $B^\ast$ induces an isomorphic embedding
of the eigenspaces
of $A(\lambda^\prime,\sigma^1)$ into the  corresponding eigenspaces of
$A(\lambda^1,\sigma^1)$ for all eigenvalues $\lambda\not= \imath(2n-4)$ of
$A(\lambda^\prime,\sigma^1)$.
At $\lambda\not= \imath(2n-4)$ one has to subtract
$$dim\{\omega\in C^\infty(M,V_M(\lambda^\prime)) | \Delta_2
\omega=16(n-1)\omega\}$$
by hand.
Since $im(B^\ast)\perp im(d)$ the assertion about the spectral singularities
follows easily.

The point $p=2n$ corresponds to harmonic forms.
Here $d$ fails to be injective and one has to subtract $1=b_0(M)$.
On $V_M(\lambda^\prime)$ there are no harmonic forms (see (\ref{noharmo})).
Thus the order of the singularity of the Selberg zeta function
$Z_S(p,\sigma^1)$
is $b_1(M)-1$. Furthermore, combining Borel-Wallach \cite{borelwallach80}
II.7.8 and VII.6.1 it follows that $b_1(M)=0$. Thus the Selberg zeta function
has
a pole of first order at $p=2n$.

At $p=-2n$ there is also a contribution from the spectrum of $A_d(1,\sigma^1)$.
Since there are no  harmonic one and two forms on $\PH^n$, we have
$m_d(2n,\sigma^1)=-1$. The order of the singularity of $Z_S(p,\sigma^1)$ at
$p=-2n$
follows.
The argument for the remaining singularities on the negative real axis is
similar
to the one of the spectral singularities.
$\Box$\newline
 \newcommand{\oc}{{\bf c}}

\section{Harmonic analysis on compact symmetric spaces}

In this section we collect some facts of representation-theoretic nature which
we need to define the theta functions. Because these results are of independent
interest we state them in greater generality than needed in the present paper.

\subsection{The Cartan-Helgason theorem for vector bundles}\label{chchc}

Let $X:=\G/\K$ be a Riemannian symmetric space of the non-compact type.
Let $\Paaa=\M\A\Naaa$ be the Langlands decomposition of a minimal parabolic
subgroup, $\G=\K\A\Naaa$ and $\g = \kaaa \oplus \aaa \oplus \naaa$ the
corresponding Iwasawa decompositions of $\G$ and its Lie algebra $\g$,
$\Phi^+(\g,\aaa)$ be the positive root system corresponding to
$\Naaa$ and $W=W(\g,\aaa)$ be the Weyl group. Let
$\rho:=\frac{1}{2}\sum_{\alpha\in\Phi^+(\g,\aaa)}m_\alpha \alpha$, where
$m_\alpha$ is the dimension of the root space corresponding to $\alpha$. We
assume that $\G$ is contained as a real form in its complexification $\G^c$. We
define $\G^d$ to be the analytic subgroup of $\G^c$ corresponding to
$\g^d:=\kaaa\oplus\paaa^d\subset\g^c$, $\paaa^d:=\imath\paaa$, $\g=\kaaa\oplus
\paaa$ being the Cartan decomposition of $\g$, and $X_d:= \G^d/\K$ the dual
compact symmetric space. We choose a maximal abelian subalgebra $\taaa$ of
$\om$. Then $\haaa := \taaa^c \oplus \aaa^c$ is a Cartan subalgebra of $\g^c$.
We choose a positive root system $\Phi^+(\g^c,\haaa)$ having the property that
for $\alpha\in\Phi(\g^c,\haaa)$ $\alpha_{|\aaa}\in \Phi^+(\g,\aaa)$ implies
$\alpha \in \Phi^+(\g^c,\haaa)$. Let
$\delta:=\frac{1}{2}\sum_{\alpha\in\Phi^+(\g^c,\haaa)}\alpha$, and set
$\rho_m:= \delta-\rho$. We choose a $W$-invariant scalar product $(.,.)$ on
$\aaa$ and define the root vector $H_\alpha \in \aaa$ for $\alpha \in
\Phi^+(\g,\aaa)$ by
\[ \lambda (H_\alpha) := \frac {(\lambda,\alpha)}{(\alpha,\alpha)}\quad \forall
\lambda \in \aaa^* \ .\]
Note that $exp (2\pi\imath H_\alpha)\in \M$.

Let $\gamma\in \hat \K$ be fixed. Then we can decompose
$\gamma_{|\M}=\bigoplus_{\sigma \in \hat \M} [\gamma:\sigma]\sigma$. We are
interested in the multipicities $[\pi_\Lambda :\gamma]$ for an irreducible
representation $\pi_\Lambda$ with highest weight $\Lambda \in \haaa^*$
occurring in the decomposition of
\[ L^2(\G^d\times_\K V_\gamma) = \bigoplus_{\pi_\Lambda \in \hat\G^d}
[\pi_\Lambda :\gamma] \pi_\Lambda\ .\]
\begin{theorem}\label{th.CH}
Let $[\pi_\Lambda :\gamma]\not=0$. Then
\begin{enumerate}
\item $\Lambda_{|\taaa}$ equals the highest weight $\mu_\sigma$ of a
representation $\sigma$ of $\M$ with $[\gamma:\sigma]\geq [\pi_\Lambda
:\gamma]$;
\item
$ \frac {(\Lambda_{|\aaa},\alpha)}{(\alpha,\alpha)} \in \varepsilon_\alpha
(\sigma) + \Naaa \quad \forall \alpha \in \Phi^+(\g,\aaa)$,
where $\varepsilon_\alpha (\sigma)\in \{0,1/2\}$ is defined by
$$ e^{2\pi \imath\varepsilon_\alpha (\sigma)}:= \sigma (exp (2\pi\imath
H_\alpha)) \in \{\pm 1\}\ ;$$
\item
$ \frac {(\Lambda_{|\aaa},\alpha)}{(\alpha,\alpha)}\geq \max_{\{\beta \in
\Phi^+(\g^c,\haaa)| \beta_{|\aaa}=\alpha\}} \frac
{(\mu_\sigma,\beta)}{(\alpha,\alpha)} =:\Phi_\alpha (\sigma)\quad \forall
\alpha \in \Phi^+(\g,\aaa)\ .$
\end{enumerate}
Moreover, for $ (\sigma,\gamma)\in \hat \M \times \hat \K$ with
$[\gamma:\sigma]\not=0$ there exist constants $C_\alpha (\sigma,\gamma)\geq
\Phi_\alpha (\sigma)$ such that for $\lambda \in \aaa^*$ with
\begin{equation}\label{f.eps}
\frac {(\lambda,\alpha)}{(\alpha,\alpha)} \in \varepsilon_\alpha (\sigma) +
\Naaa
\end{equation}
and
$$\frac {(\lambda,\alpha)}{(\alpha,\alpha)}\geq C_\alpha (\sigma,\gamma)\quad
\forall \alpha \in \Phi^+(\g,\aaa) $$
and for $\Lambda := \lambda + \mu_\sigma \in \haaa^*$, where $\mu_\sigma$ is
the highest weight of $\sigma$, the following equation holds
\[ [\pi_\Lambda :\gamma] = [\gamma :\sigma] \ .\]
\end{theorem}
{\it Proof }: Note that by holomorphic continuation $\pi_\Lambda$ extends to an
irreducible representation of $\G^c$, hence of $\G$, acting on the same vector
space $V_\pi$. Let $v_\Lambda$ be the highest weight vector of $\pi_\Lambda$.
Then $v_\Lambda$ is the unique highest weight vector for the action of $\M$ in
the $\naaa$-invariants of $V_\pi$ denoted by $V_\pi^\naaa$. This gives us the
desired irreducible representation $\sigma$ of $\M$. In fact, the Iwasawa
decomposition implies that $v_\Lambda$ is a cyclic vector for the action of
$\K$ in $V_\pi$. Therefore the mapping
\begin{equation}\label{f.F}
F: Hom_\K(V_\pi,V_\gamma)\longrightarrow Hom_\M (V_\sigma,V_\gamma)
\end{equation}
defined by
\[ F(T) := T\circ i\ ; \quad i: V_\sigma \rightarrow V_\pi^\naaa\ , \]
where $V_\sigma$ and $V_\gamma$ are representation spaces for $\sigma$ and
$\gamma$ respectively, is injective. This proves 1.

In order to prove 2., we recall that (see Helgason
\cite{helgason78},Ch.VII,Thm.8.5.)
\begin{equation}\label{f.cap1}
\K\cap exp \imath \aaa = exp \{2\pi\imath\sum_{\alpha\in\Phi^+(\g,\aaa)}
n_\alpha H_\alpha\ |\ n_\alpha \in \Naaa \}\subset \M\ .
\end{equation}
All these elements are of order at most 2. Hence
\[ \pm v_\Lambda = \sigma (exp (2\pi\imath H_\alpha))v_\Lambda = \pi(exp
(2\pi\imath H_\alpha))v_\Lambda = e^{2\pi \imath
\frac{(\Lambda_{|\aaa},\alpha)}{(\alpha,\alpha)}} v_\Lambda\ .\]

3. is a simple consequence of the $\Phi^+(\g^c,\haaa)$-dominance of $\Lambda$
and the fact that for $\beta \in \Phi^+(\g^c,\haaa)$ also $-\theta \beta$ is in
$\Phi^+(\g^c,\haaa)$. Here, $\theta$ is the Cartan involution of $\g$.

We come to the second part of the theorem. Given a pair $(\sigma,\lambda)\in
\hat \M\times \aaa^*$ with
$\frac{(\Lambda_{|\aaa},\alpha)}{(\alpha,\alpha)}\geq \Phi_\alpha(\sigma)
\:\:\forall \alpha\in\Phi^+(\g,\aaa)$ we form the $\Phi^+(\g^c,\haaa)$-dominant
element $\Lambda:= \lambda + \mu_\sigma \in \haaa^*$. We claim that if
$\lambda$ also satisfies the condition (\ref{f.eps}), then $\Lambda$ is
integral. We have to show that for $H=H_\aaa + H_\taaa \in \imath \aaa \oplus
\taaa$ with $exp(H)=1$ also $e^{\Lambda(H)}=1$. Indeed, since
$exp(H_\aaa)=exp(H_\taaa)^{-1}\in \M$, we obtain for the highest weight vector
$v_\sigma \in V_\sigma$
\[ e^{\Lambda (H)}v_\sigma=e^{\lambda(H_\aaa)} e^{\mu_\sigma(H_\taaa)}v_\sigma
= \sigma (exp(H_\aaa))\sigma(exp(H_\taaa))v_\sigma =v_\sigma\ .\]

Thus, $\Lambda$ is the highest weight of a representaion $(\pi_\lambda,V_\pi)$
of $\G^d$ such that the irreducible representation of $\M$ on $V_\pi^\naaa$ has
highest weight $\mu_\sigma$. This shows that this representation coincides with
$\sigma$ on the connected component $\M_0$ of $1\in\M$. Using (\ref{f.eps}) and
the fact that
\[ \M = (\K\cap exp\ \imath\aaa) \cdot \M_0
\qquad\mbox{(\cite{helgason78},Ch.IX,ex.A3)} \]
we see that this representation is given by $\sigma$ on the whole group $\M$.
We are done, if we can show that the map $F$ defined in (\ref{f.F}) is
surjective for $\lambda$ big enough.

The dual map
\[ F^*: Hom_\M(V_\gamma,V_\sigma)\longrightarrow Hom_\K(V_\gamma,V_\pi) \]
applied to $t \in Hom_\M(V_\gamma,V_\sigma)$ is given by
\[ F^*(t) = \int_\K \pi(k^{-1}) i\circ t \gamma(k) \ dk\ .\]
Indeed, for $T\in Hom_\K(V_\pi,V_\gamma)$ we have
\[ tr(T\circ F^*(t)=tr(T \int_\K \pi(k^{-1}) i\circ t \gamma(k) \
dk)=tr(\int_\K \gamma(k^{-1}) T\circ i\circ t \gamma(k) \ dk)= tr(F(T)\circ t)
\ .\]
In order to decide injectivity of $F^*$, we consider for $(\sigma,\lambda)\in
\hat \M\times \aaa^*$, $(\lambda,\alpha)>0$ $\forall \alpha\in\Phi^+(\g,\aaa)$
the $\oc$-function
\[\oc_\gamma(\sigma,\lambda):= P_\sigma \int_{\overline \Naaa} a(\bar
n)^{-(\lambda +\rho)}\gamma (\kappa(\bar n))\ d\bar n P_\sigma \ \in End_\M
V_\gamma \ ,\]
where $P_\sigma$ is the projection onto the $\sigma$-isotypic component of
$V_\gamma$, $\overline \Naaa := \theta \Naaa$ and $\bar n = \kappa(\bar n)
a(\bar n) n(\bar n)$ is decomposed according to the Iwasawa decomposition.
There exist constants $C^\prime_\alpha(\sigma,\lambda)$ such that
\begin{equation}\label{f.ccc}
det \oc_\gamma(\sigma,\lambda)\not=0 \qquad \mbox{whenever}\  \frac
{(\lambda,\alpha)}{(\alpha,\alpha)}\geq C^\prime_\alpha (\sigma,\gamma)\quad
\forall \alpha \in \Phi^+(\g,\aaa)\ .
\end{equation}
This can be seen using the expression of $det \oc_\gamma(\sigma,\lambda)$ in
terms of the $\Gamma$-function given for example in Wallach \cite{wallach92},
Cor.10.4.7. Alternatively one could argue that the integrand concentrates near
the identity for $\lambda$ large.

Let $\langle .,.\rangle$ be a $\G^d$-invariant hermitian scalar product on
$V_\pi$, $v\in V_\gamma$,$v_\lambda \in V_\pi^\naaa$, $t\in
Hom_\M(V_\gamma,V_\sigma)$. We compute
\begin{eqnarray*}
\langle F^*(t)v,v_\lambda\rangle &=& \langle\int_\K \pi(k^{-1}) i\circ t
\gamma(k)v\ dk,v_\lambda\rangle \\
&=& \langle \int_{\overline \Naaa \M} a(\bar n)^{-2\rho} \pi (\kappa(\bar
n)^{-1}) \pi(m)^{-1} i\circ t \gamma(m) \gamma(\kappa(\bar n))v \ d\bar n,
v_\lambda \rangle \\
&=& \int_{\overline \Naaa} \langle a(\bar n)^{-2\rho} i\circ
t\gamma(\kappa(\bar n))v , \pi (\kappa(\bar n))v_\lambda \rangle\ d\bar n \\
&=& \int_{\overline \Naaa} \langle i\circ t a(\bar n)^{-(\lambda +2\rho)}
\gamma(\kappa(\bar n))v , \pi (\bar n)v_\lambda \rangle\ d\bar n \\
&=& \langle i\circ t \oc_\gamma(\sigma,\lambda+\rho)v,v_\lambda \rangle \ .
\end{eqnarray*}
In the last step we have used that with respect to  $\langle.,.\rangle$
\[ \pi(\bar n)^* = \pi (\theta (\bar n)^{-1}) \ .\]
We see that $F^*$ is injective iff $\oc_\gamma(\sigma,\lambda+\rho)$ is
injective. This together with (\ref{f.ccc}) implies the assertion. $\Box$
\newline

 If $\gamma$ is trivial, then the \oc-function is explicitly known. In this
case we have for all $\alpha\in \Phi^+(\g,\aaa)$ $C_\alpha
=\Phi_\alpha=\varepsilon_\alpha=0$. Inserting these data into the previous
theorem we get the usual Cartan-Helgason theorem (see Helgason
\cite{helgason84},Ch.V,Thm.4.1.).\newline

Now we want to see more explicitly what Theorem \ref{th.CH} means for the
spectral decomposition of $L^2(\G^d\times_\K V_\gamma)$ with respect of the
action of the Casimir operator, if $\G^d/\K$ has rank one.

In our approach $X=\G/\K$ comes with a fixed Riemannian metric. It induces a
metric on $\aaa$ (and $\aaa^*$) considering $\aaa$ as a subspace of $T_{eK}X$.
The root system $\Phi^+(\g,\aaa)$ is of the form $\Phi^+(\g,\aaa)=\{\alpha\}$
or $\Phi^+(\g,\aaa)=\{\alpha/2, \alpha\}$. We call $\alpha$ the long root. We
identify
\[ \aaa^*\cong \R,\quad (\aaa^c)^*\cong \C \]
sending $\alpha$ to $|\alpha|$. The metric on $\aaa$ provides a normalization
of the $Ad(\g^c)$-invariant bilinear form $(.,.)$ on $(\g^c)^*\times (\g^c)^*$
(and $ \g^c \times  \g^c $) by
\[ (\alpha,\alpha)=|\alpha|^2\ .\]
Note that $(.,.)_{|\g^d\times\g^d}$ is negative definite. Let $\Omega_d$ be the
corresponding Casimir element in the universal enveloping algebra ${\cal
U}(\g^d)$.

Now we fix the necessary bookkeeping for its action on $L^2(\G^d\times_\K
V_\gamma)$.
\begin{ddd}\label{d.constants}
Let $\alpha$ be the long root as before and set $
T:= |\alpha|$.

For $\sigma\in \hat \M$ we define $\epsilon(\sigma)\in \{0,1/2\}$ by
\[ \epsilon(\sigma) \equiv \rho/T + \varepsilon_{\alpha}(\sigma)\ \ mod\:\Z \
,\]
the lattice $L(\sigma)\subset \R\cong \aaa^*$ by
$$ L(\sigma):= T(\epsilon(\sigma)+\Z) $$
and the shift constant
\[ c(\sigma):= |\delta|^2 - |\mu_\sigma + \rho_m |^2 = |\rho |^2 + |\rho_m |^2
- |\mu_\sigma + \rho_m |^2 \ .\]
Finally, for $\lambda \in (\aaa^c)^*\cong \C$ we set
\[ P(\lambda,\sigma) := \prod_{\beta \in \Phi^+(\g^c,\haaa)} \frac {(\lambda +
\mu_\sigma +\rho_m,\beta)}{(\delta,\beta)} \ .\]
\end{ddd}
Recall the definition of $\Phi_\alpha(\sigma)$ from Theorem \ref{th.CH} and
let
\[ \Phi(\sigma):= \max\{
|\alpha|\Phi_{\alpha}(\sigma)\:|\:\alpha\in\Phi^+(\g,\aaa)\}+\rho \ .\]
Consider for $\lambda\in L(\sigma), \lambda\geq\Phi(\sigma)$ and $\Lambda:=
\lambda-\rho +\mu_\sigma$ the representation $\pi_\Lambda$ of $\G^d$ with
highest weight $\Lambda$. Then Weyl's dimension formula states that
\[ dim (\pi_\Lambda) = P(\lambda,\sigma)\ \ .\]
In addition, we have
\[ \pi_\Lambda (\Omega_d) = (|\Lambda +\delta |^2 -|\delta |^2 )Id = (\lambda^2
-c(\sigma))Id \ .\]
Now, in our situation Theorem \ref{th.CH} reads as
\begin{prop}\label{p.CH}
The spectrum of $\Omega_d$ acting on $L^2(\G^d\times_\K V_\gamma)$ is contained
in
\[ \bigcup_{\{\sigma\in\hat \M|\:[\gamma:\sigma]\not=0\}} \{ \lambda^2
-c(\sigma)\:|\: \lambda \in L(\sigma), \lambda\geq \Phi(\sigma)\}\ .\]
Moreover, any eigenspace $E_\nu$ of $\Omega_d$ splits into a direct sum of
subspaces
\[ E_\nu = \bigoplus_{(\sigma,\lambda)\in Q_\gamma(\nu)} E_{\sigma,\lambda}\
,\]
where
\[ Q_\gamma(\nu):=\{(\sigma,\lambda)\:|\  [\gamma:\sigma]\not=0, \lambda\in
L(\sigma),\lambda\geq\Phi(\sigma),\lambda^2-c(\sigma)=\nu \}\ ,\]
and there exist constants $C(\sigma,\gamma)\geq\Phi(\sigma)$ such that
\[ dim E_{\sigma,\lambda} = [\gamma:\sigma] P(\lambda,\sigma)\quad \forall
\lambda\in L(\sigma),\lambda\geq C(\sigma,\gamma)\ . \Box\]
\end{prop}

\subsection{The $\G^d$-index theorem}

Let us return to the situation at the beginning of Subsection \ref{chchc}. We
assume in addition that $rank(\G^d)=rank(\K)$. Otherwise the index of any
$\G^d$-homogeneous elliptic operator on $X_d$ would vanish (see Bott
\cite{bott65}). Note that the rank condition is satisfied by symmetric spaces
of rank one of even dimension. Our goal is to realize any (irreducible)
representation of $\G^d$ as the index of a homogeneous Dirac operator operator
on $X_d$. That this is possible, even for general compact homgeneous spaces
$\G/\Haaa$ satisfying the equal rank condition, was already shown by Bott
(\cite{bott65}) under the additional hypothesis that $\pi_1(\G)$ has no
2-torsion. Instead of refering to that theorem, we prefer to include here an
elementary proof resting only on the Weyl character formula and having the
advantage that it gives explicitly the desired Dirac operator.

First, we need some more root systems and Weyl groups. We choose a Cartan
algebra $\haaa_k\subset \kaaa$
being also a Cartan algebra of $\g^d$. Let $\Haaa$ be the corresponding maximal
torus. Let $\Phi^+(\g^c,\haaa_k)$ be a system of positive
roots for $\g^d$ and $\Phi^+(\kaaa^c,\haaa_k)$ be the subsystem of positive
roots of $\kaaa$. Note that
$\sharp(\Phi^+(\g^c,\haaa_k)\setminus\Phi^+(\kaaa^c,\haaa_k))=n/2$, where $n$
is the dimension of $X_d$.
We set $$\rho_g:=\frac{1}{2}\sum_{\alpha\in\Phi^+(\g^c,\haaa_k)}
\alpha,\quad\rho_k:=\frac{1}{2}\sum_{\alpha\in\Phi^+(\kaaa^c,\haaa_k)}\alpha\
.$$
Let $W_g:=W(\Phi^+(\g^c,\haaa_k))$, $W_k:=W(\Phi^+(\kaaa^c,\haaa_k))$
be the corresponding Weyl groups and set
$$W^1=\{w\in W_g\:|\:\Phi^+(\kaaa^c,\haaa_k)\subset
w(\Phi^+(\g^c,\haaa_k))\}.$$

The following considerations and computations are completely analogous to the
ones in Parthasarathy \cite{parthasarathy72} for the noncompact case.
Let $S:=S^+\oplus S^-=S(\paaa^d)$ be the Spinor modul of $\paaa^d$
and $c:\paaa^d\rightarrow End(S)$ be the Clifford multiplication.
The Spinor modul $S$ is a $Spin(\paaa^d)$-representation. Let
$s:Spin(\paaa^d)\rightarrow GL(S)$ and
$p:Spin(\paaa^d)\rightarrow SO(\paaa^d)$ be the canonical maps.
Then
\begin{equation}\label{mmm1} s(a)c(X)s(a^{-1})=c(p(a)X)\quad X\in \paaa^d,\quad
a\in Spin(n).\end{equation}
The isotropy representation gives a homomorphism $i:\K\rightarrow
SO(\paaa^d)$.
Let $\tilde{\K}$ be the connected component of $1$ in $p^{-1}(i(\K))\subset
Spin(\paaa^d)$. We will also denote by $s$ the restriction of $s$ to
$\tilde\K$.

Let $\Lambda\in\imath\haaa_k^*$ be the highest weight of an irreducible
representation $\pi_\Lambda$ of $\G^d$. It is not difficult to see (compare
\cite{parthasarathy72}) that for all $w\in W^1$ the functional
$w(\Lambda+\rho_g)-\rho_k$ is dominant and integral with respect to $\tilde\K$,
hence it is the highest weight of an irreducible representation
$\gamma_{\Lambda,w}$ of $\tilde\K$ acting on a vector space $V_{\Lambda,w}$.
Furthermore, if the non-trivial element of $p^{-1}(i(1))$ is included in
$\tilde\K$, then it
acts on $S$ as well as on $V_{\Lambda,w}$ by multiplication with $-1$. Indeed,
this element is central in $Spin(\paaa^d)$, hence in $\tilde\K$, and of the
form $exp(H)$ for some $H\in\haaa_k$. Therefore, we have
\[ \gamma_{\Lambda,w}(exp(H))=e^{(w(\Lambda+\rho_g)-\rho_k)(H)}\cdot Id =
e^{(w\rho_g-\rho_k)(H)}\cdot Id = -Id\ ,\]
since $w\rho_g-\rho_k$ is a weight of $S$.

Thus, the representation
$$s\otimes \gamma_{\Lambda,w}:\K\rightarrow GL(S\otimes V_{\Lambda,w})$$
is well defined as a representation of $\K$.
We extend the Clifford multiplication to $S\otimes V_{\Lambda,w}$ by
$$\tilde{c}:\paaa^d\rightarrow End(S\otimes V_{\Lambda,w}),\quad
\tilde{c}:=c\otimes 1.$$
We consider the homogeneous vector bundle
$E_{\Lambda,w}:=\G^d\times_\K(S\otimes V_{\Lambda,w})$ and
define the differential operator $D_{\Lambda,w}$, the Dirac operator associated
to the
Dirac bundle structure of $E_{\Lambda,w}$, by
\begin{equation}\label{f.dirac}
D_{\Lambda,w}f(g):=\sum_{i=1}^n\frac{\partial}{\partial t}{}_{|t=0}
\tilde{c}(X_i) f(g\:exp(tX_i))\ ,
\end{equation}
where $\{X_i\}_{i=1}^n$ is an orthonormal basis of $\paaa^d$, and we have
identified
\begin{equation}\label{f.rolf}
\Gamma(E_{\Lambda,w})=\{f\in C^\infty(\G^d,S\otimes
V_{\Lambda,w})\:|\:f(gk)=(s\otimes \gamma_{\Lambda,w})(k^{-1})f(g)\}\ .
\end{equation}
The grading $S=S^+\oplus S^-$ induces a corresponding grading
$E_{\Lambda,w}=E_{\Lambda,w}^+\oplus E_{\Lambda,w}^-$. The operator
$D_{\Lambda,w}$ splits as
$D_{\Lambda,w}^\pm:\Gamma(E_{\Lambda,w}^\pm)\rightarrow
\Gamma(E_{\Lambda,w}^\mp)$.
By (\ref{mmm1}) we get
\begin{lem}\label{lllll1}
$D^\pm$ are $\G^d$-invariant differential operators.
\end{lem}
The $\G^d$-index of $D_{\Lambda,w}$ is the element in the representation ring
$R(\G^d)$ of $\G^d$ given by
\[ ind_{\G^d} D_{\Lambda,w}:= [ker D^+_{\Lambda,w}]-[ker D^-_{\Lambda,w}]\ .\]
\begin{theorem}\label{th.ind}
The $\G^d$-index of the homogeneous Dirac operator $D_{\Lambda,w}$ defined in
(\ref{f.dirac}) equals
\[ (-1)^{n/2}det(w)[\pi_\Lambda]\in R(\G^d)\ .\]
\end{theorem}
{\it Proof }: Let $\psi\in C^\infty(\G^d)$ a $\G^d$-finite function. $\psi$
acts on $\Gamma(E_{\Lambda,w})$ as a finite dimensional operator $l_\psi$ via
convolution. For $f\in \Gamma(E_{\Lambda,w})$ in the identification
(\ref{f.rolf}) this action is given by
\begin{eqnarray}
l_\psi (f)(g)&=& \int_{\G^d} \psi(x) f(x^{-1}g)\ dx\nonumber\\
&=&\int_{\G^d} \psi(gx^{-1}) f(x)\ dx\nonumber\\
&=&\int_{\G^d} \psi(gkx^{-1}) (s\otimes\gamma_{\Lambda,w})(k)f(x)\
dx\quad\mbox{(for any $k\in\K$) }\nonumber\\
&=&\int_{\G^d}\int_\K \psi(gkx^{-1}) (s\otimes\gamma_{\Lambda,w})(k)\ dk\ f(x)\
dx\ \label{f.faltung}.
\end{eqnarray}
We assume all Haar measures to be normalized having total mass 1.

It is enough to show that the supertrace of $l_\psi$ with respect to the
grading of $E_{\Lambda,w}$ is given by integration of $\psi$ against the
character $\chi_\Lambda$ of $\pi_\Lambda$ with the appropriate sign. In fact,
it is enough to show it for class functions $\psi$. In view of
(\ref{f.faltung}) $l_\psi$ is an integral operator with kernel given by
\[ l_{\psi}(x,y)=\int_\K \psi(xky^{-1}) (s\otimes\gamma_{\Lambda,w})(k)\ dk\
.\]
Thus, for class functions $\psi$ its supertrace can be obtained in the
following way
\begin{eqnarray*}
 Tr (l_{\psi})&=& \int_{X_d} \int_\K \psi(xkx^{-1}) tr
(s\otimes\gamma_{\Lambda,w})(k)\ dk\ dx \\
&=&\int_\K \psi(k) tr (s\otimes\gamma_{\Lambda,w})(k)\ dk\ .
\end{eqnarray*}
Let us denote by $\chi_{S^\pm}$ and $\chi_{\gamma}$ the characters of $S^\pm$
and $\gamma_{\Lambda,w}$, respectively. They are uniquely defined as functions
on $\tilde\K$, but their product drops down to $\K$. We introduce the Weyl
denominators
\begin{eqnarray*}
\Delta_{\G^d}(expX)&:=& \prod_{\alpha\in \Phi^+(\g^c,\haaa_k)}(e^{\alpha(X)/2}
- e^{-\alpha(X)/2})\ ,\\
\Delta_{\K}(expX)&:=& \prod_{\alpha\in \Phi^+(\kaaa^c,\haaa_k)}(e^{\alpha(X)/2}
- e^{-\alpha(X)/2})\ ,\quad X\in\haaa_k,\\
\Delta_n&:=& \Delta_{\G^d}/\Delta_{\K}\ .
\end{eqnarray*}
Now, using Weyl's integral and character formulae for $\K$ and $\G^d$ and
\[ (\chi_{S^+} - \chi_{S^-})_{|\Haaa} = \Delta_n\ \quad \mbox{(see
\cite{parthasarathy72})}\ ,\]
we can compute
\begin{eqnarray*}
Tr(l_\psi)&=& \int_\K \psi(k) tr (s\otimes\gamma_{\Lambda,w})(k)\ dk\\
&=& \int_\K \psi(k)(\chi_{S^+}-\chi_{S^-})(k)\chi_\gamma(k)\ dk\\
&=& \frac{1}{|W_k|}\int_\Haaa \psi(h)|\Delta_\K(h)|^2
(\chi_{S^+}-\chi_{S^-})(h)\chi_\gamma(h)\ dh\\
&=&\frac{1}{|W_k|}\int_\Haaa \psi(h)|\Delta_\K(h)|^2 \Delta_n(h) \sum_{s\in
W_k} \frac{det(s)e^{sw(\Lambda+\rho_g)(log(h))}}{\Delta_\K(h)}\ dh\\
&=&\frac{1}{|W_k|}\int_\Haaa \psi(h)|\Delta_{\G^d}(h)|^2 \sum_{s\in W_k}
\frac{det(s)e^{sw(\Lambda+\rho_g)(log(h))}}
{\Delta_\K(h)\overline{\Delta_n(h)}}\ dh\\
&=&\frac{(-1)^{n/2}}{|W_k||W_g|}\int_\Haaa \psi(h)
|\Delta_{\G^d}(h)|^2
\sum_{t\in W_g}\sum_{s\in W_k}
\frac{det(s)e^{tsw(\Lambda+\rho_g)(log(h))}}
{\Delta_{\G^d}(t^{-1}h)}\ dh\\
&=&\frac{(-1)^{n/2}det(w)}{|W_g|}\int_\Haaa \psi(h)
|\Delta_{\G^d}(h)|^2
\sum_{t\in W_g}  \frac{det(t)e^{t(\Lambda+\rho_g)(log(h))}}
{\Delta_{\G^d}(h)}\
dh\\
&=&\frac{(-1)^{n/2}det(w)}{|W_g|}\int_\Haaa \psi(h)
|\Delta_{\G^d}(h)|^2
\chi_\Lambda(h)\ dh \\
&=& (-1)^{n/2} det(w) \int_{\G^d}
\psi(g) \chi_\Lambda(g)\ dg \ .\Box
\end{eqnarray*}

\subsection{The restriction map $r:R(\K)\rightarrow R(\M)$}

The goal of the present subsection is to discuss the surjectivity of the
restriction map
\[ r: R(\K)\longrightarrow R(\M) \]
between the representation rings over $\Z$ of $\K$ and $\M$, respectively. We
retain all notations and assumptions from the beginning of Subsection
\ref{chchc}. We shall restrict our considerations to the case of rank one
symmetric spaces of even dimension. Recall the definitions given after
(\ref{mmm1}) of the covering group $\tilde\K$ and its spinor representation
$s$. For the present case we have
\begin{prop}\label{p.sur}
The map $r$ is surjective.

Moreover, if $\gamma\in ker\:r\subset R(\K)$, then there exists a
representation $\tau$ of $\tilde\K$, such that $s\otimes\tau$ is defined on
$\K$ and
\[ \gamma = [ s^+\otimes\tau ]- [s^-\otimes\tau ] \in R(\K)\ .\]
In particular, this implies that on the $\Z_2$-graded vector bundle
$\G^d\times_\K V_\gamma$ there exists an odd homogeneous Dirac operator.
\end{prop}
{\it Proof }: The assertion concerning the kernel of $r$ is proved in Miatello
\cite{miatello83}, Lemma 2.2.
We prove the surjectivity case by case, showing that the fundamental
representations of $\M$ lie in the image of $r$.
\newline

{\bf The real case}\newline
The two-dimensional case is clear. We may represent the real hyperbolic space
$X=\HR^{2n},n\geq2$ either as $X=Spin(1,2n)/Spin(2n)$ or as
$X=SO_0(1,2n)/SO(2n)$. Let us begin with the first case. Then
\[ \K=Spin(2n)\quad\mbox{and}\quad\M=Spin(2n-1)\ .\]
Let $\lambda_p:=\Lambda^p\R^{2n} \otimes\C$ be the exterior powers of the
standard representation of
$SO(2n)$ and $\sigma_p:=\Lambda^p\R^{2n-1} \otimes\C$ the ones of $SO(2n-1)$,
$s^+$ and $s^-$ the half spin representations of $Spin(2n)$ as before, and
$s_0$ the spin representation of $Spin(2n-1)$. Then
\[ R(\K)=\Z [\lambda_1,\dots,\lambda_{n-2},s^+,s^-]\quad\mbox{and}\quad
R(\M)=\Z [\sigma_1,\dots,\sigma_{n-2},s_0]\ ,\]
and we have
\[r(\lambda_1)=\sigma_1+1\ ;\quad r(\lambda_p)=\sigma_p+\sigma_{p-1},\
p=2,\dots,n-2\ ; \quad r(s^+)=r(s^-)=s_0 \ .\]
Preimages are now easily to be found, namely
\begin{eqnarray}\label{f.sose}
 \sigma_p &= &r(
\lambda_p-\lambda_{p-1}+\dots...+(-1)^{p-1}(\lambda_1-1))\:,\quad
p=1,\dots,n-2\;\\ s_0&=&r(s^+)=r(s^-) \ .\nonumber
\end{eqnarray}
For the second case, i.e.
\[ \K=SO(2n)\quad\mbox{and}\quad\M=SO(2n-1)\ ,\]
we remark that the representation ring of $\M$  is given by
\[R(\M)=\Z [\sigma_1,\dots,\sigma_{n-2},\sigma_{n-1}]\ ,\]
therefore all its generators arise as in (\ref{f.sose}) as restrictions of
elements of $R(SO(2n))$.
\newline

{\bf The complex case}\newline
The isometry group of the complex hyperbolic space $X=\HC^{n}, n\geq2$ is
$PU(1,n)\cong SU(1,n)/\Z_{n+1}$, whereas its maximal linear covering group is
$SU(1,n)$. Thus, for any divisor $k$ of $n+1$ we may write $X=\G/\K$, where
$$\G=SU(1,n)/\Z_k\ ,\quad \K=S(U(1)\times U(n))/\Z_k\ .$$
The corresponding group $\M$ is
\[ \M=\{(z,B)\in U(1)\times U(n-1)\:|\: z^2det(B)=1\}/\Z_k\ , \]
where the embedding $\M\hookrightarrow \K$ is given by
\[ [(z,B)]\longrightarrow [(z,\left(
\begin{array}{cc}
z&0\\ 0&B
\end{array}
\right) )]\ .\]
We consider $\K$ as a $\frac{n+1}{k}$-fold covering of $U(n)$ given by
\[ \pi: \K\ni [(w,A)]\rightarrow w^{-1}A\in U(n)\ .\]
Note that $\pi(\M)=U(n-1)$.

Let $\lambda_p:=\Lambda^p\C^{n}$ be the complex exterior powers of the standard
representation of
$U(n)$ and $\sigma_p:=\Lambda^p\C^{n-1}$ the ones of $U(n-1)$. By means of
$\pi$, we consider them as representations of $\K$ and $\M$, respectively.
Then
\[ R(\K)=\Z [\lambda_1,\dots,\lambda_{n-1},w^k,w^{-k}]\quad\mbox{and}\quad
R(\M)=\Z [\sigma_1,\dots,\sigma_{n-2},z^k,z^{-k}]\ ,\]
where $w$ denotes the one dimensional representation $(w,A)\rightarrow w$, and
$z$ is given analogously by $(z,B)\rightarrow z$.
The restriction map looks as follows
\begin{eqnarray*}
r(w^{\pm k})&=&z^{\pm k}\ ,\\
r(\lambda_1)&=&\sigma_1+1\ ,\\
r(\lambda_p)&=&\sigma_p+\sigma_{p-1},\quad p=2,\dots,n-2\ , \\
r(\lambda_{n-1})&=& z^{-(n+1)} +\sigma_{n-2}\ .
\end{eqnarray*}
Hence,
\begin{eqnarray*}
z^{\pm k}&=&r(w^{\pm k})\ ,\\
 \sigma_p &= &r(
\lambda_p-\lambda_{p-1}+\dots...+(-1)^{p-1}(\lambda_1-1))\:,\quad
p=1,\dots,n-2.
\end{eqnarray*}

{\bf The quaternionic case}\newline
We may represent the quaternionic hyperbolic space $X=\HH^{n},n\geq 2$ either
as $X=Sp(1,n)/Sp(1)\times Sp(n)$ or as $X=(Sp(1,n)/\Z_2)/Sp(1)\cdot Sp(n)$. Let
us begin with the first case. Then
\[ \K=Sp(1)\times Sp(n)\quad\mbox{and}\quad\M\cong Sp(1)\times Sp(n-1)\ ,\]
where the embedding $Sp(1)\times Sp(n-1)\hookrightarrow Sp(1)\times Sp(n)$ is
given by
\[ Sp(1)\times Sp(n-1)\ni (q,B)\longrightarrow (q, \left(
\begin{array}{cc}
q&0\\ 0&B
\end{array}
\right) )\in Sp(1)\times Sp(n)\ .\]
Let $\lambda_p:=\Lambda^p\C^{2n}$ be the complex exterior powers of the
standard representation of
$Sp(n)$ on $\Haaa^n\cong \C^{2n}$ and $\sigma_p:=\Lambda^p\C^{n-1}$ the ones of
$Sp(1)\times Sp(n-1)$. Note that for $p>1$ they are not irreducible,
nevertheless they generate the representation ring of $Sp(n)$. Extend them
trivially to $Sp(1)\times Sp(n)$ and $Sp(1)\times Sp(n-1)$, respectively. Then
\[ R(\K)=\Z [p,\lambda_1,\dots,\lambda_n]\quad\mbox{and}\quad R(\M)=\Z
[q,\sigma_1,\dots,\sigma_{n-1}]\ ,\]
where $p$ (resp. $q$) denotes the standard representation of $Sp(1)$ trivially
extended to $Sp(1)\times Sp(n)$ (resp. $Sp(1)\times Sp(n-1)$).
The restriction map looks as follows
\begin{eqnarray*}
r(p)&=& q\ ,\\
r(\lambda_1)&=&\sigma_1+q\ ,\\
r(\lambda_p)&=&\sigma_p+q\sigma_{p-1},\quad p=2,\dots,n-1\ , \\
r(\lambda_{n})&=& q\sigma_{n-1}+\sigma_{n-2}\ .
\end{eqnarray*}
Preimages of the generators can now inductively be determined by
\begin{eqnarray}\label{f.holle}
q&=&r(p)\ ,\nonumber\\
\sigma_1&=& r(\lambda_1 -p)\ ,\nonumber\\
&\vdots&\\
\sigma_p &= &r( \lambda_p)-q\sigma_{p-1}\ .\nonumber
\end{eqnarray}
In the second case, i.e.
\[ \K^\prime =Sp(1)\cdot Sp(n)\cong Sp(1)\times_{\Z_2} Sp(n)
\quad\mbox{and}\quad \M^\prime \cong Sp(1)\times_{\Z_2} Sp(n-1)\ ,\]
we see that the representation ring of $\K^\prime$  is given as the subring of
$R(\K)$ generated by
\[
\{p^2,\lambda_i,p\lambda_j,\lambda_j\lambda_k\:
|\:i\:\mbox{even},\:j,k\:\mbox{odd},\:1\leq i,j,k\leq n\}\ .\]
Analogously, the representation ring of $\M^\prime$  is given as the subring of
$R(\M)$ generated by
\begin{equation}\label{f.pech}
\{q^2,\sigma_i,q\sigma_j,\sigma_i\sigma_j\:
|\:i\:\mbox{even},\:j,k\:\mbox{odd},\:1\leq i,j,k\leq n-1\}\ .
\end{equation}
Call an element in $\hat{\K}$ ($\hat{\M}$) even, if $(-1,-Id)$ acts trivially,
and odd otherwise. Then we have
\begin{eqnarray*}
\hat\K^\prime &=& \{\gamma\in \hat{\K}\:|\:\gamma\: \mbox{even}\}\ ,\\
\hat\M^\prime &=& \{\sigma\in \hat{\M}\:|\:\sigma\: \mbox{even}\}\ .
\end{eqnarray*}
Thus
\begin{eqnarray*}
R(\K^\prime)&=& \{ \sum_{\{\gamma\in \hat{\K}\:|\:\gamma\: \mbox{even}\}}
a_\gamma \gamma\:|\:a_\gamma\in \Z\}\ ,\\
R(\M^\prime)&=&\sum_{\{\sigma\in \hat{\M}\:|\:\sigma\: \mbox{even}\}} a_\sigma
\sigma\:|\:a_\sigma\in \Z\}\ .
\end{eqnarray*}
{}From (\ref{f.holle}) we see that any even (odd) generator of $R(\M)$ can be
expressed as the restriction of a weighted sum of even (odd) elements in
$\hat{\K}$. We conclude that all generators in (\ref{f.pech}) arise as
restrictions of a weighted sum of even elements in $\hat{\K}$, hence belong to
 $r(R(\K^\prime))$.

{\bf The exceptional case}\newline
The unique presentation of the hyperbolic Cayley plane $X=H{\bf Ca}^{2}$ is
$$X=F_4^{-20}/Spin(9)\ . $$
Then
\[ \K=Spin(9)\quad\mbox{and}\quad\M=Spin(7)\ ,\]
where $Spin(7)$ is the stabilizer of a non-zero element in the spinor modul
$S_9$ of $Spin(9)$ (which is isomorphic to the isotropy representation on
$T_{e\K}H{\bf Ca}^{2}$).
Let $\lambda_p:=\Lambda^p\R^{9} \otimes\C$ be the exterior powers of the
standard representation of
$SO(9)$ and $\sigma_p:=\Lambda^p\R^{7} \otimes\C$ the ones of $SO(7)$, $s_9$
and $s_7$ the spin representations of $Spin(9)$ and $Spin(7)$, respectively.
Then
\[ R(\K)=\Z [\lambda_1,\lambda_2,\lambda_3,s_9]\quad\mbox{and}\quad R(\M)=\Z
[\sigma_1,\sigma_2,s_7]\ ,\]
and one can show that
\[r(\lambda_1)=s_7+1\ ,\quad r(\lambda_2)=\sigma_2+\sigma_1+s_7\ , \quad
r(s_9)=s_7+\sigma_1+1 \ .\]
Preimages are now easily to be found, namely
\begin{eqnarray*}
\sigma_1 &= &r(s_9 - \lambda_1)\ ,\\
\sigma_2&=&r(\lambda_2 -s_9 +1)\ , \\
s_7 &=& r(\lambda_1 -1) \ .\quad\Box
\end{eqnarray*}

\bibliographystyle{plain}

\end{document}